\documentclass[a4paper,11pt,reqno]{amsart}

\usepackage[T1]{fontenc}
\usepackage{lmodern,microtype}
\usepackage[a4paper,textwidth=160mm,textheight=235mm,centering,headheight=14pt,headsep=7mm]{geometry}
\usepackage{amsmath,amssymb,mathtools,esint}
\usepackage{xcolor}
\usepackage[
  unicode,
  colorlinks=true,
  linkcolor=blue,
  citecolor=blue,
  urlcolor=black
]{hyperref}
\usepackage[noabbrev]{cleveref}

\numberwithin{equation}{section}
\allowdisplaybreaks[1]
\theoremstyle{plain}
\newtheorem{theorem}{Theorem}[section]
\newtheorem{lemma}[theorem]{Lemma}

\newtheorem{proposition}[theorem]{Proposition}
\newtheorem{corollary}[theorem]{Corollary}
\crefname{theorem}{Theorem}{Theorems}
\crefname{lemma}{Lemma}{Lemmas}
\crefname{proposition}{Proposition}{Propositions}
\crefname{corollary}{Corollary}{Corollaries}
\crefname{section}{Section}{Sections}
\crefname{equation}{equation}{equations}
\newcommand{\dashint}{\fint}
\AtBeginDocument{%
  \setlength{\abovedisplayskip}{6pt plus 2pt minus 2pt}%
  \setlength{\belowdisplayskip}{6pt plus 2pt minus 2pt}%
  \setlength{\abovedisplayshortskip}{3pt plus 2pt}%
  \setlength{\belowdisplayshortskip}{4pt plus 2pt minus 1pt}%
}

\author{Huyuan Chen}
\address{Huyuan Chen:\ \, Center for Mathematics and Interdisciplinary Sciences, Fudan University,
Shanghai 200433, China\newline
\indent Shanghai Institute for Mathematics and Interdisciplinary Sciences, Shanghai 200433, China.}
\email{chenhuyuan@yeah.net, chenhuyuan@simis.cn}

\author{Rui Chen}
\address{Rui Chen:\ \, School of Mathematical Sciences, Fudan University, Shanghai 200433, China\newline
\indent Brandenburg University of Technology Cottbus--Senftenberg, Platz der Deutschen Einheit 1, 03046 Cottbus, Germany}
\email{chenrui23@m.fudan.edu.cn}

\author{Daniel Hauer}
\address{Daniel Hauer:\ \, Brandenburg University of Technology Cottbus--Senftenberg, Platz der Deutschen Einheit 1, 03046 Cottbus, Germany\newline
\indent School of Mathematics and Statistics, The University of Sydney, Sydney, NSW 2006, Australia.}
\email{daniel.hauer@b-tu.de}



\hypersetup{
  pdftitle={Sharp boundary regularity for the fractional-logarithmic Laplacian},
  pdfauthor={Huyuan Chen, Rui Chen, Daniel Hauer},
  pdfsubject={Sharp boundary regularity for the fractional-logarithmic Laplacian},
  pdfkeywords={Fractional-logarithmic Laplacian, boundary regularity,
    boundary asymptotics, optimal modulus of continuity,
    nonlocal Dirichlet problem}
}

\title[Sharp boundary regularity for the fractional-logarithmic Laplacian]
{Sharp boundary regularity for the fractional-logarithmic Laplacian}
\begin{document}
\raggedbottom

\begin{abstract}
We study the boundary regularity of variational solutions to
the Dirichlet problem for the fractional-logarithmic Laplacian,
with Fourier symbol $|\xi|^{2s}\ln|\xi|^2$, on bounded
$C^{1,\alpha}$ open sets $\Omega\subset\mathbb R^n$, where $0<s<1$
and $0<\alpha\leq1$.
For bounded right-hand sides, we construct an explicit boundary
factor $\Phi_s(t)\sim t^s/\sqrt{\ln(1/t)}$ via a Mellin integral
and prove that the quotient $u/\Phi_s(\delta)$ is H\"older
continuous up to the boundary, where $\delta$ denotes the
distance to the boundary. Its optimal modulus of continuity
is $r^s/\sqrt{\ln(e/r)}$ when $n=1$ or $\alpha>s$, and
$r^\alpha$ when $n\geq2$ and $\alpha\leq s$.
We establish the sharpness of both the quotient modulus
and the boundary remainder estimates, on every fixed open
set in the former case and in the class of $C^{1,\alpha}$
open sets in the latter.
Despite the asymptotic equivalence of the two boundary factors,
on every such open set we construct a solution with bounded
right-hand side for which $u\sqrt{\ln(R/\delta)}/\delta^s$
is continuous up to the boundary but fails to be H\"older
continuous of any positive order near any boundary point,
for every $R>e\operatorname{diam}(\Omega)$.
The corrected boundary trace can be chosen identically one
when $n=1$ or $\alpha>s$, whereas for $n\geq2$ and
$\alpha\leq s$ we construct domains on which a nonzero
constant trace is impossible even locally near a boundary point.
\end{abstract}

\maketitle
\section{Introduction and main results}

For $0<s<1$, we study the optimal boundary regularity
of solutions to the Dirichlet problem
\begin{equation}\label{eq:dirichlet}
    (-\Delta)^{s+\ln}u=f\quad\text{in }\Omega,\qquad
    u=0\quad\text{in }\Omega^c,
\end{equation}
where $\Omega\subset\mathbb R^n$, $n\geq1$, is a nonempty
bounded open set with $C^{1,\alpha}$ boundary,
$0<\alpha\leq1$, $\Omega^c=\mathbb R^n\setminus\Omega$,
and $f\in L^\infty(\Omega)$.

The fractional-logarithmic Laplacian $(-\Delta)^{s+\ln}$
was introduced in~\cite{ChenChenHauer2026} as the derivative
of the fractional Laplacian with respect to its order:
\[
    (-\Delta)^{s+\ln}u
    :=\left.\frac{d}{dt}(-\Delta)^tu\right|_{t=s},
    \qquad u\in\mathcal S(\mathbb R^n).
\]
Under the unitary Fourier transform
$\widehat u(\xi)=(2\pi)^{-n/2}\int_{\mathbb R^n}
e^{-ix\cdot\xi}u(x)\,dx$, its Fourier symbol is
$|\xi|^{2s}\ln|\xi|^2$. The operator admits the singular
integral representation
\[
    (-\Delta)^{s+\ln}u(x)
    =c_{n,s}\operatorname{p.v.}\int_{\mathbb R^n}
    \frac{u(x)-u(y)}{|x-y|^{n+2s}}
    \bigl(b_{n,s}-2\ln|x-y|\bigr)\,dy,
\]
where $c_{n,s}=4^s\pi^{-n/2}s\Gamma(n/2+s)/\Gamma(1-s)$,
$b_{n,s}=\partial_s\ln c_{n,s}$, and $\Gamma$ denotes
the gamma function; see~\cite[Theorem~1.1(ii)]
{ChenChenHauer2026}.
The fractional-logarithmic Laplacian retains several
properties of the fractional Laplacian and converges to the logarithmic
Laplacian as $s\downarrow0$, but lacks homogeneity under
dilations and has a sign-changing kernel.
Its extension representation, energy spaces, Poisson problem,
and Dirichlet spectrum were studied in~\cite{ChenChenHauer2026},
and its potential theory and global regularity
in~\cite{Chen2026Potentials}.
Biswas, Das, and Sen~\cite{BiswasDasSen2026} recently
established interior Harnack inequalities and H\"older
estimates for the fractional-logarithmic $p$-Laplacian.

Throughout the paper, we write
$\delta=\operatorname{dist}(\cdot,\partial\Omega)$ and
$\delta_0=\operatorname{dist}(\cdot,\Omega^c)$.
The Dirichlet energy space $\mathcal H_0^{s+\ln}(\Omega)$
and the variational formulation of~\eqref{eq:dirichlet}
are given in Section~\ref{sec:forms}.

\subsection{Main results}

For the fractional Laplacian on bounded $C^{1,1}$ domains,
Ros-Oton and Serra~\cite{RosOtonSerra2014} proved that
solutions with bounded right-hand sides and zero exterior
values satisfy $u\in C^s(\mathbb R^n)$ and that
$u/\delta^s$ is H\"older continuous up to the boundary.
The boundary decay of order $\delta^s$ is optimal.
Differentiating $(-\Delta)^s$ at $s=0$ yields the
logarithmic Laplacian, for which Hern\'andez-Santamar\'ia,
L\'opez R\'ios, and Salda\~na~
\cite{HernandezSantamariaLopezRiosSaldana2025}
established the optimal boundary decay of order
$(\ln(1/\delta))^{-1/2}$ under a uniform exterior sphere
condition. Since $(-\Delta)^{s+\ln}$ is obtained by
differentiating $(-\Delta)^s$ at a positive order $s$,
these results suggest the boundary scale
$\delta^s/\sqrt{\ln(1/\delta)}$.

The same prediction follows formally from the boundary
regularity theory for subordinate Brownian motions.
Under suitable assumptions on the Bernstein function
$\varphi$, Kim, Kim, Lee, and Lee~\cite{KimKimLeeLee2019}
proved the H\"older regularity of $u/V(\delta)$ on bounded
$C^{1,1}$ open sets, where the renewal function $V$ satisfies
$V(t)\asymp\varphi(t^{-2})^{-1/2}$ as $t\downarrow0$;
see~\cite{KimSongVondracek2012,KimMimica2014}.
Here $a\asymp b$ means $C^{-1}b\leq a\leq Cb$
for some $C\geq1$.
Formally taking $\varphi(\lambda)=\lambda^s\ln\lambda$ gives
\[
    \varphi(t^{-2})^{-1/2}
    =\frac{t^s}{\sqrt{2\ln(1/t)}}.
\]
However, this theory does not apply to $(-\Delta)^{s+\ln}$,
since $\lambda^s\ln\lambda$ is not a Bernstein function
and the kernel changes sign. These observations suggest
that the optimal boundary decay of solutions
to~\eqref{eq:dirichlet} is of order
$\delta^s/\sqrt{\ln(1/\delta)}$, which we establish on
bounded $C^{1,\alpha}$ open sets.

To state the boundary decay estimate, fix
$R>e\operatorname{diam}(\Omega)$ and set
$L_R(t)=\ln(R/t)$ for $0<t<R$, where
$\operatorname{diam}(\Omega)$ denotes the diameter of $\Omega$.
Define
\[
    \rho_{s,R}(x)=\frac{\delta(x)^s}{\sqrt{L_R(\delta(x))}},
    \qquad x\in\Omega.
\]
Although $\rho_{s,R}$ describes the optimal boundary decay,
the quotient $u/\rho_{s,R}$ need not be H\"older continuous
up to the boundary.

We therefore construct a more precise boundary factor
$\Phi_s$ and establish the optimal modulus of continuity
of $u/\Phi_s(\delta)$. The construction is based on the identity
\[
    \frac{t^s}{\sqrt{\ln(1/t)}}
    =\frac{1}{\sqrt{\pi}}\int_0^\infty a^{-1/2}t^{s+a}\,da,
    \qquad 0<t<1.
\]
This identity motivates the construction of $\Phi_s$
as an integral of the powers $t_+^{s+a}$.
The density is chosen to cancel the boundary singularities
of $(-\Delta_{\mathbb R})^{s+\ln}\Phi_s$ by integration
by parts in $a$.
Choose a nonnegative cutoff function
$\chi\in C_c^\infty((-2r_\chi,2r_\chi))$
with $\chi=1$ on $(-r_\chi,r_\chi)$, where
$r_\chi\in(0,1/4)$ is fixed throughout the paper.
For $a\in(0,1)$, define
\[
    m_s(a)=\frac{\Gamma(s+1)}
    {\Gamma(s+1+a)\sqrt{\sin(\pi a)}}.
\]
Since $m_s(a)\asymp a^{-1/2}(1-a)^{-1/2}$, we have
$m_s\in L^1(0,1)$. Define
\begin{equation}
    \Phi_s(t)=
    \begin{cases}
        \displaystyle\chi(t)t^s\int_0^1m_s(a)t^a\,da,
        &t>0,\\
        0,&t\leq0.
    \end{cases}
    \label{eq:boundary-factor}
\end{equation}
Section~\ref{sec:mellin} proves that
$(-\Delta_{\mathbb R})^{s+\ln}\Phi_s$ extends smoothly
to $t=0$ from the right and that
$\Phi_s(t)\sim t^s/\sqrt{\ln(1/t)}$ as $t\downarrow0$.
Note that $\Phi_s$ is independent of $R$.

The optimal modulus of continuity of $u/\Phi_s(\delta)$
near the boundary depends on the dimension and the
regularity of $\partial\Omega$.
For $0<r\leq1$, set $\mu_s(r)=r^s/\sqrt{\ln(e/r)}$ and define
\begin{equation}
    \omega_{n,s,\alpha}(r)=
    \begin{cases}
        \mu_s(r),&n=1\ \text{or }\alpha>s,\\
        r^\alpha,&n\geq2\ \text{and }0<\alpha\leq s.
    \end{cases}
    \label{eq:domain-modulus}
\end{equation}
We set $\mu_s(0)=\omega_{n,s,\alpha}(0)=0$.

Throughout the paper, we write
$\mathcal A=\|u\|_{L^2(\Omega)}+\|f\|_{L^\infty(\Omega)}$.
We denote by $\nu_{\rm in}$ the inward unit normal
to $\partial\Omega$ and, for each $x\in\Omega$, let
$z_x\in\partial\Omega$ be any point satisfying
$|x-z_x|=\delta(x)$.
For $r>0$ and $z\in\mathbb R^n$, write
$B_r(z)=\{x\in\mathbb R^n:|x-z|<r\}$ and $B_r=B_r(0)$.
For a measurable set $E$ of positive finite measure
and $v\in L^1(E)$, set $\dashint_E v=|E|^{-1}\int_E v$.
When $n\geq2$, we write
$x=(x',x_n)\in\mathbb R^{n-1}\times\mathbb R$.

Our first theorem establishes the H\"older continuity of
$u/\Phi_s(\delta)$ with the modulus $\omega_{n,s,\alpha}$
and gives a quantitative boundary expansion.

\begin{theorem}
\label{thm:quotient}
Let $\Omega\subset\mathbb R^n$ be a bounded $C^{1,\alpha}$
open set, where $0<s<1$ and $0<\alpha\leq1$.
Let $u\in\mathcal H_0^{s+\ln}(\Omega)$ be a variational
solution of~\eqref{eq:dirichlet} with $f\in L^\infty(\Omega)$.
There exist $\varepsilon_0\in(0,r_\chi)$ and $C>0$,
depending only on $n,s,\alpha,\Omega$ and the fixed cutoff,
such that the following assertions hold.

\begin{enumerate}
\item[\textup{(i)}] \textbf{Boundary quotient.}
The quotient $q=u/\Phi_s(\delta)$ extends continuously
to $\overline\Omega\cap\{\delta\leq\varepsilon_0\}$ and satisfies
\begin{equation}
    \|q\|_{L^\infty(\overline\Omega\cap
    \{\delta\leq\varepsilon_0\})}
    \leq C\mathcal A.
    \label{eq:quotient-bound}
\end{equation}
For all $x,y\in\overline\Omega\cap\{\delta\leq\varepsilon_0\}$
with $|x-y|\leq\varepsilon_0$,
\begin{equation}
    |q(x)-q(y)|
    \leq C\mathcal A\omega_{n,s,\alpha}(|x-y|).
    \label{eq:quotient-modulus}
\end{equation}

\item[\textup{(ii)}] \textbf{Boundary expansion.}
Let $c=q|_{\partial\Omega}$.
For every $z\in\partial\Omega$ and $0<r<\varepsilon_0$,
\begin{equation}
    \|u-c(z)\Phi_s(\delta_0)\|_{L^\infty(B_r(z))}
    \leq C\mathcal A\Phi_s(r)\omega_{n,s,\alpha}(r),
    \label{eq:boundary-supremum}
\end{equation}
where $u$ and $\delta_0$ are understood to vanish
outside $\Omega$.
\end{enumerate}
\end{theorem}

Theorem~\ref{thm:quotient} yields $C^{\min\{s,\alpha\}}$
regularity of the boundary quotient when $n\geq2$, including
the critical case $\alpha=s$, and $C^s$ regularity when $n=1$.
If $n=1$ or $\alpha>s$, the stronger modulus
$\mu_s(r)=o(r^s)$ improves the $C^s$ estimate by the
logarithmic factor $(\ln(e/r))^{-1/2}$.
The optimality of these moduli is established in
Theorem~\ref{thm:sharpness}.

Combining Theorem~\ref{thm:quotient} with the global
boundedness estimate obtained in Section~\ref{sec:forms}
and the asymptotic formula for $\Phi_s$, we obtain
the following corollary.

\begin{corollary}
\label{cor:boundedness}
Under the hypotheses of Theorem~\ref{thm:quotient},
there exist $C_0,C>0$ such that
$\|u\|_{L^\infty(\Omega)}\leq C_0\mathcal A$ and,
for a.e. $x\in\Omega$,
\begin{equation}\label{eq:boundary-decay}
 |u(x)|\leq C\mathcal A\rho_{s,R}(x).
\end{equation}
Here $C_0$ depends only on $n$ and $s$, while $C$
may also depend on $\alpha$, $\Omega$, $R$ and the fixed cutoff.
\end{corollary}

The preceding estimates involve both $\|u\|_{L^2(\Omega)}$
and $\|f\|_{L^\infty(\Omega)}$.
The $L^2$ term cannot be omitted if zero is a Dirichlet
eigenvalue, since there exist nonzero solutions with $f=0$.
If zero lies in the resolvent set, all estimates can be
expressed solely in terms of $\|f\|_{L^\infty(\Omega)}$.
We make this precise in the following corollary.

\begin{corollary}
\label{cor:resolvent}
Let $(-\Delta)^{s+\ln}_{\Omega}$ be the self-adjoint
Dirichlet realization defined in Section~\ref{sec:forms}.
Suppose that $0\in\rho((-\Delta)^{s+\ln}_{\Omega})$ and set $d_{s,\Omega}
    =\operatorname{dist}
    (0,\sigma((-\Delta)^{s+\ln}_{\Omega}))>0,$
where $\rho$ and $\sigma$ denote the resolvent set and
the spectrum, respectively.
Then, for every $f\in L^\infty(\Omega)$, problem
\eqref{eq:dirichlet} admits a unique variational solution.
Moreover, all estimates in Theorem~\ref{thm:quotient}
and Corollary~\ref{cor:boundedness} hold with $\mathcal A$
replaced by
\[
    \left(1+\frac{|\Omega|^{1/2}}{d_{s,\Omega}}\right)
    \|f\|_{L^\infty(\Omega)}.
\]
\end{corollary}

Although $\Phi_s(\delta)$ and $\rho_{s,R}$ are asymptotically
equivalent, replacing the former by the latter can destroy
the H\"older regularity of the boundary quotient. We show
this on every fixed open set by constructing a solution
with strictly positive boundary trace. We also determine
when the trace can be prescribed to be a nonzero constant
on every open set in the given regularity class.

\begin{theorem}
\label{thm:bare-quotient}
The following assertions hold.

\begin{enumerate}
\item[\textup{(i)}] \textbf{Positive boundary trace.}
For every $\Omega,s,\alpha$ as above, there exists a
variational solution $U\in\mathcal H_0^{s+\ln}(\Omega)$
of~\eqref{eq:dirichlet} with $f\in L^\infty(\Omega)$
whose boundary trace
$c_U=(U/\Phi_s(\delta))|_{\partial\Omega}$ satisfies
$\min_{\partial\Omega}c_U>0$.
For every $R>e\operatorname{diam}(\Omega)$,
the quotient $U/\rho_{s,R}$ extends continuously to
$\partial\Omega$ with trace $c_U$, but admits no
$C^\gamma$ extension near any boundary point for
any $\gamma>0$.

\item[\textup{(ii)}] \textbf{Constant boundary trace.}
If $n=1$ or $\alpha>s$, the solution $U$ in
part~\textup{(i)} can be chosen so that $c_U\equiv1$.

\item[\textup{(iii)}] \textbf{Nonexistence of constant traces.}
If $n\geq2$ and $0<\alpha\leq s<1$, there exist a bounded
connected $C^{1,\alpha}$ domain $\Omega_*\subset\mathbb R^n$
and a point $z_*\in\partial\Omega_*$ such that no variational
solution of~\eqref{eq:dirichlet} on $\Omega_*$ with bounded
right-hand side has a nonzero constant trace of
$u/\Phi_s(\delta)$ on any relative neighborhood of $z_*$
in $\partial\Omega_*$.
\end{enumerate}
\end{theorem}

Part~\textup{(i)} shows that the boundary decay in
Corollary~\ref{cor:boundedness} is optimal on every fixed
domain and at every boundary point.
Although $u/\rho_{s,R}$ and $u/\Phi_s(\delta)$ have the same
continuous boundary trace, the former need not be
H\"older continuous of any positive order, whereas the
latter satisfies the modulus estimate in
Theorem~\ref{thm:quotient}.
This follows from the expansion in
Proposition~\ref{prop:factor-asymptotics}:
\[
    \frac{\sqrt{L_R(t)}\Phi_s(t)}{t^s}
    =1+\frac{A_{1,s,R}}{L_R(t)}
    +\frac{A_{2,s,R}}{L_R(t)^2}
    +O(L_R(t)^{-3})
\]
as $t\downarrow0$, where
$A_{1,s,R}=(\ln R-\psi(s+1))/2$ and
\[
    A_{2,s,R}
    =\frac38\left((\ln R-\psi(s+1))^2
    +\frac{\pi^2}{6}-\psi_1(s+1)\right)>0.
\]
Here $\psi=\Gamma'/\Gamma$ and $\psi_1=\psi'$ denote
the digamma and trigamma functions, respectively.
Thus no choice of $R$ eliminates both inverse-logarithmic
corrections. At every boundary point where the trace is
nonzero, the first nonzero correction prevents
$u/\rho_{s,R}$ from being H\"older continuous of any
positive order.

Parts~\textup{(ii)} and~\textup{(iii)} show that, for
$n\geq2$, the condition $\alpha>s$ is sharp for guaranteeing
the existence of solutions with a nonzero constant
boundary trace on every bounded $C^{1,\alpha}$ domain.
More precisely, when $\alpha\leq s$, the domain whose
boundary is locally given by $x_n=-|x'|^{1+\alpha}$
admits no solution with bounded right-hand side and
nonzero constant trace even locally near the origin.
This does not exclude constant boundary traces on
individual $C^{1,\alpha}$ domains; for instance,
part~\textup{(ii)} applies whenever the domain is $C^{1,1}$.

\smallskip

We next prove that the quotient modulus and the boundary
expansion in Theorem~\ref{thm:quotient} are optimal.
The lower bound $\mu_s(r)$ holds on every fixed open set,
whereas the lower bound $r^\alpha$ for $n\geq2$ and
$\alpha\leq s$ is attained on the domain used in
part~\textup{(iii)} of Theorem~\ref{thm:bare-quotient}.

\begin{theorem}
\label{thm:sharpness}
The following assertions hold.

\begin{enumerate}
\item[\textup{(i)}] \textbf{Fixed-domain examples.}
Let $\Omega,s,\alpha$ be as in Theorem~\ref{thm:quotient}
and fix $z\in\partial\Omega$.
There exist $r_*,C_*,c_*>0$ such that, for every
$0<r<r_*$, there are a solution $u_r\in C_c^\infty(\Omega)$
of~\eqref{eq:dirichlet} with right-hand side $f_r$ satisfying
$\mathcal A_r:=\|u_r\|_2+\|f_r\|_\infty\leq C_*$
and a point $x_r\in\Omega$ with
$|x_r-z|=r$ and $\delta(x_r)\asymp r$ such that
\[
 |q_r(x_r)-q_r(z)|\geq c_*\mu_s(r),
\]
where $q_r=u_r/\Phi_s(\delta)$ has zero boundary trace.

Moreover, there exists a variational solution on $\Omega$
with bounded right-hand side whose quotient
$u/\Phi_s(\delta)$ is not $C^{s+\varepsilon}$ near $z$
for any $\varepsilon>0$.

\item[\textup{(ii)}] \textbf{Rough-domain examples.}
Let $n\geq2$ and $0<\alpha\leq s<1$.
There exists a bounded connected $C^{1,\alpha}$ domain
$\Omega$ with $0\in\partial\Omega$ such that every
variational solution $U$ with bounded right-hand side
and nonzero boundary trace $c_0:=q(0)$ satisfies
\[
    \liminf_{r\downarrow0}
    \frac{\dashint_{B_r}|U-c_0\Phi_s(\delta_0)|\,dx}
    {\Phi_s(r)r^\alpha}>0,
\]
where $q=U/\Phi_s(\delta)$.
Moreover, such solutions exist with $c_0>0$.
\end{enumerate}
\end{theorem}

Part~\textup{(i)} excludes every modulus $o(\mu_s)$
in estimates uniform over solutions on any fixed $\Omega$
as in Theorem~\ref{thm:quotient}.
Since $q_r(z)=0$ and $\delta(x_r)\asymp r$, we also have
\[
    \|u_r\|_{L^\infty(B_{2r}(z))}
    \geq |u_r(x_r)|
    =\Phi_s(\delta(x_r))|q_r(x_r)|
    \geq c\Phi_s(r)\mu_s(r).
\]
Thus both the quotient modulus and the boundary remainder
estimate in Theorem~\ref{thm:quotient} are optimal on every
fixed $\Omega$ when $n=1$ or $\alpha>s$.

For $n\geq2$ and $0<\alpha\leq s$, the lower bound in
part~\textup{(ii)} and the inequality
\[
    \dashint_{B_r}|U-c_0\Phi_s(\delta_0)|\,dx
    \leq\Phi_s(r)\sup_{\Omega\cap B_r}|q-c_0|
\]
give $\sup_{\Omega\cap B_r}|q-c_0|\geq c_*r^\alpha$
for sufficiently small $r>0$.
Thus the modulus $r^\alpha$ is optimal in the class of
bounded $C^{1,\alpha}$ domains.
In particular, at $\alpha=s$, the modulus $r^s$ cannot
be improved to $\mu_s(r)$.
By H\"older's inequality, the $L^1$ lower bound also
yields the corresponding lower bounds for the local
$L^2$ average and the supremum of the boundary remainder.

Unlike part~\textup{(i)}, the lower bound in
part~\textup{(ii)} cannot hold on every bounded
$C^{1,\alpha}$ domain.
Indeed, on any bounded $C^{1,1}$ domain,
Theorem~\ref{thm:quotient} gives
\[
    \|U-c(z)\Phi_s(\delta_0)\|_{L^\infty(B_r(z))}
    \leq C\mathcal A\Phi_s(r)\mu_s(r)
    =o(\Phi_s(r)r^\alpha)
\]
for every solution with bounded right-hand side and
every $z\in\partial\Omega$, since $\alpha\leq s$.

\subsection{Outline of the proofs}

\noindent\textbf{Proof of Theorem~\ref{thm:quotient}.}
Unlike the fractional Laplacian, the fractional-logarithmic
Laplacian does not annihilate the half-line power $t_+^s$.
We construct $\Phi_s$ by differentiating the fractional
Laplacian of half-line powers~\cite[(2.8)]{FallRosOton2022}
and choosing a Mellin density that converts the resulting
expression into a total derivative in the Mellin parameter.
Subtracting the operator value at a fixed positive point
removes the pole at $a=s$. Integration by parts then gives
the smoothness of $(-\Delta_{\mathbb R})^{s+\ln}\Phi_s$
up to zero from the right
(Proposition~\ref{prop:mellin-factor}).

The lack of homogeneity and the sign change of the kernel
require local estimates uniform in the rescaled operators
\[
    A_{\ell_r}=(-\Delta)^s+
    \frac{1}{2\ell_r}(-\Delta)^{s+\ln},
    \qquad \ell_r=\ln(1/r).
\]
For the regularized boundary factor
$\Psi_s=\zeta\Phi_s(D)$, a further difficulty arises
when $n\geq2$ and $\alpha=s$:
its operator image $G_s=(-\Delta)^{s+\ln}\Psi_s$
may grow as $\ln(e/\delta)^{3/2}$.
Lemma~\ref{lem:critical-image} removes the leading
contribution by subtracting a constant $b_{z,r}$:
\[
    r^{-n/2}\|G_s-b_{z,r}\|_{L^2(\Omega\cap B_r(z))}
    \leq C\sqrt{\ln(e/r)}.
\]

We adapt the least-squares boundary approximation
of~\cite{RosOtonSerra2016Stable,RosOtonSerra2017Boundary,
Fall2019NonlocalSchrodinger} to establish
\[
    \inf_{a\in\mathbb R}
    \left(\dashint_{B_r(z)}|u-a\Psi_s|^2\,dx\right)^{1/2}
    \leq C\mathcal A\Phi_s(r)\omega_{n,s,\alpha}(r).
\]
In contrast to the subcritical boundary estimates
in~\cite{RosOtonSerra2016Stable,RosOtonSerra2017Boundary},
the rescaled errors may now grow at order $2s$,
at which the fractional tail need not be integrable.
We compensate the tail by subtracting a constant
and use the half-space classification
of~\cite[Theorem~3.10 and Proposition~3.11]
{AbatangeloRosOton2020}.
The uniform interior H\"older estimates then yield
the pointwise boundary expansion and the exact
modulus of continuity of the quotient.

\medskip
\noindent\textbf{Proof of Theorem~\ref{thm:bare-quotient}.}
To construct a solution with positive boundary
trace, we adapt the higher-power barrier
of~\cite[Lemma~2.10]{RosOtonSerra2017Boundary}
by adding $\zeta(D_+)^{s+\beta}$ to $\Psi_s$, where
$0<\beta<\min\{\alpha,s\}$.
The term preserves the boundary trace and its
negative operator image dominates the
singularity of $G_s$. To account for the signed kernel,
we perform the comparison using its positive part.

The first two inverse-logarithmic corrections of $\Phi_s$
exclude every positive H\"older exponent for
$U/\rho_{s,R}$ at points with nonzero boundary trace.
When $n=1$ or $\alpha>s$, the boundedness of $G_s$
allows us to take $U=\Psi_s$ with boundary trace one.
For $n\geq2$ and $\alpha\leq s$, we show that even a
locally constant nonzero trace may be impossible.
On the domain with boundary $x_n=-|x'|^{1+\alpha}$
near zero, such a trace would give a limiting function
\[
    W(y',t)=s|y'|^{1+\alpha}t^{s-1}+E(y',t),
    \qquad |E(y',t)|\leq Ct^{s+\alpha}.
\]
Taking tangential second differences and subtracting
the corresponding Martin integral
\cite[Theorem~1.3 and~(1.12)]{AbatangeloEtAl2019HalfSpace},
we reduce the problem to the half-space Liouville
classification~\cite[Theorem~3.10]{AbatangeloRosOton2020}.
The resulting equation for the tangential coefficient
involves $(-\Delta_{\mathbb R^{n-1}})^{1/2}$
and contradicts its decay at infinity and positive maximum.

\medskip
\noindent\textbf{Proof of Theorem~\ref{thm:sharpness}.}
We distinguish between the restrictions imposed by
the operator and the boundary.
On every fixed domain, interior bumps of radius $r$
and amplitude $r^{2s}/\ln(1/r)$ have uniformly bounded
right-hand sides and attain the quotient modulus
$\mu_s(r)$. A superposition at sufficiently separated
scales gives a single solution without higher
H\"older regularity.

For $n\geq2$ and $\alpha\leq s$, the boundary geometry
requires a different argument.
On the domain used in
Theorem~\ref{thm:bare-quotient}\textup{(iii)},
a boundary remainder of order
$o(\Phi_s(r)r^\alpha)$ along a sequence of scales
would force the explicit geometric correction
\[
    |y'|^{1+\alpha}(y_n)_+^{s-1}
\]
to be fractional harmonic in the upper half-space
modulo constants.
For $\alpha<s$, a direct computation using the
strict convexity of $|y'|^{1+\alpha}$ excludes this
possibility.
At $\alpha=s$, the correction is homogeneous of
degree $2s$, and its nonzero spherical mean contradicts
Lemma~\ref{lem:critical-homogeneous}.
This proves the sharpness of $r^\alpha$, including
the critical modulus $r^s$ rather than merely the
optimality of the H\"older exponent.

\subsection{Related results on boundary regularity}

We conclude the introduction by reviewing related results
on the boundary regularity of other nonlocal operators.

For the fractional Laplacian, Bogdan~\cite{Bogdan1997}
established the boundary Harnack principle on Lipschitz domains,
and Chen and Song~\cite{ChenSong1998} obtained two-sided Green
function and Poisson kernel estimates on bounded $C^{1,1}$
domains in dimension $n\geq2$.
Further boundary Harnack estimates for stable processes
and more general jump processes were developed
in~\cite{SongWu1999,BogdanKulczyckiKwasnicki2008,
BogdanKumagaiKwasnicki2015}.
Interior Harnack inequalities and H\"older estimates for
nonlocal equations were established
in~\cite{BassLevin2002,Silvestre2006,CaffarelliSilvestre2009}.

Following the boundary quotient estimates for the fractional
Laplacian in~\cite{RosOtonSerra2014}, Ros-Oton and Serra
\cite[Theorem~1.2]{RosOtonSerra2016Stable} proved
$u/\delta^s\in C^{s-\varepsilon}(\overline\Omega)$
for uniformly elliptic symmetric stable operators on
bounded $C^{1,1}$ domains,
where $0<\varepsilon<s$.
For symmetric homogeneous kernels comparable to $|y|^{-n-2s}$,
they subsequently established $C^\alpha$ regularity of
$u/\delta^s$ on $C^{1,\alpha}$ domains when $0<\alpha<s$,
as well as a boundary Harnack principle on $C^1$ domains
\cite[Theorems~1.2 and~1.3]{RosOtonSerra2017Boundary}.
For nonsymmetric stable operators, Dipierro, Ros-Oton,
Serra, and Valdinoci~\cite{DipierroRosOtonSerraValdinoci2022}
obtained H\"older estimates for $u/\delta^\gamma$, where
$\gamma$ depends on the symbol of the operator and the
inward unit normal.
Related half-space classifications, boundary Harnack
inequalities, and gradient estimates can be found
in~\cite{RosOtonSerra2016FullyNonlinear,
DeSilvaSavin2016,RosOtonSerra2019Harnack,FallJarohs2021}.

Using the $\mu$-transmission calculus, Grubb
\cite[Theorem~2.2]{Grubb2014} established
$u/\delta^s\in C^s(\overline\Omega)$ for the fractional
Dirichlet problem with bounded right-hand side on smooth
domains when $s\ne1/2$, and $C^{s-\varepsilon}$ regularity
for every $\varepsilon\in(0,s)$ when $s=1/2$.
For extensions to more general fractional-order
pseudodifferential operators and nonsmooth domains,
see~\cite{Grubb2015,AbelsGrubb2023,Grubb2026}.

For nonhomogeneous kernels, Kassmann and Mimica
\cite{KassmannMimica2017} developed interior regularity
estimates adapted to the scaling of the kernel.
Ros-Oton and Weidner
\cite[Theorems~1.1 and~1.6]{RosOtonWeidner2024}
established optimal $C^s$ boundary regularity for
translation-invariant operators with symmetric kernels
comparable to $|y|^{-n-2s}$
and obtained boundary expansions in terms of
one-dimensional half-space solutions on $C^{1,\alpha}$
domains.
Kim and Weidner~\cite{KimWeidner2025} proved $C^s$
boundary regularity and sharp Green function estimates
for nonlocal equations in divergence form with
H\"older continuous coefficients.
More recently, Grube and Ros-Oton~\cite{GrubeRosOton2026}
established optimal $C^s$ boundary regularity for
symmetric L\'evy operators with Fourier symbols
comparable to $|\xi|^{2s}$ under a $C^1$-Dini-type
boundary condition.
Guo, Li, and Ouyang~\cite{GuoLiOuyang2026} obtained
local boundary H\"older estimates for nonnegative
fractional Dirichlet solutions on unbounded $C^{1,1}$
domains without assuming a global $L^\infty$ bound.

For the fractional $p$-Laplacian, interior Harnack
inequalities and H\"older estimates were developed
in~\cite{DiCastroKuusiPalatucci2014,
DiCastroKuusiPalatucci2016,Cozzi2017}.
Global regularity and boundary quotient estimates
were subsequently established
in~\cite{IannizzottoMosconiSquassina2016,
IannizzottoMosconiSquassina2020,IannizzottoMosconi2024,
IannizzottoMosconi2026}.

For logarithmic-type operators, Grzywny and Ryznar
\cite[Lemma~4.1]{GrzywnyRyznar2012} obtained the
renewal function estimate
$V(t)\asymp(\ln(1/t))^{-1/2}$ as $t\downarrow0$
for one-dimensional geometric stable processes
with generator $-\ln(1+(-\Delta)^\theta)$,
where $0<\theta\leq1$.
Power-logarithmic complete Bernstein functions
were also studied in~\cite[Section~2.1]
{KimSongVondracekRevisited2012}.

For the logarithmic Laplacian, Chen and Weth
\cite[Section~3 and Theorem~1.11]{ChenWeth2019}
developed the Dirichlet theory and proved boundary
decay of order $(\ln(1/\delta))^{-\tau}$ for
$0<\tau<1/2$ on bounded Lipschitz domains satisfying
a uniform exterior sphere condition.
Hern\'andez-Santamar\'ia, L\'opez R\'ios, and Salda\~na
\cite[Theorems~1.1, 1.2 and~1.4]
{HernandezSantamariaLopezRiosSaldana2025}
established the optimal decay of order
$(\ln(1/\delta))^{-1/2}$, proved its sharpness on balls,
and obtained a Hopf-type estimate under an interior
sphere condition.
The dependence of fractional Poisson solutions and
Dirichlet eigenvalues on the fractional order was
studied in~\cite{JarohsSaldanaWeth2020,
FeulefackJarohsWeth2022}.
Further regularity estimates for zero-order operators
and the logarithmic Laplacian were obtained
in~\cite{ChangLaraSaldana2024,BiswasRoySenSharma2026}.

\medskip

The remainder of the paper is organized as follows.
Section~\ref{sec:forms} establishes uniform local energy,
boundedness, and interior H\"older estimates for the rescaled
operators. Section~\ref{sec:mellin} constructs the boundary
factor $\Phi_s$ and estimates its operator image on
$C^{1,\alpha}$ domains. Section~\ref{sec:boundary-quotient}
proves Theorem~\ref{thm:quotient} and its corollaries.
Section~\ref{sec:sharpness} proves
Theorems~\ref{thm:bare-quotient} and~\ref{thm:sharpness}.

\section{Quadratic forms and local estimates}
\label{sec:forms}

We recall the quadratic form of $(-\Delta)^{s+\ln}$
and its self-adjoint Dirichlet realization.
We then derive the scaling identities and establish
local compactness, uniform local $L^\infty$ estimates,
and interior H\"older estimates for the rescaled operators
$A_\ell$. These results are used in the boundary blow-up
arguments of Section~\ref{sec:boundary-quotient}.

Throughout the paper, all functions are real-valued unless
otherwise stated. Functions satisfying the zero-exterior condition
are identified with their zero extensions to $\mathbb R^n$.
We denote by $e_n=(0,\ldots,0,1)$ the last coordinate unit vector
and by $O(n)$ the orthogonal group. We write $Q^T$ for the transpose
of a matrix $Q$, $B_R(z)$ for the radius-$R$ ball centered at $z$,
and $B_R=B_R(0)$.
For $0<\beta\leq1$, we use the H\"older seminorm
\[
    [v]_{C^{0,\beta}(E)}
    =\sup_{x,y\in E,\,x\ne y}
    \frac{|v(x)-v(y)|}{|x-y|^\beta}
\]
and the norm $\|v\|_{C^{0,\beta}(E)}=\sup_E|v|+[v]_{C^{0,\beta}(E)}$.
We write $C^\beta(E)=C^{0,\beta}(E)$ for $0<\beta<1$
and use $C^{0,1}(E)$ for Lipschitz regularity.
Unless otherwise specified, all integrals and $L^p$ norms
are taken over $\mathbb R^n$. Norms of right-hand sides
are taken over the set where the equation is imposed.
The constant $C>0$ may vary from line to line, and its dependence
is specified in the corresponding results.
The estimates are not required to be uniform as
$s\downarrow0$ or $s\uparrow1$.

\subsection{The quadratic form and its Dirichlet realization}

The kernel and the whole-space and zero-exterior energy spaces
come from~\cite[Theorem~1.1(ii) and Definitions~1.3--1.4]{ChenChenHauer2026}.
We record the form estimates used below and derive the weighted
estimate needed for the boundary factor.

\smallskip
\noindent\textbf{(i) The kernel and energy spaces.}
We use the definition of $(-\Delta)^{s+\ln}$ and the
coefficients $c_{n,s}$ and $b_{n,s}$ from the introduction.
Its integral kernel is
$K_{s+\ln}(r)=c_{n,s}(b_{n,s}-2\ln r)r^{-n-2s}$ for $r>0$.
Principal values are taken by excluding $B_\varepsilon(x)$
and letting $\varepsilon\downarrow0$.

We use the energy space of~\cite[Definition~1.3]{ChenChenHauer2026},
with the equivalent norm
\begin{equation}
 \|v\|_{\mathcal H^{s+\ln}}^2
 :=\|v\|_{L^2(\mathbb R^n)}^2
 +c_{n,s}\iint_{|x-y|<1}
 \frac{|v(x)-v(y)|^2\ln(1/|x-y|)}
 {|x-y|^{n+2s}}\,dx\,dy.
 \label{eq:energy-norm}
\end{equation}
Thus $\mathcal H^{s+\ln}(\mathbb R^n)=\{v\in L^2(\mathbb R^n):
\|v\|_{\mathcal H^{s+\ln}}<\infty\}$.
For an open set $\mathcal O\subset\mathbb R^n$, the zero-exterior
space of~\cite[Definition~1.4]{ChenChenHauer2026} is
\[
 \mathcal H_0^{s+\ln}(\mathcal O)
 :=\bigl\{v\in\mathcal H^{s+\ln}(\mathbb R^n):
          v=0\text{ a.e. in }\mathcal O^c\bigr\}.
\]
We abbreviate $\|v\|_{L^2(\mathbb R^n)}$ by $\|v\|_2$.
For $0<s<1$, we write
$H^s(\mathbb R^n)=\{v\in L^2(\mathbb R^n):[v]_{H^s}<\infty\}$
and $\|v\|_{H^s}^2=\|v\|_2^2+[v]_{H^s}^2$, where
\[
 [v]_{H^s}^2
 :=\iint_{\mathbb R^{2n}}
 \frac{|v(x)-v(y)|^2}{|x-y|^{n+2s}}\,dx\,dy.
\]
The notation $H^s_{\mathrm{loc}}(\mathbb R^n)$ means that
$\eta v\in H^s(\mathbb R^n)$ for every $\eta\in C_c^\infty(\mathbb R^n)$.

\smallskip
\noindent\textbf{(ii) The quadratic form and variational solutions.}
For $u,v\in\mathcal H^{s+\ln}(\mathbb R^n)$, define
\[
 \mathcal E_{s+\ln}(u,v)
 :=\frac12\iint_{\mathbb R^{2n}}
 (u(x)-u(y))(v(x)-v(y))
 K_{s+\ln}(|x-y|)\,dx\,dy.
\]

Choose $0<r_*<\exp(-1-|b_{n,s}|)$. For $0<r<r_*$,
\[
 c_{n,s}\frac{\ln(1/r)}{r^{n+2s}}
 \leq K_{s+\ln}(r)
 \leq 3c_{n,s}\frac{\ln(1/r)}{r^{n+2s}},
\]
whereas
$K_{s+\ln}(|\cdot|)\mathbf1_{\{|\cdot|\geq r_*\}}
\in L^1(\mathbb R^n)$.
Cauchy--Schwarz and translation invariance give
\[
 \frac12\iint_{|x-y|\geq r_*}
 |u(x)-u(y)|\,|v(x)-v(y)|\,|K_{s+\ln}(|x-y|)|\,dx\,dy
 \leq 2\|u\|_2\|v\|_2\int_{|h|\geq r_*}|K_{s+\ln}(|h|)|\,dh.
\]
Together with the short-distance kernel estimates, this proves
that $\mathcal E_{s+\ln}$ is absolutely convergent and continuous
on $\mathcal H^{s+\ln}(\mathbb R^n)\times
\mathcal H^{s+\ln}(\mathbb R^n)$.
Moreover, since the contribution of $r_*\leq|x-y|<1$
to~\eqref{eq:energy-norm} is bounded by $C\|v\|_2^2$,
there exist $\lambda>0$ and $C\geq1$, depending only on $n,s$,
such that
\begin{equation}
 C^{-1}\|v\|_{\mathcal H^{s+\ln}}^2
 \leq\mathcal E_{s+\ln}(v,v)+\lambda\|v\|_2^2
 \leq C\|v\|_{\mathcal H^{s+\ln}}^2.
 \label{eq:shifted-form}
\end{equation}
The same decomposition at $|x-y|=r_*$ gives
$\|v\|_{H^s}\leq C\|v\|_{\mathcal H^{s+\ln}}$ and hence the continuous
embedding $\mathcal H^{s+\ln}(\mathbb R^n)\hookrightarrow H^s(\mathbb R^n)$.

A function $u\in\mathcal H_0^{s+\ln}(\Omega)$ is a variational
solution of~\eqref{eq:dirichlet} if
$\mathcal E_{s+\ln}(u,\varphi)=\int_\Omega f\varphi\,dx$
for every $\varphi\in\mathcal H_0^{s+\ln}(\Omega)$.
Two elementary properties of the energy space will be used
in the local estimates. First, if $F:\mathbb R\to\mathbb R$
is Lipschitz and $F(0)=0$, then
$\|F(v)\|_{\mathcal H^{s+\ln}}
\leq\operatorname{Lip}(F)\|v\|_{\mathcal H^{s+\ln}}$
for $v\in\mathcal H^{s+\ln}(\mathbb R^n)$,
where $\operatorname{Lip}(F)$ is its Lipschitz constant.
In particular, Lipschitz truncations preserve
$\mathcal H_0^{s+\ln}(\mathcal O)$.
Second, multiplication by $\eta\in C_c^\infty(\mathbb R^n)$
is bounded on $\mathcal H^{s+\ln}(\mathbb R^n)$.

\smallskip
\noindent\textbf{(iii) The Dirichlet realization.}
The energy space is Hilbert by~\cite[Proposition~3.1]{ChenChenHauer2026}.
The zero-exterior subspace is closed, since convergence in the energy
norm implies $L^2$ convergence, and contains $C_c^\infty(\Omega)$.
Hence~\eqref{eq:shifted-form} implies that $\mathcal E_{s+\ln}$,
with form domain $\mathcal H_0^{s+\ln}(\Omega)$,
is a densely defined, closed, symmetric and lower-bounded form
on $L^2(\Omega)$. Since $\Omega$ is bounded, every bounded sequence in
$\mathcal H_0^{s+\ln}(\Omega)$ is bounded in $H^s(\mathbb R^n)$
and vanishes outside a fixed ball. The compact embedding
on this ball
\cite[Theorems~5.4 and~7.1]{DiNezzaPalatucciValdinoci2012}
therefore gives
$\mathcal H_0^{s+\ln}(\Omega)\hookrightarrow L^2(\Omega)$.

By the first representation theorem
\cite[Chapter~VI, Theorems~2.1 and~2.6]{Kato1995},
the form $\mathcal E_{s+\ln}$ defines a self-adjoint operator
$(-\Delta)^{s+\ln}_{\Omega}$ on $L^2(\Omega)$.
More precisely, a function
$u\in\mathcal H_0^{s+\ln}(\Omega)$ belongs to its operator domain
if and only if there exists $g\in L^2(\Omega)$ such that
$\mathcal E_{s+\ln}(u,\varphi)=\int_\Omega g\varphi\,dx$
for every $\varphi\in\mathcal H_0^{s+\ln}(\Omega)$.
In this case, $(-\Delta)^{s+\ln}_{\Omega}u=g$.
For this Dirichlet realization and its spectrum, see also
\cite[Remark~1.4 and Theorem~1.5]{ChenChenHauer2026}.

By~\eqref{eq:shifted-form},
$((-\Delta)^{s+\ln}_{\Omega}+\lambda I)^{-1}$
maps $L^2(\Omega)$ boundedly into
$\mathcal H_0^{s+\ln}(\Omega)$, where $I$ denotes the identity operator.
The compact embedding above implies that
$(-\Delta)^{s+\ln}_{\Omega}$ has compact resolvent.
By~\cite[Chapter~III, Theorem~6.29]{Kato1995}, its spectrum consists
of isolated eigenvalues of finite multiplicity. Self-adjointness and
the lower bound imply that they are real and tend to $+\infty$.
If $0\in\rho((-\Delta)^{s+\ln}_{\Omega})$, the resolvent norm formula
\cite[Chapter~V, (3.16)]{Kato1995} gives, with $d_{s,\Omega}$ as in
Corollary~\ref{cor:resolvent},
\begin{equation}
 \bigl\|((-\Delta)^{s+\ln}_{\Omega})^{-1}
 \bigr\|_{L^2(\Omega)\to L^2(\Omega)}
 =d_{s,\Omega}^{-1}.
 \label{eq:resolvent-distance}
\end{equation}

\smallskip
\noindent\textbf{(iv) A zero-exterior Hardy estimate.}
To control the boundary singularities arising from
the action of $(-\Delta)^{s+\ln}$ on the boundary factor
when $\alpha<s$, we use the fractional Hardy inequality
for functions with zero exterior values;
see~\cite[(1.6)--(1.8)]{ChenSong2003Hardy}.
Since $\Omega$ is a bounded Lipschitz open set, every
$v\in H^s(\mathbb R^n)$ vanishing a.e. in $\Omega^c$ satisfies
\begin{equation}
    \|\delta^{-s}v\|_{L^2(\Omega)}
    \leq C\|v\|_{H^s(\mathbb R^n)},
    \label{eq:zero-exterior-hardy}
\end{equation}
where $C$ depends only on $n,s,\Omega$.
Consequently, if $\delta^s g\in L^2(\Omega)$ and
$v\in\mathcal H_0^{s+\ln}(\Omega)$, then
\[
 \int_\Omega|gv|\,dx\leq\|\delta^s g\|_2\|\delta^{-s}v\|_2
 \leq C\|\delta^s g\|_2\|v\|_{H^s}.
\]
For $\Omega_{z,r,Q}=r^{-1}Q^T(\Omega-z)$ with $0<r\leq1$,
a change of variables in~\eqref{eq:zero-exterior-hardy}
gives
\[
    \|\delta_{\Omega_{z,r,Q}}^{-s}v\|_{L^2(\Omega_{z,r,Q})}^2
    \leq C\bigl(r^{2s}\|v\|_2^2+[v]_{H^s}^2\bigr)
    \leq C\|v\|_{H^s}^2
\]
for every $v\in H^s(\mathbb R^n)$ vanishing a.e. in
$\Omega_{z,r,Q}^c$, where
$\delta_{\Omega_{z,r,Q}}=\operatorname{dist}
(\cdot,\partial\Omega_{z,r,Q})$.
In particular, the constant in~\eqref{eq:zero-exterior-hardy}
is uniform in $z,r,Q$.

\subsection{Dilation and local compactness}
\label{sec:local-estimates}

Unlike the fractional Laplacian, $(-\Delta)^{s+\ln}$ is not
homogeneous under dilations. Rescaling the equation at a small
radius produces a linear combination of $(-\Delta)^s$ and
$(-\Delta)^{s+\ln}$. We therefore introduce the rescaled operators
and establish local estimates that are uniform with respect to
the dilation parameter.

For $\ell\geq1$, set
$A_\ell=(-\Delta)^s+(2\ell)^{-1}(-\Delta)^{s+\ln}$.
Its kernel and bilinear form are given by
\[
    J_\ell(h)
    :=\frac{c_{n,s}}{|h|^{n+2s}}
    \left(1+\frac{b_{n,s}-2\ln|h|}{2\ell}\right),
    \qquad h\in\mathbb R^n\setminus\{0\},
\]
and
\[
    \mathcal B_\ell(v,\varphi)
    :=\frac12\iint_{\mathbb R^{2n}}
    (v(x)-v(y))(\varphi(x)-\varphi(y))
    J_\ell(x-y)\,dx\,dy.
\]
To control the contribution from large distances, we define
the weight
\begin{equation}
    \omega_s(y):=
    \frac{1+\ln(1+|y|)}{(1+|y|)^{n+2s}}.
    \label{eq:tail-weight}
\end{equation}

We denote by $\mathcal H_{\mathrm{loc}}^{s+\ln}(\mathbb R^n)$
the space of functions $v$ such that
$\eta v\in\mathcal H^{s+\ln}(\mathbb R^n)$ for every
$\eta\in C_c^\infty(\mathbb R^n)$.
For
$v\in\mathcal H_{\mathrm{loc}}^{s+\ln}(\mathbb R^n)
\cap L^1(\mathbb R^n,\omega_s\,dy)$,
the integral defining $\mathcal B_\ell(v,\varphi)$
is absolutely convergent for every compactly supported
$\varphi\in\mathcal H^{s+\ln}(\mathbb R^n)$.
Indeed, choose a ball $B$ containing the closed
$1$-neighborhood of $\operatorname{supp}\varphi$
and a cutoff $\eta\in C_c^\infty(\mathbb R^n)$
with $\eta=1$ on $B$.
By the preceding kernel estimates and the embedding
$\mathcal H^{s+\ln}(\mathbb R^n)\hookrightarrow H^s(\mathbb R^n)$,
the integral over $B\times B$ is bounded by
$C\|\eta v\|_{\mathcal H^{s+\ln}}
\|\varphi\|_{\mathcal H^{s+\ln}}$.
For $x\in\operatorname{supp}\varphi$ and $y\in B^c$,
we have $|J_\ell(x-y)|\leq C\omega_s(y)$, uniformly in
$\ell\geq1$. Hence the remaining contribution is bounded by
\[
    C\|v\|_{L^2(B)}\|\varphi\|_{L^2(B)}
    +C\|\varphi\|_{L^1(B)}
    \|v\|_{L^1(\mathbb R^n,\omega_s\,dy)}.
\]

The preceding estimates allow us to define weak solutions
without requiring global finite energy.
For an open set $\mathcal O\subset\mathbb R^n$ and
$v\in\mathcal H_{\mathrm{loc}}^{s+\ln}(\mathbb R^n)
\cap L^1(\mathbb R^n,\omega_s\,dy)$,
we say that $A_\ell v=g$ weakly in $\mathcal O$ if
\[
    \mathcal B_\ell(v,\varphi)
    =\int_{\mathcal O}g\varphi\,dx
\]
for every compactly supported
$\varphi\in\mathcal H_0^{s+\ln}(\mathcal O)$,
provided the right-hand side defines a continuous
linear functional on tests with fixed compact support.
The support of $\varphi$ may meet $\partial\mathcal O$.

\smallskip

The following lemma gives the scaling identities for
$(-\Delta)^{s+\ln}$ and $A_\ell$. In particular, rescaling
at radius $r<1$ replaces $\ell$ by $\ell+\ln(1/r)$.

\begin{lemma}
\label{lem:dilation}
For $z\in\mathbb R^n$, $r>0$ and $Q\in O(n)$, define
$T_{z,r,Q}v(y)=v(z+rQy)$. Then
\[
    (-\Delta)^s(T_{z,r,Q}v)
    =r^{2s}T_{z,r,Q}((-\Delta)^sv)
\]
and
\[
    (-\Delta)^{s+\ln}(T_{z,r,Q}v)
    =r^{2s}T_{z,r,Q}((-\Delta)^{s+\ln}v)
    +2\ln r\,(-\Delta)^s(T_{z,r,Q}v).
\]
These identities hold pointwise for $v\in\mathcal S(\mathbb R^n)$
and weakly for
$v\in\mathcal H_{\mathrm{loc}}^{s+\ln}(\mathbb R^n)
\cap L^1(\mathbb R^n,\omega_s\,dy)$.
In particular, for $0<r\leq e^{-1}$ and $\ell_r=\ln(1/r)$,
\[
    A_{\ell_r}(T_{z,r,Q}v)
    =\frac{r^{2s}}{2\ell_r}
    T_{z,r,Q}((-\Delta)^{s+\ln}v).
\]
Moreover, if $A_\ell v=g$ weakly in an open set
$\mathcal O\subset\mathbb R^n$, where $\ell\geq1$,
then for $0<r<1$ and $\ell'=\ell+\ln(1/r)$,
\begin{equation}
    A_{\ell'}(T_{z,r,Q}v)
    =\frac{\ell}{\ell'}r^{2s}T_{z,r,Q}g
    \quad\text{in }r^{-1}Q^T(\mathcal O-z)
    \label{eq:normalized-dilation}
\end{equation}
in the weak sense.
\end{lemma}

\begin{proof}[\textnormal{\textbf{Proof}}]
\noindent\textbf{Step 1: Scaling identities.}
The fractional kernel is homogeneous of degree $-n-2s$,
whereas
\[
    r^{n+2s}K_{s+\ln}(r|h|)
    =K_{s+\ln}(|h|)
    -2c_{n,s}\ln r\,|h|^{-n-2s}.
\]
A change of variables in the singular integral representations,
which preserves spherical principal-value truncations up to
scaling, gives the first two identities for Schwartz functions.
For fixed $z,r,Q$, a change of variables, the comparability
$1+|\ln(r|h|)|\asymp1+|\ln|h||$, and the integrability
of the kernels away from the origin show that $T_{z,r,Q}$
is an isomorphism of $\mathcal H^{s+\ln}(\mathbb R^n)$.
By the cutoff characterization of the local energy space
and $\omega_s(z+rQy)\asymp\omega_s(y)$, it also preserves
$\mathcal H_{\mathrm{loc}}^{s+\ln}(\mathbb R^n)
\cap L^1(\mathbb R^n,\omega_s\,dy)$.

For any compactly supported
$\varphi\in\mathcal H^{s+\ln}(\mathbb R^n)$, set
$w=T_{z,r,Q}v$ and
$\psi(x)=\varphi(r^{-1}Q^T(x-z))$.
Writing $\mathcal E_s$ for the fractional Dirichlet form,
a change of variables gives
\[
    \mathcal E_s(w,\varphi)
    =r^{2s-n}\mathcal E_s(v,\psi)
\]
and
\[
    \mathcal E_{s+\ln}(w,\varphi)
    =r^{2s-n}\mathcal E_{s+\ln}(v,\psi)
    +2\ln r\,\mathcal E_s(w,\varphi).
\]
This proves the weak scaling identities.

\medskip
\noindent
\textbf{Step 2: Scaling of the normalized operators.}
Since $\ell_r+\ln r=0$, the second scaling identity yields
\[
    2\ell_r A_{\ell_r}w
    =r^{2s}T_{z,r,Q}((-\Delta)^{s+\ln}v).
\]
Similarly, using $\ell'+\ln r=\ell$, we obtain
\[
   2\ell'A_{\ell'}w
    =r^{2s}T_{z,r,Q}
    \bigl(2\ell(-\Delta)^sv+(-\Delta)^{s+\ln}v\bigr)=2\ell r^{2s}T_{z,r,Q}(A_\ell v).
\]
To verify~\eqref{eq:normalized-dilation} in the weak sense,
set $\mathcal O'=r^{-1}Q^T(\mathcal O-z)$ and take any
compactly supported
$\varphi\in\mathcal H_0^{s+\ln}(\mathcal O')$.
Then $\psi\in\mathcal H_0^{s+\ln}(\mathcal O)$ has compact
support in $\mathbb R^n$. By the form identities above
and the weak equation for $v$,
\[
    \begin{aligned}
    \mathcal B_{\ell'}(w,\varphi)
    &=\frac{\ell}{\ell'}r^{2s-n}
    \mathcal B_\ell(v,\psi)=\frac{\ell}{\ell'}r^{2s-n}
    \int_{\mathcal O}g\psi\,dx=\frac{\ell}{\ell'}r^{2s}
    \int_{\mathcal O'}(T_{z,r,Q}g)\varphi\,dy.
    \end{aligned}
\]
For tests with fixed compact support, the assumption on $g$
and the boundedness of $T_{z,r,Q}^{-1}$ give
\[
    \left|\int_{\mathcal O'}(T_{z,r,Q}g)\varphi\,dy\right|
    =r^{-n}\left|\int_{\mathcal O}g\psi\,dx\right|
    \leq C\|\psi\|_{\mathcal H^{s+\ln}}
    \leq C\|\varphi\|_{\mathcal H^{s+\ln}}.
\]
Thus~\eqref{eq:normalized-dilation} holds in the weak sense.
\end{proof}

The lower bound for $J_\ell$ near the origin allows us to
derive local $H^s$ estimates uniformly in $\ell\geq1$.
We use these estimates to obtain local compactness of
solutions on varying open sets.

\begin{lemma}
\label{lem:local-compactness}
Let $\ell_j\geq1$ and let $\mathcal O_j\subset\mathbb R^n$
be open. Suppose that
$q_j\in\mathcal H_{\mathrm{loc}}^{s+\ln}(\mathbb R^n)$,
$q_j=0$ a.e. in $\mathcal O_j^c$, and, for every $R>0$,
\[
    \sup_j\left(
    \|q_j\|_{L^2(B_R)}
    +\|q_j\|_{L^1(\mathbb R^n,\omega_s\,dy)}
    \right)<\infty.
\]
Assume that, for every $R>0$ and all sufficiently large $j$, $ \mathcal B_{\ell_j}(q_j,\varphi)
    =\int_{\mathcal O_j}g_j\varphi\,dx$
for every $\varphi\in\mathcal H_0^{s+\ln}(\mathcal O_j)$
with $\operatorname{supp}\varphi\subset B_R$, and
\[
    \left|\int_{\mathcal O_j}g_j\varphi\,dx\right|
    \leq C_R\|\varphi\|_{H^s(\mathbb R^n)},
\]
where $C_R$ is independent of $j$ and $\varphi$. Then $(q_j)$ is bounded in $H^s_{\mathrm{loc}}(\mathbb R^n)$.
Moreover, up to a subsequence, $q_j\to q$ strongly in
$L^2_{\mathrm{loc}}(\mathbb R^n)$ and a.e. in $\mathbb R^n$
for some $q\in H^s_{\mathrm{loc}}(\mathbb R^n)$.
\end{lemma}

\begin{proof}[\textnormal{\textbf{Proof}}]
\noindent\textbf{Step 1: Local energy estimates.}
Fix $\eta\in C_c^\infty(\mathbb R^n)$ and choose $R>0$
such that the closed $1$-neighborhood of
$\operatorname{supp}\eta$ is contained in $B_R$.
Then $\eta q_j\in\mathcal H^{s+\ln}(\mathbb R^n)$ and
$\eta^2q_j\in\mathcal H_0^{s+\ln}(\mathcal O_j)$.

For $r_*$ chosen as above, we have, uniformly in $\ell\geq1$,
\[
    J_\ell(h)\geq c_{n,s}|h|^{-n-2s}
    \quad\text{for }0<|h|<r_*,
\]
and
$|J_\ell(h)|\leq C(1+|\ln|h||)|h|^{-n-2s}$ for $h\neq0$.
In particular,
$\int_{|h|\geq r_*}|J_\ell(h)|\,dh\leq C$.
Since both the contribution of $|x-y|\geq r_*$ to the
$H^s$ seminorm and the corresponding part of
$\mathcal B_\ell(v,v)$ are bounded by $C\|v\|_2^2$,
we obtain
\[
    \mathcal B_\ell(v,v)
    \geq c\|v\|_{H^s(\mathbb R^n)}^2-C\|v\|_2^2,
    \qquad v\in\mathcal H^{s+\ln}(\mathbb R^n),
\]
where $c,C>0$ depend only on $n,s$.

Using the identity
$(a-b)(A^2a-B^2b)=|Aa-Bb|^2-ab|A-B|^2$
for real numbers $a,b,A,B$, we obtain
\[
 \mathcal B_{\ell_j}(q_j,\eta^2q_j)
    =\mathcal B_{\ell_j}(\eta q_j,\eta q_j)-\frac12\iint_{\mathbb R^{2n}}
    q_j(x)q_j(y)|\eta(x)-\eta(y)|^2
    J_{\ell_j}(x-y)\,dx\,dy.
\]
We estimate the last integral on $|x-y|<1$ and
$|x-y|\geq1$ separately.
On $|x-y|<1$, the integrand vanishes unless $x,y\in B_R$.
Using $2|q_j(x)q_j(y)|\leq|q_j(x)|^2+|q_j(y)|^2$,
$|\eta(x)-\eta(y)|\leq\|\nabla\eta\|_\infty|x-y|$
and the upper bound for $J_{\ell_j}$, we obtain
\[
\iint_{|x-y|<1}
    |q_j(x)q_j(y)||\eta(x)-\eta(y)|^2
    |J_{\ell_j}(x-y)|\,dx\,dy
    \leq C_\eta\|q_j\|_{L^2(B_R)}^2.
\]
For $x\in\operatorname{supp}\eta$ and $|x-y|\geq1$,
we have $|J_{\ell_j}(x-y)|\leq C_\eta\omega_s(y)$, and thus
\[
    \iint_{|x-y|\geq1}
    |q_j(x)q_j(y)||\eta(x)-\eta(y)|^2
    |J_{\ell_j}(x-y)|\,dx\,dy\leq C_\eta\|q_j\|_{L^1(B_R)}
    \|q_j\|_{L^1(\mathbb R^n,\omega_s\,dy)}.
\]
Both terms are uniformly bounded by the assumptions on
$q_j$. Together with the lower bound for
$\mathcal B_{\ell_j}(\eta q_j,\eta q_j)$ and the local
$L^2$ bound, this gives
\[
    \mathcal B_{\ell_j}(q_j,\eta^2q_j)
    \geq c\|\eta q_j\|_{H^s(\mathbb R^n)}^2-C_\eta.
\]

Testing the equation with $\eta^2q_j$ and using the
boundedness of multiplication by $\eta$ on $H^s(\mathbb R^n)$,
we obtain
\[
    \bigl|\mathcal B_{\ell_j}(q_j,\eta^2q_j)\bigr|
    \leq C_R\|\eta^2q_j\|_{H^s(\mathbb R^n)}
    \leq C_\eta\|\eta q_j\|_{H^s(\mathbb R^n)}.
\]
Combining the last two estimates and applying Young's
inequality yields $\sup_j\|\eta q_j\|_{H^s(\mathbb R^n)}<\infty$,
since the finitely many remaining terms also belong to
$H^s(\mathbb R^n)$ by the preceding embedding.

\medskip
\noindent\textbf{Step 2: Compactness.}
For each $R>0$, choose $\eta\in C_c^\infty(B_{2R})$
with $\eta=1$ on $B_R$. Step~1 yields
$\sup_j\|q_j\|_{H^s(B_R)}<\infty$.
Since $B_R$ is a bounded Lipschitz domain,
\cite[Theorems~5.4 and~7.1]{DiNezzaPalatucciValdinoci2012}
gives the compact embedding
$H^s(B_R)\hookrightarrow L^2(B_R)$.
By weak compactness in $H^s(B_R)$ and a diagonal argument
on an increasing sequence of balls, there exist
$q\in H^s_{\mathrm{loc}}(\mathbb R^n)$ and a subsequence
such that $q_j\to q$ strongly in $L^2_{\mathrm{loc}}(\mathbb R^n)$.
Passing to a further subsequence gives $q_j\to q$ a.e.
in $\mathbb R^n$.
\end{proof}

Taking $\mathcal O_j=\mathbb R^n$ gives the corresponding
interior compactness result.
The assumption on $g_j$ holds if
$\sup_j\|g_j\|_{L^2(\mathcal O_j\cap B_R)}<\infty$
for every $R>0$.
It also holds if
$\sup_j\|\delta_j^s g_j\|_{L^2(\mathcal O_j\cap B_R)}<\infty$
and the Hardy estimate~\eqref{eq:zero-exterior-hardy}
holds uniformly in $j$, where
$\delta_j=\operatorname{dist}(\cdot,\partial\mathcal O_j)$.

The energy estimate in Step~1 only requires the equation
and the right-hand side estimate for the test $\eta^2q_j$.
We use this local form of the estimate in the proof of
Proposition~\ref{prop:endpoint-excess}.

\subsection{Local boundedness}

We next derive a local $L^\infty$ estimate for $A_\ell$,
uniformly for sufficiently large $\ell$.
The estimate applies also to balls intersecting the boundary,
without any regularity assumption on the open set.
We follow the De Giorgi iteration in~\cite[Section~6.2]
{Cozzi2017}, estimating separately the contribution of
the sign-changing kernel outside the ball.

\begin{lemma}
\label{lem:local-boundedness}
There exists $\ell_0=\ell_0(n,s)\geq1$ such that the following
holds for every $\ell\geq\ell_0$.
Let $\mathcal O\subset\mathbb R^n$ be open and $q\in\mathcal H_{\mathrm{loc}}^{s+\ln}(\mathbb R^n)
    \cap L^1(\mathbb R^n,\omega_s\,dy).$
Suppose that $q=0$ a.e. in $B_2\setminus\mathcal O$
and $A_\ell q=g$ weakly in $\mathcal O\cap B_2$,
where $g\in L^\infty(\mathcal O\cap B_2)$. Then
\begin{equation}
    \|q\|_{L^\infty(B_1)}
    \leq C\left(
    \|q\|_{L^2(B_2)}
    +\|q\|_{L^1(\mathbb R^n,\omega_s\,dy)}
    +\|g\|_{L^\infty(\mathcal O\cap B_2)}
    \right),
    \label{eq:local-boundedness}
\end{equation}
where $C$ depends only on $n,s$.

The same estimate holds for $(-\Delta)^s$ in place of
$A_\ell$, with $H^s_{\mathrm{loc}}(\mathbb R^n)$ in place
of $\mathcal H_{\mathrm{loc}}^{s+\ln}(\mathbb R^n)$,
$\omega_s^0(y)=(1+|y|)^{-n-2s}$ in place of $\omega_s$,
and the corresponding weak formulation in $H^s$.
\end{lemma}

\begin{proof}[\textnormal{\textbf{Proof}}]
\noindent\textbf{Step 1: Truncation estimates.}
Choose $\ell_0=1+|b_{n,s}|+2\ln4$. Then, for
$\ell\geq\ell_0$ and $0<|h|<4$,
$J_\ell(h)\geq(c_{n,s}/2)|h|^{-n-2s}$.
Moreover, uniformly in $\ell\geq1$,
$|J_\ell(h)|\leq C(1+|\ln|h||)|h|^{-n-2s}$ for $h\neq0$.
Fix $\sigma\in\{-1,1\}$ and $k\geq0$, and set
$v=(\sigma q-k)_+$, where $a_+=\max\{a,0\}$.
For $1\leq r<R\leq2$, set $a=R-r$ and choose
$\eta\in C_c^\infty(B_{(r+R)/2})$ with $0\leq\eta\leq1$,
$\eta=1$ on $B_r$ and $|\nabla\eta|\leq C/a$.
Since truncation by $(\cdot-k)_+$ is Lipschitz and
preserves zero for $k\geq0$, we have
$v\in\mathcal H_{\mathrm{loc}}^{s+\ln}(\mathbb R^n)$.
The zero-exterior condition on $q$ then gives
$\eta^2v\in\mathcal H_0^{s+\ln}(\mathcal O\cap B_2)$.
Thus $\eta^2v$ is an admissible test for the equation
of $\sigma q$.

Since $\eta^2v=0$ outside $B_R$, by symmetry we write
$\mathcal B_\ell(\sigma q,\eta^2v)=I_{\rm loc}+I_{\rm far}$,
where
\[I_{\rm loc}:=\frac12\iint_{B_R\times B_R}
    (\sigma q(x)-\sigma q(y))
    (\eta(x)^2v(x)-\eta(y)^2v(y))
    J_\ell(x-y)\,dx\,dy,\]
    and 
\[
  I_{\rm far}
    :=\int_{B_R}\int_{B_R^c}
    (\sigma q(x)-\sigma q(y))\eta(x)^2v(x)
    J_\ell(x-y)\,dy\,dx.
\]
We first obtain a lower bound for $I_{\rm loc}$ and
then estimate $|I_{\rm far}|$.

For $p,t\in\mathbb R$, write $v_p=(p-k)_+$ and
$v_t=(t-k)_+$. Decomposing $p=v_p+\min\{p,k\}$ and
$t=v_t+\min\{t,k\}$, we obtain
\[
    (p-t)(\eta(x)^2v_p-\eta(y)^2v_t)
    \geq|\eta(x)v_p-\eta(y)v_t|^2
    -v_pv_t|\eta(x)-\eta(y)|^2.
\]
Taking $p=\sigma q(x)$ and $t=\sigma q(y)$, and using
the lower bound for $J_\ell$ on $0<|h|<4$, we obtain
\[
    \begin{aligned}
    I_{\rm loc}
    &\geq\frac{c_{n,s}}4
    \iint_{B_R\times B_R}
    \frac{|\eta(x)v(x)-\eta(y)v(y)|^2}
    {|x-y|^{n+2s}}\,dx\,dy\\
    &\quad-\frac12\iint_{B_R\times B_R}
    v(x)v(y)|\eta(x)-\eta(y)|^2
    |J_\ell(x-y)|\,dx\,dy.
    \end{aligned}
\]
To estimate the last integral, we use
$2v(x)v(y)\leq v(x)^2+v(y)^2$ and
\[
    \sup_x\int_{B_4}
    |\eta(x+h)-\eta(x)|^2|J_\ell(h)|\,dh
    \leq Ca^{-2}\int_0^4
    t^{1-2s}(1+|\ln t|)\,dt
    \leq Ca^{-2}.
\]
It follows that
\[
    I_{\rm loc}
    \geq c[\eta v]_{H^s(B_R)}^2
    -Ca^{-2}\|v\|_{L^2(B_R)}^2.
\]

We next estimate $I_{\rm far}$ without using the sign
of $J_\ell$.
For $x\in\operatorname{supp}\eta$ and $y\in B_R^c$,
we have $|x-y|\geq a/2$. Consequently,
\[
    |J_\ell(x-y)|\leq Ca^{-N}\omega_s(y),
    \qquad N=n+4.
\]
Indeed, for $|y|\leq4$, the kernel bound and $0<a\leq1$
give $|J_\ell(x-y)|\leq Ca^{-n-2s-1}\leq Ca^{-N}$,
while $\omega_s$ is bounded below on $B_4$.
For $|y|>4$, the estimate follows from
$|x-y|\asymp1+|y|$.
Since $\sigma q=v+k$ on $\{v>0\}$ and
$\omega_s\in L^1(\mathbb R^n)$, we obtain
\[
    \begin{aligned}
    |I_{\rm far}|
    &\leq Ca^{-N}
    \int_{B_R}\int_{B_R^c}
    (v(x)+k+|q(y)|)v(x)\omega_s(y)\,dy\,dx\\
    &\leq Ca^{-N}\left[
    \|v\|_{L^2(B_R)}^2+
    \left(k+\|q\|_{L^1(\mathbb R^n,\omega_s\,dy)}\right)
    \|v\|_{L^1(B_R)}
    \right].
    \end{aligned}
\]

We now use the weak equation to control the local
$H^s$ seminorm. Testing the equation for $\sigma q$
with $\eta^2v$ gives $I_{\rm loc}+I_{\rm far}
    =\int_{\mathcal O\cap B_2}\sigma g\eta^2v\,dx.$
Combining this identity with the lower bound for
$I_{\rm loc}$, we obtain
\[
    \begin{aligned}
    c[\eta v]_{H^s(B_R)}^2
    &\leq
    \left|\int_{\mathcal O\cap B_2}\sigma g\eta^2v\,dx\right|
    +|I_{\rm far}|+Ca^{-2}\|v\|_{L^2(B_R)}^2\\
    &\leq
    \|g\|_{L^\infty(\mathcal O\cap B_2)}
    \|v\|_{L^1(B_R)}
    +|I_{\rm far}|+Ca^{-2}\|v\|_{L^2(B_R)}^2.
    \end{aligned}
\]
Setting $F=\|q\|_{L^1(\mathbb R^n,\omega_s\,dy)}
+\|g\|_{L^\infty(\mathcal O\cap B_2)}$ and inserting
the estimate for $|I_{\rm far}|$ gives
\[
    [\eta v]_{H^s(B_R)}^2
    \leq Ca^{-N}\left(
    \|v\|_{L^2(B_R)}^2
    +(F+k)\|v\|_{L^1(B_R)}
    \right).
\]

Finally, to obtain the full $H^s(\mathbb R^n)$ norm,
we use the zero values of $\eta v$ outside $B_R$ to write
\[
    [\eta v]_{H^s(\mathbb R^n)}^2
    =[\eta v]_{H^s(B_R)}^2
    +2\int_{B_R}\int_{B_R^c}
    \frac{|\eta(x)v(x)|^2}{|x-y|^{n+2s}}\,dy\,dx.
\]
Since $\operatorname{dist}(\operatorname{supp}\eta,B_R^c)
\geq a/2$, the last integral is bounded by
$Ca^{-2s}\|v\|_{L^2(B_R)}^2$.
Combining these estimates with
$\|\eta v\|_2\leq\|v\|_{L^2(B_R)}$ and
$1+a^{-2s}\leq2a^{-N}$, we conclude that
\begin{equation}
    \|\eta v\|_{H^s(\mathbb R^n)}^2
    \leq Ca^{-N}\left(
    \|v\|_{L^2(B_R)}^2
    +(F+k)\|v\|_{L^1(B_R)}
    \right).
    \label{eq:local-truncation-energy}
\end{equation}

\medskip
\noindent\textbf{Step 2: Sobolev estimates.}
Set $\theta=\min\{s/2,n/4\}$ and
$\varkappa=n/(n-2\theta)>1$, so that
$0<\theta<s$ and $2\theta<n$.
The embedding $H^s(\mathbb R^n)\hookrightarrow
H^\theta(\mathbb R^n)$ and the fractional Sobolev
inequality~\cite[Proposition~2.1 and Theorem~6.5]
{DiNezzaPalatucciValdinoci2012} give $ \|w\|_{L^{2\varkappa}(\mathbb R^n)}^2
    \leq C\|w\|_{H^s(\mathbb R^n)}^2$
for every compactly supported $w\in H^s(\mathbb R^n)$.
Applying this inequality to $w=\eta v$ and using
\eqref{eq:local-truncation-energy}, we obtain
\[
    \|v\|_{L^{2\varkappa}(B_r)}^2
    \leq\frac{C}{(R-r)^N}
    \left(
    \|v\|_{L^2(B_R)}^2
    +(F+k)\|v\|_{L^1(B_R)}
    \right).
\]

\medskip
\noindent\textbf{Step 3: De Giorgi iteration.}
Let $K>0$ with $K\geq F$. For $j\geq0$, set
$R_j=1+2^{-j}$, $k_j=K(1-2^{-j})$,
$v_j=(\sigma q-k_j)_+$ and
$Y_j=\|v_j\|_{L^2(B_{R_j})}^2$.
Since $k_{j+1}-k_j=K2^{-j-1}$, on $\{v_{j+1}>0\}$
we have $v_j\geq K2^{-j-1}$ and $v_{j+1}\leq v_j$.
Consequently,
\[
    |\{v_{j+1}>0\}\cap B_{R_j}|
    \leq\frac{4^{j+1}}{K^2}Y_j,
    \qquad
    \|v_{j+1}\|_{L^1(B_{R_j})}
    \leq\frac{2^{j+1}}{K}Y_j.
\]
Applying the estimate in Step~2 with $r=R_{j+1}$,
$R=R_j$ and $v=v_{j+1}$, and using
$F+k_{j+1}\leq2K$, we obtain
\[
    \|v_{j+1}\|_{L^{2\varkappa}(B_{R_{j+1}})}^2
    \leq C2^{(N+1)j}Y_j.
\]
Let $\vartheta=1-1/\varkappa>0$. By H\"older's inequality,
\[ Y_{j+1}
    \leq
    \|v_{j+1}\|_{L^{2\varkappa}(B_{R_{j+1}})}^2
    |\{v_{j+1}>0\}\cap B_{R_{j+1}}|^\vartheta\leq Cb^jK^{-2\vartheta}Y_j^{1+\vartheta},\]
where $b=2^{N+1+2\vartheta}>1$. Thus $Z_j=K^{-2}Y_j$ satisfies
$Z_{j+1}\leq Cb^jZ_j^{1+\vartheta}$.

Choose $\varepsilon>0$ such that
$C\varepsilon^\vartheta\leq b^{-1/\vartheta}$.
If $Z_0\leq\varepsilon$, induction gives
$Z_j\leq\varepsilon b^{-j/\vartheta}$.
Since $Y_0\leq\|q\|_{L^2(B_2)}^2$, choosing
$K=C(\|q\|_{L^2(B_2)}+F)$ with $C$ sufficiently large
ensures both $K\geq F$ and $Z_0\leq\varepsilon$.
If $\|q\|_{L^2(B_2)}+F=0$, the conclusion holds directly.
As $B_1\subset B_{R_j}$ and $k_j\uparrow K$, we have
$\|(\sigma q-K)_+\|_{L^2(B_1)}^2\leq Y_j\to0$.
Applying this argument with $\sigma=1$ and $\sigma=-1$
proves~\eqref{eq:local-boundedness}.

For the fractional Laplacian, the kernel
$c_{n,s}|h|^{-n-2s}$ satisfies the same local lower bound
and the tail estimate with $\omega_s^0$ in place of $\omega_s$.
The preceding argument applies with $H^s$ in place of
$\mathcal H^{s+\ln}$ and yields the asserted estimate.
\end{proof}

Applying Lemma~\ref{lem:local-boundedness} at a fixed scale,
we obtain the global $L^\infty$ estimate in
Corollary~\ref{cor:boundedness}, with a constant depending
only on $n,s$.

\begin{proof}[\textnormal{\textbf{Proof of the boundedness estimate in Corollary~\ref{cor:boundedness}}}]
Fix $r=e^{-\ell_0}$, where $\ell_0$ is given by
Lemma~\ref{lem:local-boundedness}.
For $x\in\mathbb R^n$, set $v_x(y)=u(x+ry)$ and
$\Omega_x=r^{-1}(\Omega-x)$.
By Lemma~\ref{lem:dilation},
$v_x\in\mathcal H_0^{s+\ln}(\Omega_x)$ and
\[
    A_{\ell_0}v_x
    =\frac{r^{2s}}{2\ell_0}f(x+r\,\cdot)
    \quad\text{weakly in }\Omega_x.
\]
In particular, the equation holds for all tests in
$\mathcal H_0^{s+\ln}(\Omega_x\cap B_2)$.
By a change of variables and the Cauchy--Schwarz inequality,
\[
    \|v_x\|_{L^2(B_2)}
    \leq r^{-n/2}\|u\|_2,\qquad
    \|v_x\|_{L^1(\mathbb R^n,\omega_s\,dy)}
    \leq r^{-n/2}\|\omega_s\|_2\|u\|_2,
\]
where $\omega_s\in L^2(\mathbb R^n)$.
Applying Lemma~\ref{lem:local-boundedness}, we obtain
\[
   \|u\|_{L^\infty(B_r(x))}
    =\|v_x\|_{L^\infty(B_1)}\leq C\left(
    r^{-n/2}(1+\|\omega_s\|_2)\|u\|_2
    +\frac{r^{2s}}{2\ell_0}\|f\|_{L^\infty(\Omega)}
    \right)
    \leq C_0\mathcal A,
\]
where $C_0$ depends only on $n,s$, since $r=e^{-\ell_0}$
is fixed. Covering $\overline\Omega$ by finitely many balls
of radius $r$ yields
$\|u\|_{L^\infty(\Omega)}\leq C_0\mathcal A$.
\end{proof}

\subsection{Uniform interior estimates}

We next establish an interior H\"older estimate uniformly
in $\ell\geq\ell_0$. We adapt the least-squares blow-up
argument of~\cite[Propositions~3.2--3.4]{RosOtonSerra2016Stable}.
By Lemma~\ref{lem:dilation}, rescaling at radius
$r<1$ replaces $\ell$ by $\ell+\ln(1/r)$.
The local compactness estimate and suitable tail bounds
allow us to pass to the fractional Laplacian in the limit
and apply the whole-space Liouville theorem.

\begin{lemma}
\label{lem:interior-holder}
Let $0<\beta<\min\{1,2s\}$ and $\ell\geq\ell_0$.
Suppose that $v\in\mathcal H_{\mathrm{loc}}^{s+\ln}(\mathbb R^n)\cap L^1(\mathbb R^n,\omega_s\,dy)$
satisfies $A_\ell v=g$ weakly in $B_4$, where
$g\in L^\infty(B_4)$. Then $v$ admits a $C^\beta$
representative in $B_1$ and
\begin{equation}
    \|v\|_{C^\beta(B_1)}
    \leq C\left(
    \|v\|_{L^2(B_4)}
    +\|v\|_{L^1(\mathbb R^n,\omega_s\,dy)}
    +\|g\|_{L^\infty(B_4)}
    \right),
    \label{eq:interior-holder}
\end{equation}
where $C$ depends only on $n,s,\beta$.
\end{lemma}

\begin{proof}[\textnormal{\textbf{Proof}}]
By linearity, we may assume that
\[
    \|v\|_{L^2(B_4)}
    +\|v\|_{L^1(\mathbb R^n,\omega_s\,dy)}
    +\|g\|_{L^\infty(B_4)}=1,
\]
since the case in which this sum vanishes is immediate.
As $\omega_s(\cdot+y)\leq C\omega_s$ uniformly for
$y\in B_2$, Lemma~\ref{lem:local-boundedness} applied
on translated balls gives $\|v\|_{L^\infty(B_2)}\leq C$.

\medskip
\noindent\textbf{Step 1: Mean oscillation and blow-up sequences.}
We first prove that, for every $y\in B_1$ and $0<h\leq1/4$,
\begin{equation}
    \inf_{m\in\mathbb R}
    \left(\dashint_{B_h(y)}|v-m|^2\,dx\right)^{1/2}
    \leq Ch^\beta,
    \label{eq:interior-mean-oscillation}
\end{equation}
where $C$ is independent of $\ell$ and $v$.

Suppose that~\eqref{eq:interior-mean-oscillation} fails.
Then there exist normalized solutions $v_j$ of
$A_{\ell_j}v_j=g_j$ in $B_4$, with $\ell_j\geq\ell_0$,
points $y_j\in B_1$, and radii $\tau_j\in(0,1/4)$ such that $M_j:=\sup_{\tau_j\leq h\leq1/4}
    \frac{E_j(h)}{h^\beta}\rightarrow\infty,$
where
\[
    E_j(h):=
    \left(\dashint_{B_h(y_j)}
    |v_j-(v_j)_{B_h(y_j)}|^2\,dx\right)^{1/2},
    \qquad
    (v_j)_{B_h(y_j)}:=\dashint_{B_h(y_j)}v_j\,dx.
\]
Here we used that the average $(v_j)_{B_h(y_j)}$
minimizes the squared distance to constants in $L^2(B_h(y_j))$. Choose $r_j\in[\tau_j,1/4]$ such that
$\Theta_j:=E_j(r_j)/r_j^\beta\geq M_j/2$.
Since $\|v_j\|_{L^\infty(B_2)}\leq C$ and
$\Theta_j\to\infty$, we have $r_j\to0$.
Set $\varepsilon_j=E_j(r_j)=\Theta_jr_j^\beta$,
$m_j=(v_j)_{B_{r_j}(y_j)}$, and define
\[
    w_j(Y):=\frac{v_j(y_j+r_jY)-m_j}{\varepsilon_j}.
\]
By construction,
\begin{equation}
    \dashint_{B_1}w_j\,dY=0,
    \qquad
    \dashint_{B_1}|w_j|^2\,dY=1.
    \label{eq:interior-blowup-normalization}
\end{equation}
Moreover, the choice of $r_j$ gives
$E_j(h)\leq2\Theta_jh^\beta$ for every
$r_j\leq h\leq1/4$.

\medskip
\noindent\textbf{Step 2: Growth and tail estimates.}
We first transfer the mean oscillation estimates to
the rescaled functions $w_j$.
For $r_j\leq h\leq k\leq\min\{2h,1/4\}$,
the Cauchy--Schwarz inequality gives
\[
    \left|(v_j)_{B_h(y_j)}-(v_j)_{B_k(y_j)}\right|
    \leq\left(\frac{k}{h}\right)^{n/2}E_j(k)
    \leq C\Theta_jk^\beta.
\]
For $1\leq T\leq1/(4r_j)$, summing over the radii
$r_j,2r_j,\ldots,2^mr_j,Tr_j$, where
$2^m\leq T<2^{m+1}$, yields
\[
    \left|(v_j)_{B_{Tr_j}(y_j)}-m_j\right|
    \leq C\Theta_j(Tr_j)^\beta,
\]
since $\beta>0$. Together with the bound for $E_j(Tr_j)$,
this gives
\[
    \left(\dashint_{B_{Tr_j}(y_j)}
    |v_j-m_j|^2\,dx\right)^{1/2}
    \leq C\varepsilon_jT^\beta.
\]
Changing variables, we obtain
\begin{equation}
    \|w_j\|_{L^2(B_T)}
    \leq CT^{n/2+\beta},
    \qquad
    1\leq T\leq\frac1{4r_j}.
    \label{eq:interior-blowup-growth}
\end{equation}

To pass to the limiting equation, we also need control
of the weighted tails of $w_j$.
For $k\geq0$ with $2^k<1/(4r_j)$, set
$T=\min\{2^{k+1},1/(4r_j)\}$.
By the Cauchy--Schwarz inequality
and~\eqref{eq:interior-blowup-growth},
\[
    \int_{B_T\setminus B_{2^k}}
    |w_j(Y)|\omega_s(Y)\,dY
    \leq C(1+k)2^{-k(2s-\beta)}.
\]
The right-hand side is summable in $k$ since $\beta<2s$.

For the remaining part of the tail, observe that
$|x-y_j|\asymp1+|x|$ whenever $|x-y_j|\geq1/4$,
uniformly for $y_j\in B_1$. Hence
\[
    r_j^{-n}\omega_s\left(\frac{x-y_j}{r_j}\right)
    \leq Cr_j^{2s}(1+\ln(1/r_j))\omega_s(x).
\]
Since $|m_j|\leq C$, $\omega_s\in L^1(\mathbb R^n)$
and $\|v_j\|_{L^1(\mathbb R^n,\omega_s\,dx)}\leq1$,
a change of variables gives
\begin{align*}
    \int_{|Y|>1/(4r_j)}|w_j(Y)|\omega_s(Y)\,dY
    &\leq
    \frac{Cr_j^{2s}(1+\ln(1/r_j))}{\varepsilon_j}
    \int_{\mathbb R^n}(|v_j|+|m_j|)\omega_s\,dx\\
    &\leq
    \frac{Cr_j^{2s-\beta}(1+\ln(1/r_j))}{\Theta_j}
    \longrightarrow0.
\end{align*}
Combining these estimates with the normalization
on $B_1$, we obtain
\begin{equation}
    \sup_j\|w_j\|_{L^1(\mathbb R^n,\omega_s\,dY)}<\infty,
    \qquad
    \lim_{M\to\infty}\limsup_{j\to\infty}
    \int_{|Y|>M}|w_j|\omega_s\,dY=0.
    \label{eq:interior-blowup-tails}
\end{equation}

\medskip
\noindent\textbf{Step 3: Compactness and the limiting equation.}
Set $\ell_j'=\ell_j+\ln(1/r_j)$ and
$\mathcal O_j'=r_j^{-1}(B_4-y_j)$.
By~\eqref{eq:normalized-dilation} and the fact that
$A_\ell$ annihilates constants,
\[
    A_{\ell_j'}w_j(Y)
    =\frac{\ell_j}{\ell_j'}
    \frac{r_j^{2s}}{\varepsilon_j}
    g_j(y_j+r_jY)
    \quad\text{weakly in }\mathcal O_j'.
\]
Since $r_j\to0$ and $y_j\in B_1$, every fixed ball is
contained in $\mathcal O_j'$ for all sufficiently large $j$.
Moreover, $\ell_j'\to\infty$ and
\[
    \frac{\ell_j}{\ell_j'}
    \frac{r_j^{2s}}{\varepsilon_j}
    \leq\frac{r_j^{2s-\beta}}{\Theta_j}
    \rightarrow0.
\]
Thus the right-hand sides, extended by zero outside
$\mathcal O_j'$, tend to zero in
$L^\infty_{\mathrm{loc}}(\mathbb R^n)$.

By Lemma~\ref{lem:dilation}, each $w_j$ belongs to
$\mathcal H_{\mathrm{loc}}^{s+\ln}(\mathbb R^n)$.
Estimates~\eqref{eq:interior-blowup-growth}
and~\eqref{eq:interior-blowup-tails}, together with
the local $L^\infty$ bounds on the right-hand sides,
allow us to apply Lemma~\ref{lem:local-compactness}
with the open sets taken to be $\mathbb R^n$.
After passing to a subsequence,
$w_j\to w$ strongly in $L^2_{\mathrm{loc}}(\mathbb R^n)$
and a.e. in $\mathbb R^n$, where
$w\in H^s_{\mathrm{loc}}(\mathbb R^n)$.
Moreover, for every $T\geq1$, $\|w\|_{L^2(B_T)}\leq CT^{n/2+\beta}.$

We next identify the equation satisfied by $w$.
For every $\varphi\in C_c^\infty(\mathbb R^n)$,
the singular integral representations give
\begin{equation}
    |(-\Delta)^s\varphi(Y)|\leq C_\varphi\omega_s^0(Y),
    \qquad
    |(-\Delta)^{s+\ln}\varphi(Y)|
    \leq C_\varphi\omega_s(Y).
    \label{eq:test-function-tail}
\end{equation}
On bounded sets, these estimates follow from the bound
$|\varphi(Y+h)+\varphi(Y-h)-2\varphi(Y)|
\leq C_\varphi|h|^2$ and the integrability of
$t^{1-2s}(1+|\ln t|)$ near zero.
Away from the support of $\varphi$, they follow from
the decay of the corresponding kernels.
The same estimates hold uniformly for the integrals
truncated at $|h|=\varepsilon$, with $0<\varepsilon<1$.

By symmetry of the truncated kernels, followed by
dominated convergence using~\eqref{eq:test-function-tail}
and the absolute convergence of the forms, we obtain
\[
    \mathcal B_{\ell_j'}(w_j,\varphi)
    =\int_{\mathbb R^n}w_j(-\Delta)^s\varphi\,dY
    +\frac1{2\ell_j'}
    \int_{\mathbb R^n}w_j(-\Delta)^{s+\ln}\varphi\,dY.
\]
The first integral converges to
$\int_{\mathbb R^n}w(-\Delta)^s\varphi\,dY$
by strong local $L^2$ convergence
and~\eqref{eq:interior-blowup-tails}.
The second term tends to zero
by~\eqref{eq:test-function-tail}, the uniform weighted
$L^1$ bound and $\ell_j'\to\infty$.
Since the rescaled right-hand sides tend to zero locally,
we conclude that
\[
    \int_{\mathbb R^n}w(-\Delta)^s\varphi\,dY=0
    \qquad
    \text{for all }\varphi\in C_c^\infty(\mathbb R^n).
\]
By Fatou's lemma,
$w\in L^1(\mathbb R^n,\omega_s\,dY)
\subset L^1(\mathbb R^n,\omega_s^0\,dY)$.
Since $w\in H^s_{\mathrm{loc}}(\mathbb R^n)$,
the fractional form is continuous on $H^s$ tests with
fixed compact support, by the same local and tail
estimates used for $\mathcal B_\ell$.
Approximating compactly supported $H^s$ tests by
mollification therefore gives $(-\Delta)^sw=0$
in the corresponding weak sense in $\mathbb R^n$.

\medskip
\noindent\textbf{Step 4: The Liouville theorem.}
To apply the Liouville theorem, we first upgrade
the $L^2$ growth estimate for $w$ to an $L^\infty$
growth estimate.
For $T\geq1$, set $w_T(X)=T^{-\beta}w(TX)$.
Then $(-\Delta)^sw_T=0$ weakly in $\mathbb R^n$ and
\[
    \|w_T\|_{L^2(B_R)}
    \leq CR^{n/2+\beta}
    \qquad\text{for every }R\geq1.
\]
By the Cauchy--Schwarz inequality,
\[
    \int_{B_{2^{k+1}}\setminus B_{2^k}}
    |w_T|\omega_s^0\,dX
    \leq C2^{-k(2s-\beta)}
    \qquad\text{for every }k\geq0.
\]
Since $\beta<2s$, summing over $k$ yields a uniform
bound for $\|w_T\|_{L^1(\mathbb R^n,\omega_s^0\,dX)}$.
The fractional version of
Lemma~\ref{lem:local-boundedness} gives
$\|w_T\|_{L^\infty(B_1)}\leq C$, and hence
$\|w\|_{L^\infty(B_T)}\leq CT^\beta$ for every $T\geq1$.

The operator $-(-\Delta)^s$ is a symmetric stable operator
with constant positive spectral density and satisfies
the ellipticity conditions of~\cite{RosOtonSerra2016Stable}.
As $\beta<2s$, the whole-space Liouville theorem
\cite[Theorem~2.1]{RosOtonSerra2016Stable}
implies that $w$ is a polynomial of degree at most
$\lfloor\beta\rfloor=0$.
On the other hand, strong convergence in $L^2(B_1)$
and~\eqref{eq:interior-blowup-normalization} give
\[
    \dashint_{B_1}w\,dY=0,
    \qquad
    \dashint_{B_1}|w|^2\,dY=1.
\]
This contradicts the constancy of $w$ and
proves~\eqref{eq:interior-mean-oscillation}.

\medskip
\noindent\textbf{Step 5: H\"older regularity.}
By~\eqref{eq:interior-mean-oscillation},
\[
    \inf_{m\in\mathbb R}
    \int_{B_h(y)}|v-m|^2\,dx
    \leq Ch^{n+2\beta}
\]
for every $y\in B_1$ and $0<h\leq1/4$.
Together with $\|v\|_{L^\infty(B_2)}\leq C$,
the Campanato characterization of H\"older spaces
\cite[Theorem~I.2]{Campanato1963}, with $p=2$ and
$\lambda=n+2\beta$, gives $\|v\|_{C^\beta(B_1)}\leq C$.
Undoing the normalization proves~\eqref{eq:interior-holder}.
\end{proof}

\section{The boundary factor and its operator image}
\label{sec:mellin}

Using the half-line power formula for $(-\Delta)^s$,
we prove that $(-\Delta_{\mathbb R})^{s+\ln}\Phi_s$
extends smoothly to zero from the right and establish
the asymptotic expansion and dilation estimates for $\Phi_s$.
We then construct a regularized distance $D$ satisfying
$\nabla D=\nu_{\rm in}$ on $\partial\Omega$ and set
$\Psi_s=\zeta\Phi_s(D)$.
We estimate $(-\Delta)^{s+\ln}\Psi_s$ in the three cases
$\alpha<s$, $\alpha=s$ and $\alpha>s$, with an additional
local $L^2$ oscillation estimate when $\alpha=s$.
Finally, we prove that $\Psi_s\in\mathcal H_0^{s+\ln}(\Omega)$
and identify its quadratic form pairing with its pointwise
operator image.

\subsection{The one-dimensional boundary factor}

Throughout this section, set $b=r_\chi/4$, where $r_\chi$
and $\chi$ are as in the construction of $\Phi_s$
in~\eqref{eq:boundary-factor}.
For $\gamma>-1$, write $t_+^\gamma=t^\gamma$ if $t>0$
and $t_+^\gamma=0$ otherwise. We first recall the action of the fractional Laplacian
on half-line powers. For $0<q<1$, $-1<\gamma<2q$ and $t>0$,
Fall and Ros-Oton~\cite[Lemma~2.4, (2.8)]{FallRosOton2022}
give
\begin{equation}
    (-\Delta_{\mathbb R})^q t_+^\gamma
    =F(q,\gamma)t^{\gamma-2q},
    \qquad
    F(q,\gamma)
    =-\frac{\Gamma(1+\gamma)\Gamma(2q-\gamma)}{\pi}
    \sin\bigl(\pi(\gamma-q)\bigr).
    \label{eq:half-line-coefficient}
\end{equation}
The coefficient is rewritten using the gamma reflection
formula~\cite[(5.5.3)]{NISTDLMF}.
The conditions $\gamma>-1$ and $\gamma<2q$ ensure
integrability near the origin and at infinity, respectively,
while $t_+^\gamma$ is smooth near each $t>0$.
Thus~\eqref{eq:half-line-coefficient} holds in the usual
principal-value sense.

Differentiating~\eqref{eq:half-line-coefficient} with respect
to $q$, keeping $\gamma$ fixed, gives the coefficients
used below to construct the Mellin density $m_s$.

\begin{lemma}
\label{lem:half-line-powers}
For $0<a<s$, set $d_s(a)=F(s,s+a)$ and
\[
    e_s(a)=d_s(a)
    \bigl(2\psi(s-a)-\pi\cot(\pi a)\bigr),
\]
where $\psi=\Gamma'/\Gamma$.
Then, for every $t>0$,
\[
    (-\Delta_{\mathbb R})^{s+\ln}t_+^{s+a}
    =t^{a-s}\bigl(e_s(a)-2d_s(a)\ln t\bigr).
\]
\end{lemma}

\begin{proof}[\textnormal{\textbf{Proof}}]
Fix $\gamma=s+a$ and $t>0$, and choose
$0<q_-<s<q_+<1$ with $\gamma<2q_-$.
For $q\in[q_-,q_+]$, differentiating the fractional kernel gives
\[
    \partial_q\bigl(c_{1,q}|h|^{-1-2q}\bigr)
    =c_{1,q}|h|^{-1-2q}(b_{1,q}-2\ln|h|).
\]
Using the symmetric second difference of $t_+^\gamma$,
which is $O(|h|^2)$ near $h=0$, the differentiated
integrand is bounded by
$C|h|^{1-2q_+}(1+|\ln|h||)$ near $h=0$ and by
$C|h|^{\gamma-1-2q_-}(1+\ln|h|)$ at infinity.
Both bounds are integrable since $q_+<1$ and $\gamma<2q_-$.
On the remaining bounded intervals, the integrand is
integrable since $\gamma>-1$.
Hence differentiation under the integral sign is justified,
and~\eqref{eq:half-line-coefficient} gives
\[
  (-\Delta_{\mathbb R})^{s+\ln}t_+^{s+a}=\left.\partial_q
    \bigl(F(q,\gamma)t^{\gamma-2q}\bigr)\right|_{q=s}=t^{a-s}\left(
    \partial_qF(s,s+a)-2F(s,s+a)\ln t
    \right).
\]
Differentiating the expression for $F$ in
\eqref{eq:half-line-coefficient}, we obtain
\[
    \partial_qF(s,s+a)
    =F(s,s+a)
    \bigl(2\psi(s-a)-\pi\cot(\pi a)\bigr)
    =e_s(a).
\]
The conclusion follows from the definition of $d_s(a)$.
\end{proof}

We now explain the choice of the Mellin density $m_s$
in~\eqref{eq:boundary-factor}. For $a\in(0,1)\setminus\{s\}$,
define
\[
    \mathfrak G_s(a)
    :=-\frac{\Gamma(s+1)\Gamma(s-a)}{\pi}
    \sqrt{\sin(\pi a)}.
\]
For $0<a<s$, we have $\mathfrak G_s=d_sm_s$ and
\[
    \mathfrak G_s'(a)
    =\mathfrak G_s(a)
    \left(-\psi(s-a)+\frac{\pi}{2}\cot(\pi a)\right).
\]
Thus, by the expression for $e_s$ in
Lemma~\ref{lem:half-line-powers},
\begin{equation}
    2\mathfrak G_s'(a)+e_s(a)m_s(a)=0,
    \qquad 0<a<s.
    \label{eq:mellin-density-identity}
\end{equation}
In particular, for $0<a<s$ and $t>0$,
\[
    m_s(a)(-\Delta_{\mathbb R})^{s+\ln}t_+^{s+a}
    =-2\partial_a\bigl(\mathfrak G_s(a)t^{a-s}\bigr).
\]
This identity is the reason for the choice of $m_s$.

The density satisfies $m_s(a)\asymp[a(1-a)]^{-1/2}$.
Since $\mathfrak G_s(a)/\sqrt a$ and
$\mathfrak G_s(a)/\sqrt{1-a}$ extend smoothly to $a=0$
and $a=1$, respectively, we have
$\mathfrak G_s(a)=O(a^{1/2})$ and
$\mathfrak G_s'(a)=O(a^{-1/2})$ as $a\downarrow0$,
and $\mathfrak G_s(a)=O((1-a)^{1/2})$ and
$\mathfrak G_s'(a)=O((1-a)^{-1/2})$ as $a\uparrow1$.
The only pole of $\mathfrak G_s$ in $(0,1)$ is at $a=s$.
Using $\Gamma(z)=z^{-1}+O(1)$ as $z\to0$, we obtain
\begin{equation}
    \mathfrak G_s(a)
    =\frac{\Gamma(s+1)\sqrt{\sin(\pi s)}}{\pi(a-s)}
    +O(1)
    \quad\text{as }a\to s.
    \label{eq:mellin-interior-pole}
\end{equation}
The residue of $\mathfrak G_s(a)t^{a-s}$ at $a=s$
is independent of $t>0$. We therefore consider differences
at two positive points to remove this pole before
integrating in $a$.

\begin{proposition}
\label{prop:mellin-factor}
The function $\Phi_s$ defined in~\eqref{eq:boundary-factor}
is nonnegative, compactly supported, and belongs to
$C^s(\mathbb R)\cap C^\infty((0,\infty))$.
It vanishes on $(-\infty,0]$ and is positive on $(0,r_\chi)$.
Moreover, there exists $H_s\in C^\infty((-b,b))$ such that
\begin{equation}
    (-\Delta_{\mathbb R})^{s+\ln}\Phi_s(t)=H_s(t)
    \qquad\text{for }0<t<b.
    \label{eq:factor-flat-image}
\end{equation}

As $t\downarrow0$, we have $\Phi_s(t)\sim\frac{t^s}{\sqrt{\ln(1/t)}}.$
For every integer $j\geq0$, there exists $C_j>0$ such that,
for $0<t<b$,
\begin{equation}
    |\Phi_s^{(j)}(t)|
    \leq C_j\frac{t^{s-j}}{\sqrt{\ln(1/t)}}.
    \label{eq:factor-derivatives}
\end{equation}
Furthermore, there exist $r_\Phi\in(0,b)$ and $C>0$
such that, for every $x,y\in\mathbb R$ with
$0<|x-y|<r_\Phi$,
\begin{equation}
    |\Phi_s(x)-\Phi_s(y)|
    \leq C\Phi_s(|x-y|)
    \leq C\frac{|x-y|^s}{\sqrt{\ln(e/|x-y|)}}.
    \label{eq:factor-small-modulus}
\end{equation}
\end{proposition}

\begin{proof}[\textnormal{\textbf{Proof}}]
The support and sign properties follow from
\eqref{eq:boundary-factor}. Differentiation under the
Mellin integral also gives
$\Phi_s\in C^\infty((0,\infty))$.

\medskip
\noindent\textbf{Step 1: Cancellation for the truncated powers.}
For $0<a<1$, set $p_a(t)=\chi(t)t_+^{s+a}$ and fix
$t_0=b/2$. Write
$K^{(1)}_{s+\ln}(r)
=c_{1,s}(b_{1,s}-2\ln|r|)|r|^{-1-2s}$ for $r\neq0$.
Since $\chi=1$ on $(-r_\chi,r_\chi)$,
Lemma~\ref{lem:half-line-powers} gives, for $0<a<s$
and $0<t<b$,
\[ (-\Delta_{\mathbb R})^{s+\ln}p_a(t)
    =t^{a-s}\bigl(e_s(a)-2d_s(a)\ln t\bigr)+\int_0^\infty
    (1-\chi(y))y^{s+a}K^{(1)}_{s+\ln}(y-t)\,dy.\]
The last integral is absolutely convergent since
$1-\chi$ vanishes near zero and $a<s$.

For $|t|<b$, define
\[
    T(a,t):=\int_0^\infty(1-\chi(y))y^{s+a}
    \bigl(K^{(1)}_{s+\ln}(y-t)
    -K^{(1)}_{s+\ln}(y-t_0)\bigr)\,dy.
\]
Subtracting the operator values at $t$ and $t_0$
and using~\eqref{eq:mellin-density-identity}, we obtain
\begin{equation}
     m_s(a)\bigl((-\Delta_{\mathbb R})^{s+\ln}p_a(t)
    -(-\Delta_{\mathbb R})^{s+\ln}p_a(t_0)\bigr)=-2\partial_a\left[
    \mathfrak G_s(a)(t^{a-s}-t_0^{a-s})\right]
    +m_s(a)T(a,t).
    \label{eq:mellin-difference-identity}
\end{equation}
We next extend this identity from $0<a<s$ to $0<a<1$.

Since $1-\chi(y)=0$ for $0<y<r_\chi$ and
$b=r_\chi/4$, the mean value theorem gives, for
$y\geq r_\chi$ and $|t|<b$,
\[
    \left|K^{(1)}_{s+\ln}(y-t)
    -K^{(1)}_{s+\ln}(y-t_0)\right|
    \leq C|t-t_0|y^{-2-2s}(1+|\ln y|).
\]
Thus, for fixed $t\in(-b,b)$, the integral defining
$T(a,t)$ converges locally uniformly in the complex
half-plane $\{\operatorname{Re}a<s+1\}$.
Indeed, at infinity its integrand is bounded by
$Cy^{\operatorname{Re}a-s-2}(1+\ln y)$.
Hence $T(\cdot,t)$ is holomorphic in this half-plane.

Furthermore, for every integer $k\geq1$,
\[
    |\partial_t^kK^{(1)}_{s+\ln}(y-t)|
    \leq C_k y^{-1-2s-k}(1+|\ln y|)
\]
for $y\geq r_\chi$ and $|t|<b$.
Together with the kernel difference estimate, this gives
\[
    \sup_{0\leq a\leq1}
    \|T(a,\cdot)\|_{C^k((-b,b))}\leq C_k
    \qquad\text{for every integer }k\geq0.
\]

For each fixed $t>0$, the map
$a\mapsto(-\Delta_{\mathbb R})^{s+\ln}p_a(t)$
is holomorphic on $\{\operatorname{Re}a>-1-s\}$.
Indeed, $p_a$ has fixed compact support, and its behavior
near zero is integrable when $\operatorname{Re}a>-1-s$.
Near the evaluation point $t$, its symmetric second
differences are $O(|h|^2)$, locally uniformly in $a$.
These estimates provide an integrable majorant on
every compact subset of the parameter half-plane.

Finally, by~\eqref{eq:mellin-interior-pole} and
\[
    t^{a-s}-t_0^{a-s}
    =(a-s)\ln(t/t_0)+O((a-s)^2),
\]
the function
$\mathfrak G_s(a)(t^{a-s}-t_0^{a-s})$
has a removable singularity at $a=s$.
Since $m_s$ is real analytic on $(0,1)$, both sides of
\eqref{eq:mellin-difference-identity} are real analytic
on $(0,1)$ after removing this singularity.
The identity theorem therefore gives
\eqref{eq:mellin-difference-identity} for every
$a\in(0,1)$ and $0<t<b$.

\medskip
\noindent\textbf{Step 2: Smooth extension of
$(-\Delta_{\mathbb R})^{s+\ln}\Phi_s$.}
By definition, $\Phi_s(t)=\int_0^1m_s(a)p_a(t)\,da$.
For fixed $t>0$, the second derivatives of $p_a$ are
uniformly bounded near $t$ for $0\leq a\leq1$. Hence
\[
    |2p_a(t)-p_a(t+h)-p_a(t-h)|
    \leq C_t|h|^2
    \qquad\text{for }|h|<t/2.
\]
Since $m_s\in L^1(0,1)$ and $0<s<1$, this bound
gives integrability against
$m_s(a)|K^{(1)}_{s+\ln}(h)|$ near $h=0$.
For $|h|\geq t/2$, the uniform boundedness of $p_a$
and the integrability of the kernel give another
integrable majorant. Fubini's theorem therefore yields
\[
    (-\Delta_{\mathbb R})^{s+\ln}\Phi_s(t)
    =\int_0^1m_s(a)
    (-\Delta_{\mathbb R})^{s+\ln}p_a(t)\,da.
\]

For every fixed $t>0$, the endpoint estimates for
$\mathfrak G_s$ give
\[
    \lim_{a\downarrow0}\mathfrak G_s(a)
    (t^{a-s}-t_0^{a-s})
    =
    \lim_{a\uparrow1}\mathfrak G_s(a)
    (t^{a-s}-t_0^{a-s})=0.
\]
The derivative of this function is integrable at
both endpoints and, by Step~1, extends analytically
across $a=s$. Integrating
\eqref{eq:mellin-difference-identity} over $(0,1)$
therefore gives, for $0<t<b$,
\[
    (-\Delta_{\mathbb R})^{s+\ln}\Phi_s(t)
    -(-\Delta_{\mathbb R})^{s+\ln}\Phi_s(t_0)
    =\int_0^1m_s(a)T(a,t)\,da.
\]
Define, for $|t|<b$,
\[
    H_s(t):=(-\Delta_{\mathbb R})^{s+\ln}\Phi_s(t_0)
    +\int_0^1m_s(a)T(a,t)\,da.
\]
The uniform derivative bounds for $T$ and
$m_s\in L^1(0,1)$ imply that $H_s\in C^\infty((-b,b))$.
This proves~\eqref{eq:factor-flat-image}.

\medskip
\noindent\textbf{Step 3: Asymptotics and derivative estimates.}
For $0<t<b$, we have $\chi(t)=1$.
Setting $\tau=\ln(1/t)$, we write
\[
    \Phi_s(t)=t^s\int_0^1m_s(a)e^{-a\tau}\,da.
\]
By the definition of $m_s$,
$\lim_{a\downarrow0}\sqrt a\,m_s(a)=1/\sqrt\pi$.
Moreover, since $m_s\in L^1(1/2,1)$,
\[
    \sqrt\tau\int_{1/2}^1m_s(a)e^{-a\tau}\,da
    \leq C\sqrt\tau e^{-\tau/2}\longrightarrow0
    \quad\text{as }\tau\to\infty.
\]
On $(0,1/2)$, the change of variables $u=a\tau$ gives
\[
    \sqrt\tau\int_0^{1/2}m_s(a)e^{-a\tau}\,da
    =\int_0^{\tau/2}
    \frac{m_s(u/\tau)}{\sqrt\tau}e^{-u}\,du.
\]
Since $m_s(a)\leq Ca^{-1/2}$ on $(0,1/2)$,
the integrand, extended by zero for $u\geq\tau/2$,
is bounded by $Cu^{-1/2}e^{-u}$.
By dominated convergence,
\[
    \lim_{\tau\to\infty}
    \sqrt\tau\int_0^1m_s(a)e^{-a\tau}\,da
    =\frac1{\sqrt\pi}\int_0^\infty u^{-1/2}e^{-u}\,du
    =1.
\]
Consequently,
$\Phi_s(t)\sim t^s/\sqrt{\ln(1/t)}$ as $t\downarrow0$.

For each integer $j\geq0$, differentiation of $t^{s+a}$
gives a factor $t^{-j}$ and a polynomial in $s+a$,
bounded uniformly for $0\leq a\leq1$.
The preceding asymptotic estimate therefore gives
\[
    |\Phi_s^{(j)}(t)|
    \leq C_jt^{s-j}\int_0^1m_s(a)e^{-a\tau}\,da
    \leq C_j\frac{t^{s-j}}{\sqrt{\ln(1/t)}}
\]
for all sufficiently small $t>0$.
By smoothness on $(0,\infty)$, the estimate holds on
$(0,b)$ after increasing $C_j$.
This proves~\eqref{eq:factor-derivatives}.

\medskip
\noindent\textbf{Step 4: Modulus of continuity.}
Let $h=|x-y|$ with $0<h<r_\chi$.
We may assume that at least one of $x$ and $y$ belongs
to $\operatorname{supp}\chi$, since otherwise the
difference vanishes. Thus both points lie in a fixed
bounded interval.

For $0<a\leq1-s$, subadditivity gives
$|x_+^{s+a}-y_+^{s+a}|\leq h^{s+a}$.
For $1-s<a<1$, the functions $t_+^{s+a}$ are uniformly
Lipschitz on this interval.
Using also the smoothness of $\chi$ and
$m_s\in L^1(0,1)$, we obtain
\[
    |\Phi_s(x)-\Phi_s(y)|
    \leq C\int_0^{1-s}m_s(a)h^{s+a}\,da+Ch
    \leq C\bigl(\Phi_s(h)+h\bigr),
\]
where the last inequality follows from $\chi(h)=1$.
Since $s<1$, the asymptotic formula in Step~3 gives
$h=o(\Phi_s(h))$ as $h\downarrow0$.
Hence, for some $r_\Phi\in(0,b)$,
\[
    |\Phi_s(x)-\Phi_s(y)|
    \leq C\Phi_s(h)
    \leq C\frac{h^s}{\sqrt{\ln(e/h)}}
    \qquad\text{for }0<h<r_\Phi.
\]
This proves~\eqref{eq:factor-small-modulus}.
Since $\Phi_s$ is bounded, the same estimate also yields
$\Phi_s\in C^s(\mathbb R)$.
\end{proof}

We next reduce the operator acting on functions of the
normal coordinate to the one-dimensional operator.
This reduction will be used to compare the boundary factor
on curved domains with its half-space counterpart.

For $n\geq2$, $0<q<1$ and $r\neq0$, polar integration
and the beta integral~\cite[(5.12.1) and (5.12.3)]{NISTDLMF}
give
\[
 c_{n,q}\int_{\mathbb R^{n-1}}
    (|z'|^2+r^2)^{-(n+2q)/2}\,dz'=c_{n,q}\pi^{(n-1)/2}
    \frac{\Gamma(q+1/2)}{\Gamma(n/2+q)}
    |r|^{-1-2q}=c_{1,q}|r|^{-1-2q}.
\]
For fixed $r\neq0$, differentiation under the integral
sign is justified by the integrability of
$(1+|z'|)^{-n-2q_-}(1+\ln(1+|z'|))$ on $\mathbb R^{n-1}$,
where $q_-$ is the lower endpoint of a compact subinterval
of $(0,1)$ containing $s$ in its interior.
Differentiating at $q=s$, we obtain
\[
    \int_{\mathbb R^{n-1}}
    K_{s+\ln}\left(\sqrt{|z'|^2+r^2}\right)\,dz'
    =K^{(1)}_{s+\ln}(r).
\]
Moreover, the change of variables $z'=|r|\zeta'$ gives
\[
    \int_{\mathbb R^{n-1}}
    \left|K_{s+\ln}\left(\sqrt{|z'|^2+r^2}\right)\right|\,dz'
    \leq C\frac{1+|\ln|r||}{|r|^{1+2s}}.
\]

Let $v$ be $C^2$ near $t$ and satisfy
\[
    \int_{\mathbb R}
    |v(u)|\frac{1+\ln(1+|u|)}{(1+|u|)^{1+2s}}\,du<\infty.
\]
We claim that
\begin{equation}
    (-\Delta)^{s+\ln}\bigl(v(x_n)\bigr)(te_n)
    =(-\Delta_{\mathbb R})^{s+\ln}v(t).
    \label{eq:tangential-reduction}
\end{equation}
Indeed, the second difference
$2v(t)-v(t+r)-v(t-r)$ is $O(r^2)$ near $r=0$.
Together with the absolute kernel estimate above, this gives
integrability near $r=0$, since
$\int_0^1r^{1-2s}(1+|\ln r|)\,dr<\infty$.
Away from zero, integrability follows from the weighted
$L^1$ assumption on $v$. Fubini's theorem and the kernel
identity therefore yield
\[
    \begin{aligned}
    (-\Delta)^{s+\ln}\bigl(v(x_n)\bigr)(te_n)&=\frac12\int_{\mathbb R}
    \bigl(2v(t)-v(t+r)-v(t-r)\bigr)
    \int_{\mathbb R^{n-1}}
    K_{s+\ln}\left(\sqrt{|z'|^2+r^2}\right)\,dz'\,dr\\
    \quad&=\frac12\int_{\mathbb R}
    \bigl(2v(t)-v(t+r)-v(t-r)\bigr)
    K^{(1)}_{s+\ln}(r)\,dr=(-\Delta_{\mathbb R})^{s+\ln}v(t).
    \end{aligned}
\]
For $n=1$, no tangential integration is needed.

\subsection{Asymptotic expansion and dilation}

We derive the first two logarithmic correction terms
in the asymptotic expansion of $\Phi_s$.
They will be used to show that normalization by
$t^s/\sqrt{\ln(R/t)}$ need not yield a H\"older continuous
quotient. We also establish dilation estimates for $\Phi_s$
and its derivatives for the boundary blow-up argument.

\begin{proposition}
\label{prop:factor-asymptotics}
For any fixed $R>0$, set $L_R(t)=\ln(R/t)$.
As $t\downarrow0$,
\begin{equation}
    \Phi_s(t)
    =\frac{t^s}{\sqrt{L_R(t)}}
    \left(
    1+\frac{A_{1,s,R}}{L_R(t)}
    +\frac{A_{2,s,R}}{L_R(t)^2}
    +O(L_R(t)^{-3})
    \right),
    \label{eq:factor-watson}
\end{equation}
where $A_{1,s,R}=(\ln R-\psi(s+1))/2$ and
\[
    A_{2,s,R}
    =\frac38\left(
    (\ln R-\psi(s+1))^2
    +\frac{\pi^2}{6}-\psi_1(s+1)
    \right)>0,
\]
with $\psi_1=\psi'$ denoting the trigamma function.

After reducing $r_\Phi\in(0,b)$ if necessary, for every
$\varepsilon>0$ there exists $C_\varepsilon>0$ such that
\begin{equation}
    T^s\leq\frac{\Phi_s(Tr)}{\Phi_s(r)}
    \leq C_\varepsilon T^{s+\varepsilon}
    \qquad\text{for }T\geq1,\quad 0<Tr\leq r_\Phi.
    \label{eq:factor-potter}
\end{equation}
Moreover, for every $M>0$,
\begin{equation}
    \lim_{r\downarrow0}
    \sup_{0\leq t\leq M}
    \left|\frac{\Phi_s(rt)}{\Phi_s(r)}-t^s\right|=0.
    \label{eq:factor-ratio-limit}
\end{equation}
For every integer $k\geq1$,
\begin{equation}
    \frac{r^k\Phi_s^{(k)}(rt)}{\Phi_s(r)}
    \longrightarrow
    s(s-1)\cdots(s-k+1)t^{s-k}
    \label{eq:factor-derivative-ratio}
\end{equation}
as $r\downarrow0$, locally uniformly in $(0,\infty)$.
\end{proposition}

\begin{proof}[\textnormal{\textbf{Proof}}]
\noindent\textbf{Step 1: Asymptotic expansion.}
Set $\lambda=L_R(t)$ and define
\[
    h(a):=\frac{\Gamma(s+1)R^a}{\Gamma(s+1+a)}
    \sqrt{\frac{\pi a}{\sin(\pi a)}}.
\]
For $0<t<b$, we have $\chi(t)=1$ and hence
\[
    \lambda^{1/2}t^{-s}\Phi_s(t)
    =\frac{\lambda^{1/2}}{\sqrt\pi}
    \int_0^1a^{-1/2}h(a)e^{-a\lambda}\,da.
\]
The function $h$ is analytic near zero, with $h(0)=1$.
Fix a sufficiently small $a_*\in(0,1/2)$.
Since $a^{-1/2}h(a)\in L^1(0,1)$, the contribution
of $(a_*,1)$ to the last expression is
$O(\lambda^{1/2}e^{-a_*\lambda})$ as $\lambda\to\infty$.

On $(0,a_*)$, Taylor's formula gives
$h(a)=1+h'(0)a+\frac12h''(0)a^2+O(a^3)$.
By Watson's lemma~\cite[Section~2.3(ii)]{NISTDLMF}
and the gamma integral
\[
    \int_0^\infty a^{j-1/2}e^{-a\lambda}\,da
    =\Gamma(j+1/2)\lambda^{-j-1/2},
    \qquad j\geq0,
\]
we obtain
\[
    \lambda^{1/2}t^{-s}\Phi_s(t)
    =1+\frac{h'(0)}{2\lambda}
    +\frac{3h''(0)}{8\lambda^2}
    +O(\lambda^{-3}).
\]
Indeed, extending each polynomial term from $(0,a_*)$
to $(0,\infty)$ introduces only an exponentially small
error, while the Taylor remainder contributes at most
\[
    C\lambda^{1/2}\int_0^\infty
    a^{5/2}e^{-a\lambda}\,da
    =O(\lambda^{-3}).
\]

To compute the coefficients, we use
\[
    \ln h(a)=a\ln R+\ln\Gamma(s+1)
    -\ln\Gamma(s+1+a)
    -\frac12\ln\frac{\sin(\pi a)}{\pi a}.
\]
Since $\sin(\pi a)/(\pi a)=1-\pi^2a^2/6+O(a^4)$,
we have $(\ln h)'(0)=\ln R-\psi(s+1)$ and
$(\ln h)''(0)=\pi^2/6-\psi_1(s+1)$.
Using $h(0)=1$, we obtain
\[
h'(0)=\ln R-\psi(s+1),\quad
    h''(0)=(\ln R-\psi(s+1))^2
    +\frac{\pi^2}{6}-\psi_1(s+1).
\]
This proves~\eqref{eq:factor-watson} with the stated
coefficients.

Moreover, the series representation of the trigamma
function~\cite[(5.15.1)]{NISTDLMF} gives
\begin{equation}
    \frac{\pi^2}{6}-\psi_1(s+1)
    =\sum_{j=1}^\infty
    \left(\frac1{j^2}-\frac1{(j+s)^2}\right)>0.
    \label{eq:trigamma-positive}
\end{equation}
Consequently, $A_{2,s,R}>0$ for every $R>0$.

\medskip
\noindent\textbf{Step 2: Dilation estimates.}
After reducing $r_\Phi$, we have
$\ln(1/r_\Phi)\geq1$ and
$\Phi_s(r)\asymp r^s/\sqrt{\ln(1/r)}$
for $0<r\leq r_\Phi$.
For $T\geq1$ and $0<Tr\leq r_\Phi$, positivity of $m_s$
gives
\[
    \Phi_s(Tr)
    =T^sr^s\int_0^1m_s(a)T^ar^a\,da
    \geq T^s\Phi_s(r).
\]
On the other hand, the two-sided estimate for $\Phi_s$
yields
\[
    \begin{aligned}
    \frac{\Phi_s(Tr)}{\Phi_s(r)}
    &\leq CT^s
    \left(\frac{\ln(1/r)}{\ln(1/(Tr))}\right)^{1/2}=CT^s
    \left(1+\frac{\ln T}{\ln(1/(Tr))}\right)^{1/2}\\
    &\leq CT^s(1+\ln T)^{1/2}
    \leq C_\varepsilon T^{s+\varepsilon}.
    \end{aligned}
\]

We next prove~\eqref{eq:factor-ratio-limit}.
For every $0<d<M$, the expansion
\eqref{eq:factor-watson} with $R=1$ gives, uniformly
for $t\in[d,M]$,
\[
    \frac{\Phi_s(rt)}{\Phi_s(r)}
    =t^s\left(
    \frac{\ln(1/r)}{\ln(1/r)+\ln(1/t)}
    \right)^{1/2}
    \left(1+O\left(\frac1{\ln(1/r)}\right)\right).
\]
Since $\ln(1/t)$ is bounded on $[d,M]$, the right-hand
side converges uniformly to $t^s$ as $r\downarrow0$.

To include $t=0$, observe that, for $0\leq t\leq1$
and $0<r\leq r_\Phi$, the integral representation gives
$0\leq\Phi_s(rt)/\Phi_s(r)\leq t^s$.
Thus, for $0<d<\min\{1,M\}$,
\[
    \sup_{0\leq t\leq d}
    \left|\frac{\Phi_s(rt)}{\Phi_s(r)}-t^s\right|
    \leq d^s.
\]
Combining this estimate with the uniform convergence
on $[d,M]$ and letting $d\downarrow0$ proves
\eqref{eq:factor-ratio-limit}.

\medskip
\noindent\textbf{Step 3: Convergence of the rescaled derivatives.}
For each integer $k\geq1$, set
$P_k(a)=(s+a)(s+a-1)\cdots(s+a-k+1)$.
For $0<u\leq r_\Phi$, differentiation under the
Mellin integral gives
\[
    u^k\Phi_s^{(k)}(u)
    =u^s\int_0^1P_k(a)m_s(a)u^a\,da.
\]
Since $P_k(a)-P_k(0)=O(a)$ uniformly for $a\in[0,1]$,
we obtain
\[
    \left|
    \frac{u^k\Phi_s^{(k)}(u)}{\Phi_s(u)}-P_k(0)
    \right|
    \leq C_k
    \frac{\int_0^1am_s(a)u^a\,da}
    {\int_0^1m_s(a)u^a\,da}.
\]
Writing $\tau=\ln(1/u)$ and using
$m_s(a)\asymp a^{-1/2}$ on $(0,1/2)$, we have
\[
    \int_0^1am_s(a)u^a\,da
    \leq C\int_0^{1/2}a^{1/2}e^{-a\tau}\,da
    +Ce^{-\tau/2}
    \leq C\tau^{-3/2}.
\]
Since $\int_0^1m_s(a)u^a\,da\asymp\tau^{-1/2}$,
it follows that
\[
    \frac{u^k\Phi_s^{(k)}(u)}{\Phi_s(u)}
    =P_k(0)+O\left(\frac1{\ln(1/u)}\right).
\]

Taking $u=rt$, we write
\[
    \frac{r^k\Phi_s^{(k)}(rt)}{\Phi_s(r)}
    =t^{-k}
    \frac{(rt)^k\Phi_s^{(k)}(rt)}{\Phi_s(rt)}
    \frac{\Phi_s(rt)}{\Phi_s(r)}.
\]
The preceding estimate and~\eqref{eq:factor-ratio-limit}
give~\eqref{eq:factor-derivative-ratio}, locally uniformly
for $t>0$.
\end{proof}

\subsection{Regularized distance and operator estimates}

When $0<\alpha<1$, the distance function $\delta$
need not be differentiable throughout an interior
neighborhood of $\partial\Omega$.
We therefore construct a regularized distance $D$ that
is smooth off the boundary and satisfies
$\nabla D=\nu_{\rm in}$ on $\partial\Omega$.
Using the signed affine distances to the tangent planes
in the Whitney construction~\cite[Chapter~VI, Sections~1--2]
{Stein1970}, we obtain $D-\delta=O(\delta^{1+\alpha})$
inside $\Omega$. This normalization also ensures that
$\Phi_s(D)/\Phi_s(\delta)$ is $C^{0,\alpha}$ up to the
boundary, allowing us to recover the desired regularity
of $u/\Phi_s(\delta)$ from that of $u/\Phi_s(D)$.

\begin{lemma}
\label{lem:normalized-distance}
Let $\Omega\subset\mathbb R^n$ be a bounded $C^{1,\alpha}$
open set, with $0<\alpha\leq1$.
There exist $r_D>0$ and
\[
    D\in C^{1,\alpha}(\mathcal U)
    \cap C^\infty(\mathcal U\setminus\partial\Omega),
    \qquad
    \mathcal U=\{x:\delta(x)<r_D\},
\]
such that $D=0$ and $\nabla D=\nu_{\rm in}$ on
$\partial\Omega$. Moreover, $D>0$ in $\mathcal U\cap\Omega$
and $D<0$ in $\mathcal U\setminus\overline\Omega$.

For every $x\in\mathcal U\cap\Omega$, with $t=\delta(x)$
and any nearest boundary point $z_x$, we have
\begin{equation}
    |D(x)-t|+t|\nabla D(x)-\nu_{\rm in}(z_x)|
    \leq Ct^{1+\alpha}.
    \label{eq:normalized-distance}
\end{equation}
In particular, $D(x)\asymp t$ and
$|\nabla^kD(x)|\leq C_kt^{1+\alpha-k}$
for every integer $k\geq2$.

The quotient $\Phi_s(D)/\Phi_s(\delta)$ extends to a
$C^{0,\alpha}$ function on a sufficiently small closed
interior boundary collar, with boundary value one.
More precisely, almost everywhere in the interior collar,
\begin{equation}
    \left|\frac{\Phi_s(D)}{\Phi_s(\delta)}-1\right|
    +\delta\left|
    \nabla\left(\frac{\Phi_s(D)}{\Phi_s(\delta)}\right)
    \right|
    \leq C\delta^\alpha.
    \label{eq:distance-factor-ratio}
\end{equation}
If $n=1$, we may take $D$ to be the signed distance
near $\partial\Omega$, in which case
$\Phi_s(D)/\Phi_s(\delta)=1$ in the interior collar.
\end{lemma}

\begin{proof}[\textnormal{\textbf{Proof}}]
If $n=1$, the compact $C^1$ boundary consists of finitely
many isolated endpoints. Taking $D$ to be the signed
distance in disjoint neighborhoods of these endpoints
gives all the assertions. We henceforth assume $n\geq2$.

\noindent\textbf{Step 1: Construction of the regularized distance.}
Since $\partial\Omega$ is compact and of class $C^{1,\alpha}$,
the boundary charts give, for every $z,w\in\partial\Omega$,
\[
    |\nu_{\rm in}(z)-\nu_{\rm in}(w)|
    \leq C|z-w|^\alpha,
    \qquad
    |\nu_{\rm in}(z)\cdot(w-z)|
    \leq C|z-w|^{1+\alpha}.
\]
These are the Whitney compatibility conditions for the
affine functions $\ell_z(x)=(x-z)\cdot\nu_{\rm in}(z)$.

Take a Whitney decomposition of
$\mathbb R^n\setminus\partial\Omega$ and a subordinate
smooth partition of unity $\{\theta_Q\}_Q$;
see~\cite[Chapter~VI, Sections~1--2]{Stein1970}.
For each Whitney cube $Q$, choose a nearest boundary
point $z_Q$ to its center. Define
\[
    D(x)=\sum_Q\theta_Q(x)\ell_{z_Q}(x)
    \quad\text{for }x\notin\partial\Omega,
    \qquad D=0\quad\text{on }\partial\Omega.
\]
The sum is locally finite, so $D$ is smooth off
$\partial\Omega$.

Fix $x\notin\partial\Omega$, set $t=\delta(x)$, and
let $z$ be a nearest boundary point. The cubes with
$x\in\operatorname{supp}\theta_Q$ have bounded overlap,
side lengths comparable to $t$, and satisfy
$|z_Q-z|\leq Ct$ and $|\nabla^k\theta_Q(x)|\leq C_kt^{-k}$.
By the compatibility conditions,
\[
    \begin{aligned}
    |\ell_{z_Q}(x)-\ell_z(x)|
    &\leq |x-z|
    |\nu_{\rm in}(z_Q)-\nu_{\rm in}(z)|
    +|(z-z_Q)\cdot\nu_{\rm in}(z_Q)|\leq Ct^{1+\alpha},
    \end{aligned}
\]
and $|\nabla\ell_{z_Q}-\nabla\ell_z|\leq Ct^\alpha$.
Since $\sum_Q\theta_Q=1$, we write
$D-\ell_z=\sum_Q\theta_Q(\ell_{z_Q}-\ell_z)$.
Differentiating with $z$ fixed gives
\[
    |D(x)-\ell_z(x)|
    +t|\nabla D(x)-\nu_{\rm in}(z)|
    \leq Ct^{1+\alpha}.
\]
For $k\geq2$, the affine functions have no derivatives
of order greater than one. Hence
\[
    |\nabla^kD(x)|
    \leq C_k\bigl(t^{-k}t^{1+\alpha}
    +t^{1-k}t^\alpha\bigr)
    \leq C_kt^{1+\alpha-k}.
\]
For related regularized distance estimates, see
\cite[Theorem~1.3, (1.7b)]{Lieberman1985}.

We next verify the regularity of $D$ across the boundary.
For $w\in\partial\Omega$ and $x\notin\partial\Omega$,
a nearest boundary point $z$ satisfies
$t\leq|x-w|$ and $|z-w|\leq2|x-w|$.
The preceding estimates and the compatibility conditions
give $|D(x)-\ell_w(x)|\leq C|x-w|^{1+\alpha}$.
The same estimate holds for $x\in\partial\Omega$.
Thus $D$ is differentiable at $w$, with
$\nabla D(w)=\nu_{\rm in}(w)$.
Moreover,
\[
    |\nabla D(x)-\nu_{\rm in}(w)|
    \leq Ct^\alpha+C|z-w|^\alpha
    \leq C|x-w|^\alpha.
\]

To prove the $\alpha$-H\"older continuity of $\nabla D$,
let $x,y\notin\partial\Omega$, set $h=|x-y|$, and put
$t=\min\{\delta(x),\delta(y)\}$.
If $h<t/2$, the segment joining $x$ and $y$ has
boundary distance at least $t/2$. The Hessian estimate
then gives
$|\nabla D(x)-\nabla D(y)|\leq Cht^{\alpha-1}\leq Ch^\alpha$.
If $h\geq t/2$, then $\delta(x),\delta(y)\leq3h$ and
$|z_x-z_y|\leq7h$. Consequently,
\[
 |\nabla D(x)-\nabla D(y)|
    \leq C\bigl(\delta(x)^\alpha+\delta(y)^\alpha
    +|z_x-z_y|^\alpha\bigr)\leq Ch^\alpha.
\]
Together with the estimates at boundary points, this
proves $D\in C^{1,\alpha}(\mathcal U)$.

Since $\partial\Omega$ is $C^1$, the segment from a
nearest boundary point $z$ to $x$ is normal to
$\partial\Omega$ at $z$. Its direction is inward if
$x\in\Omega$ and outward if $x\notin\overline\Omega$.
Hence $\ell_z(x)=t$ in the former case and
$\ell_z(x)=-t$ in the latter.
Reducing $r_D$ so that $Ct^\alpha\leq1/2$ for
$t<r_D$, the estimate for $D-\ell_z$ gives the sign
and comparability assertions, as well as
\eqref{eq:normalized-distance}.

\medskip
\noindent\textbf{Step 2: Estimates for the quotient
$\Phi_s(D)/\Phi_s(\delta)$.}
By Propositions~\ref{prop:mellin-factor}
and~\ref{prop:factor-asymptotics}, for sufficiently
small $t>0$,
\[
    |(\ln\Phi_s)'(t)|\leq Ct^{-1},
    \qquad
    |(\ln\Phi_s)''(t)|\leq Ct^{-2}.
\]
Since $D\asymp\delta$, the mean value theorem
and~\eqref{eq:normalized-distance} give
\[
    \left|\ln\frac{\Phi_s(D(x))}{\Phi_s(\delta(x))}\right|
    \leq C\frac{|D(x)-\delta(x)|}{\delta(x)}
    \leq C\delta(x)^\alpha.
\]
Exponentiating, we obtain
\[
    \left|\frac{\Phi_s(D(x))}{\Phi_s(\delta(x))}-1\right|
    \leq C\delta(x)^\alpha.
\]
In particular, the quotient extends continuously to
$\partial\Omega$ with value one.

Since $\delta$ is Lipschitz, it is differentiable a.e.
At every point of differentiability in $\Omega$,
the relation $\delta(\cdot)\leq|\cdot-z_x|$, with
equality at $x$, gives
$\nabla\delta(x)=(x-z_x)/\delta(x)=\nu_{\rm in}(z_x)$.
Thus~\eqref{eq:normalized-distance} implies
$|\nabla D-\nabla\delta|\leq C\delta^\alpha$
a.e. in the interior collar.

Differentiating the logarithm of the quotient yields
\[
\nabla\ln\frac{\Phi_s(D)}{\Phi_s(\delta)}
    =(\ln\Phi_s)'(D)(\nabla D-\nabla\delta)+\bigl((\ln\Phi_s)'(D)
    -(\ln\Phi_s)'(\delta)\bigr)\nabla\delta.
\]
Using $D\asymp\delta$, $|\nabla\delta|=1$ and the
bounds for the first two derivatives of $\ln\Phi_s$,
we obtain
\[
    \left|\nabla\ln\frac{\Phi_s(D)}{\Phi_s(\delta)}\right|
    \leq C\left(
    \frac{|\nabla D-\nabla\delta|}{\delta}
    +\frac{|D-\delta|}{\delta^2}
    \right)
    \leq C\delta^{\alpha-1}.
\]
Since $\Phi_s(D)/\Phi_s(\delta)$ is bounded above
and below by positive constants, this proves
\eqref{eq:distance-factor-ratio}.

\medskip
\noindent\textbf{Step 3: H\"older continuity of the quotient.}
Write $Q=\Phi_s(D)/\Phi_s(\delta)$ in the interior
collar and extend $Q$ by one on $\partial\Omega$.
Choose a smaller closed collar so that the estimates
in Step~2 hold on a larger collar.

Let $x,y$ belong to the smaller collar, set $h=|x-y|$,
and assume that $t=\delta(x)\leq\delta(y)$.
If $h<t/2$, then $B_{t/2}(x)$ lies inside the larger
interior collar, with $\delta\asymp t$ on this ball.
Since $Q$ is locally Lipschitz, estimate
\eqref{eq:distance-factor-ratio} gives
$\|\nabla Q\|_{L^\infty(B_{t/2}(x))}
\leq Ct^{\alpha-1}$.
The $W^{1,\infty}$ estimate on the ball therefore yields
\[
    |Q(x)-Q(y)|
    \leq Cht^{\alpha-1}\leq Ch^\alpha.
\]
If $h\geq t/2$, then $\delta(x),\delta(y)\leq3h$.
By~\eqref{eq:distance-factor-ratio},
\[
    |Q(x)-Q(y)|
    \leq |Q(x)-1|+|Q(y)-1|
    \leq C\bigl(\delta(x)^\alpha+\delta(y)^\alpha\bigr)
    \leq Ch^\alpha.
\]
The latter estimate also applies when either point
lies on $\partial\Omega$. This completes the proof.
\end{proof}

Choose $\zeta\in C_c^\infty(\mathcal U)$ with
$0\leq\zeta\leq1$ and $\zeta=1$ in a smaller neighborhood
of $\partial\Omega$. Extending the product by zero outside
$\mathcal U$, we define
\begin{equation}
    \Psi_s=\zeta\Phi_s(D),\qquad
    G_s=(-\Delta)^{s+\ln}\Psi_s.
    \label{eq:regularized-factor}
\end{equation}
In particular, $\Psi_s=0$ in $\Omega^c$.
We write $L(t)=\ln(e/t)$ for $0<t<1$.

For $z\in\partial\Omega$, set
$\ell_z(y)=(y-z)\cdot\nu_{\rm in}(z)$,
$P_z=\Phi_s(\ell_z)$ and $E_z=\Psi_s-P_z$.
By~\eqref{eq:tangential-reduction}
and~\eqref{eq:factor-flat-image}, we have
\[
    (-\Delta)^{s+\ln}P_z(z+t\nu_{\rm in}(z))
    =H_s(t),\qquad 0<t<b.
\]
We estimate $G_s$ by comparing $\Psi_s$ with $P_z$,
following the affine comparison used in the proofs of
\cite[Propositions~2.3 and~2.6]{RosOtonSerra2017Boundary}.
For the Mellin factor $\Phi_s$, we first estimate its
translated differences in the normal direction and then
sum the resulting kernel integrals over dyadic annuli.

\begin{lemma}
\label{lem:distance-lift}
There exist $r_0>0$ and $C>0$ such that, for every
$z\in\partial\Omega$ and $0<R<r_0$,
\begin{equation}
    \int_{B_{2R}(z)}|E_z|\,dy
    +\int_{B_{2R}(z)}|\Phi_s(\delta_0)-P_z|\,dy
    \leq CR^{n+\alpha}\Phi_s(R).
    \label{eq:geometry-annulus}
\end{equation}
The principal value defining $G_s$ exists at every point
of $\Omega$, and $G_s\in L^\infty_{\rm loc}(\Omega)$.
For every $x\in\Omega$ with $0<\delta(x)<r_0$,
\begin{equation}
    |G_s(x)|\leq C
    \begin{cases}
    \delta(x)^{\alpha-s}\sqrt{L(\delta(x))},
        &\alpha<s,\\
    L(\delta(x))^{3/2},
        &\alpha=s,\\
    1,  &\alpha>s.
    \end{cases}
    \label{eq:lift-image-bound}
\end{equation}
If $n=1$, then $G_s\in L^\infty(\Omega)$ for every
$0<\alpha\leq1$.
\end{lemma}

\begin{proof}[\textnormal{\textbf{Proof}}]
Since $\Psi_s$ is bounded on $\mathbb R^n$ and smooth
in $\Omega$, the symmetric second-difference representation
gives the existence of the principal value defining $G_s$
and its local boundedness in $\Omega$. If $n=1$, Lemma~\ref{lem:normalized-distance} allows
us to take $D$ equal to the signed distance near
each boundary point $z$. Consequently,
$\Psi_s=\Phi_s(\delta_0)=P_z$ in a fixed neighborhood
of $z$. In particular, both integrals in
\eqref{eq:geometry-annulus} vanish for sufficiently small $R$.
For $x\in\Omega$ sufficiently close to $z$, the integral
defining $(-\Delta)^{s+\ln}E_z(x)$ is separated from
the kernel singularity and is uniformly bounded.
Since
$(-\Delta)^{s+\ln}P_z(x)=H_s(\delta(x))$
and $\partial\Omega$ consists of finitely many points,
we obtain $G_s\in L^\infty(\Omega)$.
We henceforth assume $n\geq2$.

\medskip
\noindent\textbf{Step 1: Integral comparison with $P_z$.}
Fix $z\in\partial\Omega$.
Since $D(z)=0$ and $\nabla D(z)=\nu_{\rm in}(z)$,
Lemma~\ref{lem:normalized-distance} gives
\[
    |D(y)-\ell_z(y)|\leq CR^{1+\alpha}
    \qquad\text{for }y\in B_{2R}(z).
\]
We choose $r_0>0$ sufficiently small so that the estimates
in this step hold for $0<R<4r_0$, uniformly in $z$,
and $\zeta=1$ on $B_{2R}(z)$.

We claim that, with $\varepsilon=CR^{1+\alpha}$ and
$C$ sufficiently large,
\[
    (\ell_z(y)-\varepsilon)_+
    \leq\delta_0(y)
    \leq(\ell_z(y)+\varepsilon)_+
    \qquad\text{for }y\in B_{2R}(z).
\]
After a translation and a rotation, we may assume that
$z=0$, $\nu_{\rm in}(z)=e_n$, and $\Omega$ locally lies
above the graph of a $C^{1,\alpha}$ function $\varphi$
satisfying $\varphi(0)=0$ and $\nabla\varphi(0)=0$.
In particular, $|\varphi(y')|\leq CR^{1+\alpha}$
for $|y'|\leq4R$.

For $y\in B_{2R}\cap\Omega$, vertical projection onto
the graph gives
$\delta(y)\leq y_n-\varphi(y')\leq y_n+\varepsilon$.
On the other hand, any nearest boundary point $w$ to $y$
belongs to $B_{4R}$, since $\delta(y)\leq|y|$.
Thus $|w_n|\leq\varepsilon$ and
$\delta(y)=|y-w|\geq y_n-\varepsilon$.
If $y\in B_{2R}\setminus\Omega$,
$\delta_0(y)=0$ and $y_n\leq\varphi(y')\leq\varepsilon$.
This proves the claim.

The integral representation of $\Phi_s$ shows that
it is increasing on $(0,r_\chi)$, and
$\Phi_s=0$ on $(-\infty,0]$.
After reducing $r_0$, we may assume that
$\varepsilon\leq R$ and $3R\leq r_\Phi$ for $R<4r_0$.
Hence both $|E_z(y)|$ and
$|\Phi_s(\delta_0(y))-P_z(y)|$ are bounded by
\[
    \Phi_s(\ell_z(y)+\varepsilon)
    -\Phi_s(\ell_z(y)-\varepsilon)
    \qquad\text{for }y\in B_{2R}(z).
\]
Integrating in the normal direction gives
\[
  \int_{-2R}^{2R}
    \bigl(\Phi_s(t+\varepsilon)
    -\Phi_s(t-\varepsilon)\bigr)\,dt=\int_{2R-\varepsilon}^{2R+\varepsilon}
    \Phi_s(t)\,dt
    \leq2\varepsilon\Phi_s(2R+\varepsilon)
    \leq C\varepsilon\Phi_s(R),
\]
where the last inequality follows from
\eqref{eq:factor-potter}.
Since the tangential sections of $B_{2R}(z)$ have
$(n-1)$-dimensional measure at most $CR^{n-1}$,
Fubini's theorem yields
\[
 \int_{B_{2R}(z)}|E_z|\,dy
    +\int_{B_{2R}(z)}|\Phi_s(\delta_0)-P_z|\,dy\leq CR^{n-1}\varepsilon\Phi_s(R)
    \leq CR^{n+\alpha}\Phi_s(R).
\]

\medskip
\noindent\textbf{Step 2: Estimate near the evaluation point.}
Let $x\in\Omega$, set $t=\delta(x)<r_0$, and choose
a nearest boundary point $z\in\partial\Omega$.
Since $x=z+t\nu_{\rm in}(z)$, we have
\[
    G_s(x)=H_s(t)+(-\Delta)^{s+\ln}E_z(x).
\]
We may assume that $r_0<b/2$, so that $H_s$ is uniformly
bounded on $(0,r_0)$.

In $B_{t/4}(x)$, we have $\delta\asymp t$,
$D\asymp t$ and $\ell_z\asymp t$.
Since $B_{t/4}(x)\subset B_{5t/4}(z)$,
the $C^{1,\alpha}$ regularity of $D$ and
Lemma~\ref{lem:normalized-distance} give
\[
    |D-\ell_z|+t|\nabla D-\nu_{\rm in}(z)|
    \leq Ct^{1+\alpha},
    \qquad |\nabla^2D|\leq Ct^{\alpha-1}
\]
throughout $B_{t/4}(x)$.
Since $\zeta=1$ there, the mean value theorem and
\eqref{eq:factor-derivatives} give
\[
    |E_z(x)|\leq Ct^{1+\alpha}
    \sup_{u\asymp t}|\Phi_s'(u)|
    \leq C\frac{t^{s+\alpha}}{\sqrt{L(t)}}.
\]

Writing $\nu=\nu_{\rm in}(z)$, we have
\[
\nabla^2E_z
    =\bigl(\Phi_s''(D)-\Phi_s''(\ell_z)\bigr)
    \nabla D\otimes\nabla D+\Phi_s''(\ell_z)
    \bigl(\nabla D\otimes\nabla D-\nu\otimes\nu\bigr)+\Phi_s'(D)\nabla^2D.
\]
By the preceding estimates and
\eqref{eq:factor-derivatives}, in $B_{t/4}(x)$,
\[
  |\nabla^2E_z|
    \leq C\left(
    t^{1+\alpha}\sup_{u\asymp t}|\Phi_s'''(u)|
    +t^\alpha\sup_{u\asymp t}|\Phi_s''(u)|
    +t^{\alpha-1}\sup_{u\asymp t}|\Phi_s'(u)|
    \right)\leq C\frac{t^{s+\alpha-2}}{\sqrt{L(t)}}.
\]
Consequently,
\begin{equation}
    |E_z(x)|+t^2
    \sup_{B_{t/4}(x)}|\nabla^2E_z|
    \leq C\frac{t^{s+\alpha}}{\sqrt{L(t)}}.
    \label{eq:flat-error-local}
\end{equation}

By symmetry of the kernel and Taylor's formula,
\[
    \begin{aligned}
    &\left|\operatorname{p.v.}
    \int_{|h|<t/4}
    (E_z(x)-E_z(x+h))K_{s+\ln}(|h|)\,dh\right|\\
    \qquad\leq& \,C
    \sup_{B_{t/4}(x)}|\nabla^2E_z|
    \int_0^{t/4}r^{1-2s}L(r)\,dr\\
    \qquad\leq &\,C
    \frac{t^{s+\alpha-2}}{\sqrt{L(t)}}
    t^{2-2s}L(t)
    =Ct^{\alpha-s}\sqrt{L(t)}.
    \end{aligned}
\]

\medskip
\noindent\textbf{Step 3: Annular estimates.}
By~\eqref{eq:flat-error-local} and the kernel bound,
\[
  |E_z(x)|\int_{|h|\geq t/4}
    |K_{s+\ln}(|h|)|\,dh
    \leq C\frac{t^{s+\alpha}}{\sqrt{L(t)}}
    t^{-2s}L(t)=Ct^{\alpha-s}\sqrt{L(t)}.
\]
It remains to estimate the contribution of $E_z(x+h)$.

For $t/4\leq R<r_0$, we have $B_{2R}(x)\subset B_{6R}(z)$.
Applying~\eqref{eq:geometry-annulus} at scale $3R$
and using~\eqref{eq:factor-potter}, we obtain
\[
    \int_{R<|h|<2R}|E_z(x+h)|\,dh
    \leq CR^{n+\alpha}\Phi_s(R).
\]
Together with the kernel bound and
$\Phi_s(R)\asymp R^s/\sqrt{L(R)}$, this gives
\[
     \int_{R<|h|<2R}
    |E_z(x+h)||K_{s+\ln}(|h|)|\,dh\leq CR^{-n-2s}L(R)
    R^{n+\alpha}\Phi_s(R)
    \leq CR^{\alpha-s}\sqrt{L(R)}.
\]
The contribution of $|h|\geq r_0$ is uniformly bounded,
since $E_z$ is uniformly bounded and the kernel is
absolutely integrable away from zero.

Summing over the dyadic annuli with radii
$R_j=2^jt/4<r_0$, and observing that
\[
    R_j^{\alpha-s}\sqrt{L(R_j)}
    \leq C\int_{R_j}^{2R_j}
    R^{\alpha-s-1}\sqrt{L(R)}\,dR,
\]
we obtain
\[
    \int_{t/4<|h|<r_0}
    |E_z(x+h)||K_{s+\ln}(|h|)|\,dh
    \leq C\int_{t/4}^{2r_0}
    R^{\alpha-s-1}\sqrt{L(R)}\,dR.
\]
Combining these estimates with Step~2 and the boundedness
of $H_s$, we conclude that
\[
    |G_s(x)|\leq C\left(
    1+t^{\alpha-s}\sqrt{L(t)}
    +\int_{t/4}^{2r_0}
    R^{\alpha-s-1}\sqrt{L(R)}\,dR
    \right).
\]

If $\alpha<s$, then
\[
 \int_{t/4}^{2r_0}
    R^{\alpha-s-1}\sqrt{L(R)}\,dR
    \leq \sqrt{L(t/4)}
    \int_{t/4}^{2r_0}R^{\alpha-s-1}\,dR\leq Ct^{\alpha-s}\sqrt{L(t)}.
\]
If $\alpha=s$, then
\[
    \int_{t/4}^{2r_0}\frac{\sqrt{L(R)}}{R}\,dR
    =\frac23\bigl(L(t/4)^{3/2}
    -L(2r_0)^{3/2}\bigr)
    \leq CL(t)^{3/2}.
\]
Finally, if $\alpha>s$, then
$R^{\alpha-s-1}\sqrt{L(R)}$ is integrable at zero,
and $t^{\alpha-s}\sqrt{L(t)}$ is bounded.
These three estimates prove~\eqref{eq:lift-image-bound}.
\end{proof}

At $\alpha=s$, Lemma~\ref{lem:distance-lift} gives
$|G_s(x)|\leq CL(\delta(x))^{3/2}$.
To control the contribution of $G_s$ against zero-mean
test functions in the critical blow-up argument,
we estimate its local $L^2$ oscillation after
subtracting a suitable constant.

\begin{lemma}
\label{lem:critical-image}
Suppose that $\alpha=s$. After decreasing $r_0$ if
necessary, for every $z\in\partial\Omega$ and $0<r<r_0$,
there exists $b_{z,r}\in\mathbb R$, with
$|b_{z,r}|\leq CL(r)^{3/2}$, such that, for every
$x\in\Omega\cap B_r(z)$,
\begin{equation}
    |G_s(x)-b_{z,r}|
    \leq C\sqrt{L(r)}
    \left(1+\ln\frac{r}{\delta(x)}\right)^{3/2}.
    \label{eq:critical-image-compensation}
\end{equation}
In particular,
\begin{equation}
    \left(r^{-n}\int_{\Omega\cap B_r(z)}
    |G_s-b_{z,r}|^2\,dx\right)^{1/2}
    \leq C\sqrt{L(r)}.
    \label{eq:critical-image-l2}
\end{equation}
\end{lemma}

\begin{proof}[\textnormal{\textbf{Proof}}]
If $n=1$, Lemma~\ref{lem:distance-lift} gives
$G_s\in L^\infty(\Omega)$, so both assertions hold
with $b_{z,r}=0$. We henceforth assume $n\geq2$.

Write $K(h)=K_{s+\ln}(|h|)$ for $h\neq0$.
Fix $\rho>0$ sufficiently small that
\eqref{eq:geometry-annulus} holds for $0<R<4\rho$
and $12\rho<r_\Phi$.
Since $\alpha=s$, we have, for every $p\in\partial\Omega$
and $0<R<4\rho$,
\[
    \int_{B_{2R}(p)}|E_p(y)|\,dy
    \leq C\frac{R^{n+2s}}{\sqrt{L(R)}}.
\]
Moreover, the expression for the kernel gives
\[
    |K(h)|+|h||\nabla K(h)|
    \leq C|h|^{-n-2s}L(|h|)
    \qquad\text{for }0<|h|<4\rho.
\]
Fix $\kappa=8$ and decrease $r_0$ so that
$4\kappa r_0<\rho$ and $r_0<b/2$.

For $p\in\partial\Omega$ and $0<r<r_0$, define
\[
    b_{p,r}:=-\int_{|y-p|>\kappa r}E_p(y)K(y-p)\,dy.
\]
The integral is absolutely convergent since $E_p$
is bounded and the kernel is integrable away from zero.
For $\kappa r\leq R<\rho$, the preceding estimates give
\[
    \int_{R<|y-p|<2R}|E_p(y)K(y-p)|\,dy
    \leq C\sqrt{L(R)}.
\]
Summing over the dyadic radii $R=2^j\kappa r<\rho$
and using the uniform boundedness of the contribution
from $|y-p|\geq\rho$, we obtain
\[
    |b_{p,r}|
    \leq C\left(1+
    \sum_{2^j\kappa r<\rho}\sqrt{L(2^j\kappa r)}
    \right)
    \leq CL(r)^{3/2}.
\]

\medskip
\noindent\textbf{Step 1: Compensation at a nearest boundary point.}
Fix $z\in\partial\Omega$, $0<r<r_0$, and
$x\in\Omega\cap B_r(z)$.
Set $t=\delta(x)$ and choose a nearest boundary point $w$.
Then $t\leq r$, $|w-z|\leq2r$ and
$x=w+t\nu_{\rm in}(w)$.
By~\eqref{eq:tangential-reduction}
and~\eqref{eq:factor-flat-image},
\[
    G_s(x)=H_s(t)+(-\Delta)^{s+\ln}E_w(x).
\]
We first estimate $G_s(x)-b_{w,r}$.

By Step~2 in the proof of Lemma~\ref{lem:distance-lift},
the principal-value integral over $|y-x|<t/4$
is bounded in absolute value by $C\sqrt{L(t)}$.
Moreover, by~\eqref{eq:flat-error-local},
\[
    |E_w(x)|\int_{|y-x|\geq t/4}|K(y-x)|\,dy
    \leq C\sqrt{L(t)}.
\]
It remains to estimate the terms containing $E_w(y)$.

On the region $|y-w|\leq\kappa r$ and
$|y-x|\geq t/4$, we have $|y-x|\leq(\kappa+1)r$.
For $t/4\leq R\leq(\kappa+1)r$, the inclusion
$B_{2R}(x)\subset B_{6R}(w)$ and
\eqref{eq:geometry-annulus} give
\[
    \begin{aligned}
    \int_{R<|y-x|<2R}|E_w(y)K(y-x)|\,dy
    &\leq CR^{-n-2s}L(R)
    \frac{R^{n+2s}}{\sqrt{L(R)}}\leq C\sqrt{L(R)}.
    \end{aligned}
\]
There are at most $C(1+\ln(r/t))$ dyadic annuli
with radii $R=2^jt/4\leq(\kappa+1)r$.
Since $\sqrt{L(R)}\leq C\sqrt{L(t)}$ on these scales,
we obtain
\[
    \int_{\substack{|y-w|\leq\kappa r\\|y-x|\geq t/4}}
    |E_w(y)K(y-x)|\,dy
    \leq C\sqrt{L(t)}
    \left(1+\ln\frac rt\right).
\]

On the remaining region $|y-w|>\kappa r$,
subtracting $b_{w,r}$ leaves the integral
\[
    -\int_{|y-w|>\kappa r}
    E_w(y)\bigl(K(y-x)-K(y-w)\bigr)\,dy.
\]
Since $|x-w|=t\leq r$, the mean value theorem and
the kernel derivative estimate give, for
$R<|y-w|<2R$ and $\kappa r\leq R<\rho$,
\[
    |K(y-x)-K(y-w)|
    \leq CtR^{-n-2s-1}L(R).
\]
Thus the contribution of each annulus is bounded by
\[
    CtR^{-n-2s-1}L(R)
    \int_{B_{2R}(w)}|E_w(y)|\,dy
    \leq CtR^{-1}\sqrt{L(R)}.
\]
Summing over the dyadic radii $R=2^j\kappa r<\rho$,
we obtain
\[
    Ct\sum_{2^j\kappa r<\rho}
    \frac{\sqrt{L(2^j\kappa r)}}{2^j\kappa r}
    \leq C\frac tr\sqrt{L(r)}
    \leq C\sqrt{L(r)}.
\]
For $|y-w|\geq\rho$, we have $|y-x|\geq\rho/2$.
Since $E_w$ is uniformly bounded and $K$ is absolutely
integrable away from zero, the remaining contribution
is uniformly bounded.

Combining these estimates with the boundedness of $H_s$
gives
\[
    |G_s(x)-b_{w,r}|
    \leq C\sqrt{L(t)}
    \left(1+\ln\frac rt\right).
\]

\medskip
\noindent\textbf{Step 2: Comparison of $b_{w,r}$ and $b_{z,r}$.}
We show that $|b_{w,r}-b_{z,r}|\leq C\sqrt{L(r)}$.
Since $|w-z|\leq2r$ and $\partial\Omega$ is of class
$C^{1,s}$, the estimates in the proof of
Lemma~\ref{lem:normalized-distance} give
\[
    |\nu_{\rm in}(w)-\nu_{\rm in}(z)|\leq Cr^s,
    \qquad
    |(w-z)\cdot\nu_{\rm in}(w)|\leq Cr^{1+s}.
\]
Consequently, for $\kappa r\leq R<\rho$ and
$y\in B_{2R}(z)$,
\[
    \begin{aligned}
    |\ell_w(y)-\ell_z(y)|
    &\leq |y-z|
    |\nu_{\rm in}(w)-\nu_{\rm in}(z)|
    +|(w-z)\cdot\nu_{\rm in}(w)|\\
    &\leq C(Rr^s+r^{1+s})
    \leq Cr^sR.
    \end{aligned}
\]
Set $\varepsilon=Cr^sR$.
After decreasing $r_0$, we may assume that
$\varepsilon\leq R$.
By monotonicity of $\Phi_s$ near zero,
\[
    |P_w(y)-P_z(y)|
    \leq\Phi_s(\ell_z(y)+\varepsilon)
    -\Phi_s(\ell_z(y)-\varepsilon).
\]
Integrating in the direction $\nu_{\rm in}(z)$
as in Step~1 of Lemma~\ref{lem:distance-lift}
and using~\eqref{eq:factor-potter}, we obtain
\[
  \int_{B_{2R}(z)}|P_w-P_z|\,dy
    \leq CR^{n-1}\varepsilon\Phi_s(R)\leq Cr^sR^n\Phi_s(R)
    \leq C\frac{r^sR^{n+s}}{\sqrt{L(R)}}.
\]

To compare the integrals defining $b_{w,r}$ and $b_{z,r}$,
first consider their common integration region
\[
    \mathcal F=\{y:|y-w|>\kappa r,\ |y-z|>\kappa r\}.
\]
On $\mathcal F$, we write
\[
    \begin{aligned}
    &E_w(y)K(y-w)-E_z(y)K(y-z)\\
    &\qquad=(P_z(y)-P_w(y))K(y-z)
    +E_w(y)\bigl(K(y-w)-K(y-z)\bigr).
    \end{aligned}
\]
For $\kappa r\leq R<\rho$, the contribution of the
first term on $\mathcal F\cap\{R<|y-z|<2R\}$
is bounded by
\[
    CR^{-n-2s}L(R)
    \frac{r^sR^{n+s}}{\sqrt{L(R)}}
    =C\left(\frac rR\right)^s\sqrt{L(R)}.
\]

For the second term, the mean value theorem gives
\[
    |K(y-w)-K(y-z)|
    \leq CrR^{-n-2s-1}L(R)
\]
on the same annulus.
Since $B_{2R}(z)\subset B_{4R}(w)$,
\eqref{eq:geometry-annulus} yields
\[
    \begin{aligned}
    &\int_{\mathcal F\cap\{R<|y-z|<2R\}}
    |E_w(y)||K(y-w)-K(y-z)|\,dy\\
    &\qquad\leq CrR^{-n-2s-1}L(R)
    \frac{R^{n+2s}}{\sqrt{L(R)}}
    =C\frac rR\sqrt{L(R)}.
    \end{aligned}
\]
Summing over $R=2^j\kappa r<\rho$ gives
\[
    C\sqrt{L(r)}
    \sum_{j=0}^\infty(2^{-js}+2^{-j})
    \leq C\sqrt{L(r)}.
\]
On $\mathcal F\cap\{|y-z|\geq\rho\}$, both kernels
are separated from their singularities. The boundedness
of $E_w$ and $E_z$ therefore gives a uniform bound
for the remaining contribution.

Finally, the parts of the two integration regions
not contained in $\mathcal F$ lie in the annuli
$\kappa r<|y-p|\leq(\kappa+2)r$, with $p=w$ or $p=z$.
Since $(\kappa+2)r<2\kappa r$, the annular estimate
for $E_pK$ gives
\[
    \int_{\kappa r<|y-p|\leq(\kappa+2)r}
    |E_p(y)K(y-p)|\,dy
    \leq C\sqrt{L(r)}.
\]
Combining the common and noncommon parts, we obtain
$|b_{w,r}-b_{z,r}|\leq C\sqrt{L(r)}$.

Together with Step~1, this gives
\[
    |G_s(x)-b_{z,r}|
    \leq C\sqrt{L(t)}
    \left(1+\ln\frac rt\right).
\]
Since $L(t)=L(r)+\ln(r/t)$ and $L(r)\geq1$,
we have
\[
    \sqrt{L(t)}
    \leq\sqrt{L(r)}
    \left(1+\ln\frac rt\right)^{1/2}.
\]
This proves~\eqref{eq:critical-image-compensation}.

\medskip
\noindent\textbf{Step 3: The local $L^2$ estimate.}
By the uniform Lipschitz boundary charts, for
$0<\tau\leq r<r_0$,
\[
    |\{x\in\Omega\cap B_r(z):\delta(x)\leq\tau\}|
    \leq Cr^{n-1}\tau.
\]
Indeed, in coordinates in which the boundary near $z$
is given by $x_n=\varphi(x')$, a nearest boundary
point $w$ to such an $x$ belongs to $B_{2r}(z)$.
Consequently,
\[
    0<x_n-\varphi(x')
    \leq |x_n-w_n|+|\varphi(w')-\varphi(x')|
    \leq C\tau.
\]
The tangential projection of $B_r(z)$ has measure at
most $Cr^{n-1}$, proving the measure estimate.

For each integer $k\geq0$, the layer
\[
    \{x\in\Omega\cap B_r(z):
    re^{-k-1}<\delta(x)\leq re^{-k}\}
\]
has measure at most $Cr^ne^{-k}$.
Moreover, on this layer,
$1+\ln(r/\delta(x))\leq k+2$.
Squaring~\eqref{eq:critical-image-compensation}
and summing over these layers yields
\[
    \int_{\Omega\cap B_r(z)}|G_s-b_{z,r}|^2\,dx
    \leq Cr^nL(r)\sum_{k=0}^\infty(k+2)^3e^{-k}
    \leq Cr^nL(r).
\]
This proves~\eqref{eq:critical-image-l2}.
\end{proof}

To use $\Psi_s$ in the weak formulation, we verify that
it belongs to the Dirichlet energy space and identify
its quadratic form pairing with the pointwise operator.
The zero-exterior Hardy inequality allows us to treat
the possibly unbounded function $G_s$.

\begin{lemma}
\label{lem:factor-energy}
The function $\Psi_s$ defined in~\eqref{eq:regularized-factor}
belongs to $\mathcal H_0^{s+\ln}(\Omega)$ and satisfies
$\delta^sG_s\in L^2(\Omega)$. Moreover, for every
$\varphi\in\mathcal H_0^{s+\ln}(\Omega)$,
\begin{equation}
    \mathcal E_{s+\ln}(\Psi_s,\varphi)
    =\int_\Omega G_s\varphi\,dx,
    \label{eq:factor-lift-weak}
\end{equation}
where the integral on the right-hand side is absolutely
convergent.
\end{lemma}

\begin{proof}[\textnormal{\textbf{Proof}}]
\noindent\textbf{Step 1: Finite energy.}
Since $D$ is locally Lipschitz in $\mathcal U$,
$\zeta\in C_c^\infty(\mathcal U)$ and
$\Phi_s\in C^s(\mathbb R)$, the zero extension of
$\Psi_s=\zeta\Phi_s(D)$ belongs to $C^s(\mathbb R^n)$.
Moreover, $\Psi_s$ is smooth in $\Omega$ and vanishes
in $\Omega^c$. Fix a sufficiently small $t_0\in(0,1/2)$.
For $x\in\Omega$ with $t=\delta(x)<t_0$, we have
$\zeta=1$ and $D\asymp t$ in $B_{t/2}(x)$.
By~\eqref{eq:normalized-distance}, the derivative
estimates for $D$, and~\eqref{eq:factor-derivatives},
the chain rule gives
\[
    \sup_{B_{t/2}(x)}|\nabla\Psi_s|
    +t\sup_{B_{t/2}(x)}|\nabla^2\Psi_s|
    \leq Ct^{s-1}.
\]
Thus $|\Psi_s(x+h)-\Psi_s(x)|\leq Ct^{s-1}|h|$
for $|h|<t/2$.
For $t/2\leq|h|<1$, we use the global $C^s$ estimate.
It follows that
\[
    \begin{aligned}
    \int_{|h|<1}
    \frac{|\Psi_s(x+h)-\Psi_s(x)|^2}{|h|^{n+2s}}
    \ln\frac1{|h|}\,dh&\leq
    Ct^{2s-2}\int_0^{t/2}r^{1-2s}\ln\frac1r\,dr
    +C\int_{t/2}^1\frac{\ln(1/r)}r\,dr\\
    &\leq C\bigl(L(t)+L(t)^2\bigr)
    \leq CL(t)^2.
    \end{aligned}
\]
Here the first integral is bounded by
$Ct^{2-2s}L(t)$ since $s<1$, while the second equals
$\frac12(\ln(2/t))^2$.
For $\delta(x)\geq t_0$, the same energy density
is uniformly bounded by the interior smoothness and
global boundedness of $\Psi_s$.

Since $\Omega$ is bounded and Lipschitz, we have
$|\{x\in\Omega:\delta(x)\leq\tau\}|\leq C\tau$
for $0<\tau\leq t_0$.
Summing over the layers
$t_0e^{-k-1}<\delta(x)\leq t_0e^{-k}$ gives
\[
    \int_{\{x\in\Omega:\delta(x)<t_0\}}
    L(\delta(x))^2\,dx
    \leq Ct_0\sum_{k=0}^\infty
    (L(t_0)+k+1)^2e^{-k}<\infty.
\]
As $\Psi_s=0$ in $\Omega^c$, its full logarithmic
energy is at most twice the integral of the preceding
energy density over $\Omega$.
Together with $\Psi_s\in L^2(\mathbb R^n)$,
this proves $\Psi_s\in\mathcal H_0^{s+\ln}(\Omega)$.

\medskip
\noindent\textbf{Step 2: Weighted convergence of the
truncated operators.}
For $\varepsilon>0$, define
\[
    G_{s,\varepsilon}(x)
    :=\int_{|h|>\varepsilon}
    \bigl(\Psi_s(x)-\Psi_s(x+h)\bigr)
    K_{s+\ln}(|h|)\,dh.
\]
By Lemma~\ref{lem:distance-lift},
$G_{s,\varepsilon}(x)\to G_s(x)$ for every $x\in\Omega$
as $\varepsilon\downarrow0$.

For $t=\delta(x)<t_0$, the Hessian estimate in Step~1
gives
\[
    |2\Psi_s(x)-\Psi_s(x+h)-\Psi_s(x-h)|
    \leq Ct^{s-2}|h|^2
    \qquad\text{for }|h|<t/2.
\]
For $t/2\leq|h|<1$, the same second difference
is bounded by $C|h|^s$.
Using the symmetry of the kernel and its absolute
integrability on $\{|h|\geq1\}$, we obtain
\[
    |G_{s,\varepsilon}(x)|
    \leq C\left(
    1+t^{s-2}\int_0^{t/2}r^{1-2s}L(r)\,dr
    +\int_{t/2}^1r^{-1-s}L(r)\,dr
    \right).
\]
Since $0<s<1$, we have
\[
    \int_0^{t/2}r^{1-2s}L(r)\,dr
    \leq Ct^{2-2s}L(t),
    \qquad
    \int_{t/2}^1r^{-1-s}L(r)\,dr
    \leq Ct^{-s}L(t).
\]
Consequently,
$|G_{s,\varepsilon}(x)|\leq Ct^{-s}L(t)$
uniformly in $\varepsilon$.
For $\delta(x)\geq t_0$, the truncations are uniformly
bounded by the same argument with a fixed radius $t_0/2$.

The integrability of $L(\delta)^2$ in the boundary
collar, established in Step~1, now gives
$\delta^sG_s\in L^2(\Omega)$ and, by dominated convergence,
\[
    \lim_{\varepsilon\downarrow0}
    \|\delta^s(G_{s,\varepsilon}-G_s)\|_{L^2(\Omega)}
    =0.
\]

\medskip
\noindent\textbf{Step 3: Identification of the quadratic
form pairing.}
Let $\varphi\in\mathcal H_0^{s+\ln}(\Omega)$.
By the zero-exterior Hardy inequality
\eqref{eq:zero-exterior-hardy} and the continuous embedding
$\mathcal H^{s+\ln}(\mathbb R^n)\hookrightarrow H^s(\mathbb R^n)$,
we have
\[
    \|\delta^{-s}\varphi\|_{L^2(\Omega)}
    \leq C\|\varphi\|_{H^s(\mathbb R^n)}
    \leq C\|\varphi\|_{\mathcal H^{s+\ln}}.
\]
Since $\delta^sG_s\in L^2(\Omega)$, it follows that
$G_s\varphi\in L^1(\Omega)$.
Moreover, Step~2 gives
\[
    \begin{aligned}
    \int_\Omega|(G_{s,\varepsilon}-G_s)\varphi|\,dx
    &\leq
    \|\delta^s(G_{s,\varepsilon}-G_s)\|_{L^2(\Omega)}
    \|\delta^{-s}\varphi\|_{L^2(\Omega)}\rightarrow0
    \qquad\text{as }\varepsilon\downarrow0.
    \end{aligned}
\]

For each $\varepsilon>0$, the truncated kernel belongs
to $L^1(\mathbb R^n)$. Since
$\Psi_s,\varphi\in L^2(\mathbb R^n)$,
Fubini's theorem and the symmetry of the kernel yield
\[
 \int_\Omega G_{s,\varepsilon}\varphi\,dx=\frac12\iint_{|x-y|>\varepsilon}
    (\Psi_s(x)-\Psi_s(y))(\varphi(x)-\varphi(y))
    K_{s+\ln}(|x-y|)\,dx\,dy.
\]
By Step~1 and the quadratic form estimates in
Section~\ref{sec:forms}, the full integrand on the
right-hand side is absolutely integrable.
Dominated convergence therefore gives
\eqref{eq:factor-lift-weak}.
\end{proof}

\section{Boundary quotients and expansions}
\label{sec:boundary-quotient}

We establish a localization estimate modulo constants
and a half-space Liouville theorem allowing growth of order $2s$.
Together with the local compactness estimates
in Section~\ref{sec:forms}, these results bound the $L^2$ error
in approximating $u$ by multiples of $\Psi_s$.
The approximation coefficients converge to the boundary trace.
Applying the interior H\"older estimate at the boundary-distance
scale yields the pointwise expansion and modulus of continuity
of $u/\Phi_s(\delta)$.

\subsection{Localization and a half-space Liouville theorem}

At the critical growth order $2s$, the weighted tail
in~\eqref{eq:tail-weight} may diverge.
Subtracting the kernel at the origin gives one additional
power of decay at infinity and allows us to localize
the equation up to a constant. This corresponds to the
case of constant polynomials in the definition of the
fractional Laplacian modulo polynomials
\cite[Definition~1.1]{DipierroSavinValdinoci2019};
see also~\cite[Definition~3.1 and Lemma~3.6]
{AbatangeloRosOton2020}.
We establish a localization estimate uniform in
$\ell\in[\ell_0,\infty]$.

Set
\[
    \widetilde\omega_s(y)
    :=\frac{1+\ln(1+|y|)}{(1+|y|)^{n+2s+1}}.
\]
For $R\geq1$, choose $\chi_R\in C_c^\infty(B_{4R})$
with $0\leq\chi_R\leq1$ and $\chi_R=1$ on $B_{3R}$.
We use the convention $A_\infty=(-\Delta)^s$ and
$J_\infty(h)=c_{n,s}|h|^{-n-2s}$.

\begin{lemma}
\label{lem:compensated-localization}
Let $R\geq1$, $\ell\in[\ell_0,\infty]$, and let
$\mathcal O\subset\mathbb R^n$ be open, with
$B_\kappa(x_*)\Subset\mathcal O\cap B_R$.
Suppose that $v\in L^1(\mathbb R^n,\omega_s\,dy)$ and
$v\in\mathcal H_{\rm loc}^{s+\ln}(\mathbb R^n)$
if $\ell<\infty$, or $v\in H^s_{\rm loc}(\mathbb R^n)$
if $\ell=\infty$. Assume that
\[
    A_\ell v=a+g
    \qquad\text{in }\mathcal O\cap B_{2R}
\]
in the sense of distributions, where $a\in\mathbb R$
and $g\in L^1_{\rm loc}(\mathcal O\cap B_{2R})$.
Set
\[
    \mathcal M:=\|v\|_{L^2(B_{4R})}
    +\int_{\mathbb R^n}|v|\widetilde\omega_s\,dy
    +\|g\|_{L^1(B_\kappa(x_*))}.
\]
Then $w=\chi_Rv$ belongs to
$\mathcal H^{s+\ln}(\mathbb R^n)$ if $\ell<\infty$,
and to $H^s(\mathbb R^n)$ if $\ell=\infty$.
Moreover, there exist $b\in\mathbb R$ and
$F\in C^\infty(B_{2R})$ such that
\begin{equation}
    A_\ell w=g+b+F
    \qquad\text{in }\mathcal O\cap B_{2R}.
    \label{eq:localized-equation}
\end{equation}
For every integer $k\geq0$,
\[
    |b|+\|F\|_{C^k(B_{2R})}\leq C_k\mathcal M,
\]
where $C_k$ depends on $n,s,R,\kappa,x_*$ and the
cutoffs, but is independent of $\ell$ and $a$.

If the original equation holds for compactly supported
zero-exterior energy tests in $B_{2R}$, then so does
\eqref{eq:localized-equation}. Moreover, if $v=0$ a.e.
in $\mathcal O^c$, then $w=0$ a.e. in $\mathcal O^c$.
\end{lemma}

\begin{proof}[\textnormal{\textbf{Proof}}]
\noindent\textbf{Step 1: Localization of the equation.}
By the local energy assumption and
$\chi_R\in C_c^\infty(B_{4R})$, the function $w=\chi_Rv$
belongs to the corresponding global energy space.
Define
\[
    b:=a+\int_{\mathbb R^n}
    (1-\chi_R(y))v(y)J_\ell(-y)\,dy
\]
and, for $x\in B_{2R}$,
\[
    F(x):=\int_{\mathbb R^n}
    (1-\chi_R(y))v(y)
    \bigl(J_\ell(x-y)-J_\ell(-y)\bigr)\,dy.
\]
The integral defining $b$ is absolutely convergent
since $v\in L^1(\mathbb R^n,\omega_s\,dy)$ and
$1-\chi_R=0$ on $B_{3R}$.

From the expression for $J_\ell$, for every multi-index
$\gamma$, we have
\[
    |\partial^\gamma J_\ell(h)|
    \leq C_\gamma
    |h|^{-n-2s-|\gamma|}(1+|\ln|h||),
    \qquad h\neq0,
\]
uniformly for $\ell\in[\ell_0,\infty]$.
For $x\in B_{2R}$, $|y|\geq3R$ and $0\leq\theta\leq1$,
we have $|y-\theta x|\asymp|y|$.
The mean value theorem therefore gives
\[
    |J_\ell(x-y)-J_\ell(-y)|
    \leq C_R\widetilde\omega_s(y).
\]
Similarly, for $|\gamma|\geq1$, $|\partial_x^\gamma J_\ell(x-y)|
    \leq C_{R,\gamma}\widetilde\omega_s(y).$
Consequently, differentiation under the integral sign
gives $F\in C^\infty(B_{2R})$ and, for every integer
$k\geq0$,
\[
    \|F\|_{C^k(B_{2R})}
    \leq C_k\int_{\mathbb R^n}|v|\widetilde\omega_s\,dy.
\]

Since $w-v=-(1-\chi_R)v$ vanishes on $B_{3R}$,
for every $\varphi\in C_c^\infty(\mathcal O\cap B_{2R})$,
the separated-support integral gives
\[
\langle A_\ell(w-v),\varphi\rangle
    =\int_{B_{2R}}\varphi(x)
    \int_{\mathbb R^n}(1-\chi_R(y))v(y)
    J_\ell(x-y)\,dy\,dx=\int_{B_{2R}}(b-a+F)\varphi\,dx.
\]
Together with the equation for $v$, this proves
\eqref{eq:localized-equation}.

The inner integral in the first equality is bounded
on $B_{2R}$ by
$C_R\|v\|_{L^1(\mathbb R^n,\omega_s\,dy)}$.
Hence the calculation is valid for compactly
supported zero-exterior energy tests in $B_{2R}$,
whenever the original equation holds for such tests.
Multiplication by $\chi_R$ preserves the
zero-exterior condition.

\medskip
\noindent\textbf{Step 2: Estimate of the constant $b$.}
Choose $\varphi_0\in C_c^\infty(B_\kappa(x_*))$
with $\int_{\mathbb R^n}\varphi_0\,dx=1$.
Testing~\eqref{eq:localized-equation} with $\varphi_0$
and using the symmetry of the kernel, we obtain
\[
    b=\int_{\mathbb R^n}wA_\ell\varphi_0\,dx
    -\int_{B_\kappa(x_*)}(g+F)\varphi_0\,dx.
\]
Indeed, the symmetric second difference of $\varphi_0$
is bounded by $C|h|^2$ for $|h|<1$ and by $C$ for
$|h|\geq1$. Thus
\[
    \|A_\ell\varphi_0\|_{L^\infty(\mathbb R^n)}
    \leq C\left(
    \int_0^1r^{1-2s}(1+|\ln r|)\,dr
    +\int_1^\infty r^{-1-2s}(1+\ln r)\,dr
    \right)
    \leq C,
\]
uniformly for $\ell\in[\ell_0,\infty]$.
The same bound holds for the truncated operators.
Since $w$ is compactly supported and belongs to $L^2$,
Fubini's theorem for the truncated integrals and
dominated convergence justify the preceding pairing.

Using $\operatorname{supp}w\subset B_{4R}$ and
$|w|\leq|v|$, we obtain
\[
    |b|\leq C\left(
    \|v\|_{L^2(B_{4R})}
    +\|g\|_{L^1(B_\kappa(x_*))}
    +\|F\|_{L^\infty(B_{2R})}
    \right).
\]
Together with the estimate for $F$ in Step~1, this gives
$|b|+\|F\|_{C^k(B_{2R})}\leq C_k\mathcal M$.
\end{proof}

The preceding lemma bounds the sum of the original
constant and the constant introduced by truncation,
without estimating them separately.
Since $w$ is supported in $B_{4R}$, its weighted tail
is bounded by $C_R\|v\|_{L^2(B_{4R})}$.
Thus the local estimates in Section~\ref{sec:forms}
apply whenever the localized right-hand side satisfies
their hypotheses.

The half-space limit of the rescaled equation may
have growth of order $2s$ at infinity.
We therefore formulate the limit equation using
zero-mean test functions; see also
\cite[Definition~5.7]{DipierroRosOtonValdinociWeidner2025}
for a related formulation of the Neumann problem.

For every $\varphi\in C_c^\infty(\mathbb R^n)$ with
$\int_{\mathbb R^n}\varphi\,dx=0$, we have
\begin{equation}
    |(-\Delta)^s\varphi(y)|
    \leq C_\varphi(1+|y|)^{-n-2s-1}.
    \label{eq:zero-mean-test-tail}
\end{equation}
The corresponding estimate for the fractional-logarithmic operator is
\begin{equation}
    |(-\Delta)^{s+\ln}\varphi(y)|
    \leq C_\varphi\widetilde\omega_s(y).
    \label{eq:zero-mean-log-test-tail}
\end{equation}
Indeed, for sufficiently large $|y|$, the zero-mean
condition gives
\[
    (-\Delta)^s\varphi(y)
    =-c_{n,s}\int_{\mathbb R^n}\varphi(x)
    \left(
    |y-x|^{-n-2s}-|y|^{-n-2s}
    \right)\,dx.
\]
The mean value theorem gives
\eqref{eq:zero-mean-test-tail}, while the same argument
with $K_{s+\ln}$ and its first derivative gives
\eqref{eq:zero-mean-log-test-tail}.
On bounded sets, both estimates follow from
\eqref{eq:test-function-tail}.

We next adapt the half-space classification of
Abatangelo and Ros-Oton
\cite[Theorem~3.10 and Proposition~3.11]{AbatangeloRosOton2020}
to the $L^2$ growth condition and the distributional
formulation used in our blow-up argument.

\begin{corollary}
\label{lem:critical-liouville}
Let $H=\{x_n>0\}$ and
$W\in H^s_{\rm loc}(\mathbb R^n)$, with $W=0$ a.e.
in $H^c$. Suppose that
\begin{equation}
    \|W\|_{L^2(B_T)}\leq C_0T^{n/2+2s}
    \qquad\text{for every }T\geq1,
    \label{eq:critical-growth}
\end{equation}
and that
\begin{equation}
    \int_{\mathbb R^n}W(-\Delta)^s\varphi\,dx=0
    \label{eq:zero-mean-equation}
\end{equation}
for every $\varphi\in C_c^\infty(H)$ with
$\int_{\mathbb R^n}\varphi\,dx=0$.
Then $W=\kappa(x_n)_+^s$ a.e. in $\mathbb R^n$
for some $\kappa\in\mathbb R$.
\end{corollary}

\begin{proof}[\textnormal{\textbf{Proof}}]
By Cauchy--Schwarz and~\eqref{eq:critical-growth},
for every $R\geq1$,
\[
    \int_{B_{2R}\setminus B_R}|W(y)|\,dy
    \leq CC_0R^{n+2s}.
\]
Consequently,
\[
    \int_{B_{2^{j+1}}\setminus B_{2^j}}
    \frac{|W(y)|}{(1+|y|)^{n+2s+1}}\,dy
    \leq CC_0\,2^{-j}
    \qquad\text{for }j\geq0.
\]
Together with~\eqref{eq:zero-mean-test-tail},
this proves the absolute convergence of the integral
in~\eqref{eq:zero-mean-equation}.

\medskip
\noindent\textbf{Step 1: Localization of the rescaled equation.}
For $T\geq1$, set $W_T(x)=T^{-2s}W(Tx)$.
By the homogeneity of $(-\Delta)^s$,
$W_T$ satisfies~\eqref{eq:zero-mean-equation}.
Moreover, by~\eqref{eq:critical-growth},
\[
    \|W_T\|_{L^2(B_R)}
    =T^{-n/2-2s}\|W\|_{L^2(B_{TR})}
    \leq C_0R^{n/2+2s}
    \qquad\text{for }R\geq1.
\]
The preceding annular estimate applied to $W_T$ also gives
\[
    \int_{\mathbb R^n}|W_T|\widetilde\omega_s\,dy
    \leq CC_0\left(
    1+\sum_{j=0}^\infty(1+j)2^{-j}
    \right)
    \leq CC_0.
\]

Set $w_T=\chi_4W_T$ and, for $x\in B_8$, define
\[
    F_T(x):=\int_{\mathbb R^n}
    (1-\chi_4(y))W_T(y)
    \bigl(J_\infty(x-y)-J_\infty(-y)\bigr)\,dy.
\]
Since $\chi_4=1$ on $B_{12}$, the kernel difference
is bounded by $C(1+|y|)^{-n-2s-1}$ on the integration
region. Hence $\|F_T\|_{L^\infty(B_8)}\leq CC_0$.

The function $w_T$ belongs to $H^s(\mathbb R^n)$
and vanishes a.e. in $H^c$.
For every $\varphi\in C_c^\infty(H\cap B_8)$ with
$\int_{\mathbb R^n}\varphi\,dx=0$,
we obtain from~\eqref{eq:zero-mean-equation}
\[
  \int_{\mathbb R^n}w_T(-\Delta)^s\varphi\,dx
    =-\int_{\mathbb R^n}
    (1-\chi_4)W_T(-\Delta)^s\varphi\,dx=\int_{H\cap B_8}F_T\varphi\,dx.
\]

Choose a fixed $\varphi_0\in C_c^\infty(H\cap B_1)$
with $\int_{\mathbb R^n}\varphi_0\,dx=1$.
For an arbitrary $\varphi\in C_c^\infty(H\cap B_8)$,
applying the preceding identity to
$\varphi-(\int\varphi\,dx)\varphi_0$ gives
\[
    (-\Delta)^sw_T=b_T+F_T
    \qquad\text{in }H\cap B_8
\]
in the sense of distributions, where
\[
    b_T:=\int_{\mathbb R^n}w_T(-\Delta)^s\varphi_0\,dx
    -\int_{H\cap B_8}F_T\varphi_0\,dx.
\]
Since $\operatorname{supp}w_T\subset B_{16}$ and
$(-\Delta)^s\varphi_0$ is bounded, we have
\[
    |b_T|
    \leq C\left(
    \|W_T\|_{L^2(B_{16})}
    +\|F_T\|_{L^\infty(B_8)}
    \right)
    \leq CC_0.
\]

\medskip
\noindent\textbf{Step 2: Pointwise growth and continuity.}
The equation $(-\Delta)^sw_T=b_T+F_T$ also holds
for compactly supported $H^s$ tests in $B_8$
vanishing a.e. in $H^c$.
Indeed, translating such a test a distance
$\varepsilon>0$ in the direction $e_n$ and mollifying
with radius $\varepsilon/2$ gives functions in
$C_c^\infty(H\cap B_8)$ for sufficiently small
$\varepsilon$. These functions converge to the original
test in $H^s(\mathbb R^n)$.
The identity extends by continuity of the fractional
quadratic form and the boundedness of $b_T+F_T$.

Since $\operatorname{supp}w_T\subset B_{16}$,
the ordinary fractional tail satisfies
\[
    \int_{\mathbb R^n}
    \frac{|w_T(y)|}{(1+|y|)^{n+2s}}\,dy
    \leq C\|W_T\|_{L^2(B_{16})}
    \leq CC_0.
\]
Applying the fractional version of
Lemma~\ref{lem:local-boundedness}, we obtain
\[
    \|W_T\|_{L^\infty(B_1)}
    =\|w_T\|_{L^\infty(B_1)}
    \leq CC_0.
\]
Since $T\geq1$ is arbitrary, this yields
$\|W\|_{L^\infty(B_T)}\leq CC_0T^{2s}$ for every
$T\geq1$. Moreover, using
$W_T(x)=16^{2s}W_{16T}(x/16)$, we obtain
\[
    \|w_T\|_{L^\infty(\mathbb R^n)}
    \leq\|W_T\|_{L^\infty(B_{16})}
    \leq CC_0.
\]

The half-ball estimate of Ros-Oton and Serra
\cite[Proposition~4.6]{RosOtonSerra2016Stable},
applied to the equation for $w_T$ in $H\cap B_8$,
gives
\[
    \|W_T\|_{C^s(B_4)}
    =\|w_T\|_{C^s(B_4)}
    \leq C\left(
    \|w_T\|_{L^\infty(\mathbb R^n)}
    +|b_T|+\|F_T\|_{L^\infty(B_8)}
    \right)
    \leq CC_0.
\]
Thus $W$ has a continuous representative satisfying
$W=0$ in $H^c$ and
\[
    |W(x)|\leq CC_0(1+|x|)^{2s}
    \qquad\text{for every }x\in\mathbb R^n.
\]

\medskip
\noindent\textbf{Step 3: Half-space classification.}
Fix $M\geq1$. For $R>2M$, define
\[
    F_R(x):=\int_{|y|>R}
    W(y)\bigl(J_\infty(x-y)-J_\infty(-y)\bigr)\,dy,
    \qquad x\in B_M.
\]
As in Step~1, testing against zero-mean functions
and then subtracting a fixed test function of integral
one gives
\[
    (-\Delta)^s(\mathbf1_{B_R}W)=b_R+F_R
    \qquad\text{in }H\cap B_M
\]
in the sense of distributions, for some constant $b_R$.

Since $R>2M$, the mean value theorem and the annular
estimate at the beginning of the proof give
\[
    \begin{aligned}
    \|F_R\|_{L^\infty(B_M)}
    &\leq C_M\int_{|y|>R}
    \frac{|W(y)|}{(1+|y|)^{n+2s+1}}\,dy\\
    &\leq C_MC_0\sum_{j=0}^\infty(2^jR)^{-1}
    \leq\frac{C_MC_0}{R}.
    \end{aligned}
\]
Therefore $F_R\to0$ uniformly on $B_M$ as $R\to\infty$.
Since $M$ is arbitrary and
\[
    \int_{\mathbb R^n}
    \frac{|W(y)|}{1+|y|^{n+2s+1}}\,dy<\infty,
\]
this verifies $(-\Delta)^sW=0$ modulo constants
in $H$, in the sense of
\cite[Definition~3.1]{AbatangeloRosOton2020} with $k=1$.

The fractional Laplacian has a constant positive
angular kernel and hence satisfies the kernel assumptions
of~\cite[Theorem~3.10]{AbatangeloRosOton2020}.
By Step~2, its growth assumption holds with their
parameter $\alpha=s/2$, since $2s<1+s+s/2$.
Applying their Theorem~3.10 if $n\geq2$, and their
Proposition~3.11 if $n=1$, we obtain
\[
    W(x)=(a+b\cdot x)(x_n)_+^s
\]
for some $a\in\mathbb R$ and $b\in\mathbb R^n$.
If $b\neq0$, choose $\xi\in H$ with $b\cdot\xi\neq0$.
Then $|W(T\xi)|\geq cT^{s+1}$ for sufficiently large
$T$, contradicting the growth bound in Step~2 since
$s+1>2s$. Thus $b=0$, and the conclusion follows.
\end{proof}

\subsection{A unified boundary approximation}

Write $\omega(r)=\omega_{n,s,\alpha}(r)$, where
$\omega_{n,s,\alpha}$ is defined in~\eqref{eq:domain-modulus}.
Throughout this subsection, we further decrease $r_0$
from Lemma~\ref{lem:distance-lift} so that
$\zeta=1$ on $\{x:\delta(x)<4r_0\}$, where $\zeta$
is the cutoff in~\eqref{eq:regularized-factor},
$4r_0<r_\Phi$ with $r_\Phi$ as in
Proposition~\ref{prop:factor-asymptotics}, and
$\ell_r=\ln(1/r)\geq\ell_0$ for $0<r\leq r_0$,
where $\ell_0$ is given in Lemma~\ref{lem:local-boundedness}.

By the definition of $\omega$ and
Proposition~\ref{prop:factor-asymptotics}, we have
$\omega(r)\geq C^{-1}\Phi_s(r)$.
Moreover, for $0<\tau\leq1$,
$\omega(\tau r)\leq\tau^{\min\{s,\alpha\}}\omega(r)$.
Consequently,
$\sum_{k=0}^\infty\omega(2^{-k}r)\leq C\omega(r)$,
and the sum of $\omega$ over increasing dyadic radii
is bounded by a constant times its final term.
For every $\varepsilon>0$, we also have
\begin{equation}
    \frac{\omega(Tr)}{\omega(r)}
    \leq C_\varepsilon T^{s+\varepsilon}
    \qquad\text{for }T\geq1,\quad 0<Tr\leq r_0.
    \label{eq:modulus-dilation}
\end{equation}
Indeed, if $n=1$ or $\alpha>s$, then
\[
    \frac{\omega(Tr)}{\omega(r)}
    =T^s\left(1+\frac{\ln T}{L(Tr)}\right)^{1/2}
    \leq C_\varepsilon T^{s+\varepsilon}.
\]
Otherwise, the ratio equals $T^\alpha\leq T^s$.
In particular, for every fixed $T\geq1$,
\[
    \lim_{r\downarrow0}\frac{\omega(Tr)}{\omega(r)}
    \leq T^s.
\]

To compare the operator estimates in Section~\ref{sec:mellin}
with the rescaled equations, define
\begin{equation}
    S(r):=\frac{r^{2s}}{2\ell_r\Phi_s(r)\omega(r)}
    \asymp
    \begin{cases}
        1,&n=1\ \text{or }\alpha>s,\\
        r^{s-\alpha}L(r)^{-1/2},
            &n\geq2,\ \alpha<s,\\
        L(r)^{-1/2},&n\geq2,\ \alpha=s.
    \end{cases}
    \label{eq:source-scaling-factor}
\end{equation}

We approximate $u$ near each boundary point by its local
$L^2$ projection onto the span of $\Psi_s$, following
the distance-factor approximation of Ros-Oton and
Serra~\cite[Proposition~3.2 and Lemma~3.3]
{RosOtonSerra2017Boundary} and Fall
\cite[Section~5, (5.3)--(5.4)]{Fall2019NonlocalSchrodinger}.
As in the proof of Lemma~\ref{lem:interior-holder},
we select a scale at which the normalized approximation
error is nearly maximal and apply a Liouville theorem.
The rescaled errors may now have growth of order $2s$,
so we use compensated tail estimates and the
half-space classification modulo constants.

\begin{proposition}
\label{prop:endpoint-excess}
Let $u$ satisfy the hypotheses of Theorem~\ref{thm:quotient}.
For $z\in\partial\Omega$ and $0<r\leq r_0$, let
$c_{z,r}$ minimize
$\dashint_{B_r(z)}|u-a\Psi_s|^2\,dx$ over $a\in\mathbb R$,
and set
\[
    \mathcal E_z(r):=\frac1{\Phi_s(r)}
    \left(\dashint_{B_r(z)}
    |u-c_{z,r}\Psi_s|^2\,dx\right)^{1/2}.
\]
Then $c_{z,r}$ is uniquely defined and
\begin{equation}
    \mathcal E_z(r)\leq C\mathcal A\omega(r)
    \qquad\text{for }z\in\partial\Omega,\quad 0<r\leq r_0.
    \label{eq:excess-decay}
\end{equation}
\end{proposition}

\begin{proof}[\textnormal{\textbf{Proof}}]
\noindent\textbf{Step 1: Estimates for the approximation coefficients.}
By the uniform boundary graph representation,
for sufficiently small $r_0$, the ball
$B_{r/4}(z+(r/2)\nu_{\rm in}(z))$ is contained in
$\Omega\cap B_r(z)$.
On the concentric ball of radius $r/8$, we have
$\delta\asymp r$.
Since $\Psi_s=\Phi_s(D)$ near $\partial\Omega$,
$D\asymp\delta$, and $\Psi_s=0$ in $\Omega^c$,
the monotonicity and dilation estimates for $\Phi_s$ give
\begin{equation}
    C^{-1}\Phi_s(r)
    \leq\left(\dashint_{B_r(z)}\Psi_s^2\,dx\right)^{1/2}
    \leq C\Phi_s(r).
    \label{eq:factor-ball-norm}
\end{equation}
For the upper bound, we used $\delta\leq r$
on $\Omega\cap B_r(z)$. In particular, the minimizer is unique and is given by
\[
    c_{z,r}
    =\frac{\int_{B_r(z)}u\Psi_s\,dx}
    {\int_{B_r(z)}\Psi_s^2\,dx}.
\]
By the Cauchy--Schwarz inequality,~\eqref{eq:factor-ball-norm}
and the competitor $a=0$, we obtain
\[
    |c_{z,r}|+\mathcal E_z(r)
    \leq C\frac{r^{-n/2}\|u\|_2}{\Phi_s(r)}.
\]

For $r\leq h\leq\min\{2r,r_0\}$, the triangle inequality
on $B_r(z)$ gives
\[
 |c_{z,h}-c_{z,r}|
    \left(\dashint_{B_r(z)}\Psi_s^2\,dx\right)^{1/2}\leq\Phi_s(r)\mathcal E_z(r)
    +(h/r)^{n/2}\Phi_s(h)\mathcal E_z(h).
\]
Using~\eqref{eq:factor-ball-norm}
and~\eqref{eq:factor-potter}, we conclude that
\begin{equation}
    |c_{z,h}-c_{z,r}|
    \leq C\bigl(\mathcal E_z(r)+\mathcal E_z(h)\bigr).
    \label{eq:endpoint-adjacent-coefficients}
\end{equation}

\medskip
\noindent\textbf{Step 2: Selection of scales and normalization.}
By linearity, we may assume that $\mathcal A\leq1$.
Suppose that~\eqref{eq:excess-decay} fails.
Then there exist solutions $u_j$ with right-hand sides $f_j$
satisfying $\|u_j\|_2+\|f_j\|_\infty\leq1$,
boundary points $z_j$, and radii $t_j\downarrow0$ such that
\[
    \sup_{t_j\leq r\leq r_0}
    \frac{\mathcal E_{z_j}(r)}{\omega(r)}
    \longrightarrow\infty,
\]
where $\mathcal E_{z_j}$ and $c_{z_j,r}$ refer to $u_j$.

Choose $r_j\in[t_j,r_0]$ such that
\[
    \Theta_j:=\frac{\mathcal E_{z_j}(r_j)}{\omega(r_j)}
    \geq\frac12\sup_{t_j\leq r\leq r_0}
    \frac{\mathcal E_{z_j}(r)}{\omega(r)}.
\]
By Step~1, $\mathcal E_{z_j}(r)/\omega(r)$ is uniformly
bounded on every fixed interval $[\tau,r_0]$, with $\tau>0$.
Consequently, $r_j\to0$, $\Theta_j\to\infty$ and
\begin{equation}
    \mathcal E_{z_j}(h)\leq2\Theta_j\omega(h)
    \qquad\text{for }r_j\leq h\leq r_0.
    \label{eq:maximal-excess}
\end{equation}

Set $a_j=c_{z_j,r_j}$.
For $r_j\leq h\leq r_0$, choose $m\geq0$ such that
$2^mr_j\leq h<2^{m+1}r_j$.
Summing~\eqref{eq:endpoint-adjacent-coefficients}
over the radii $r_j,2r_j,\ldots,2^mr_j,h$
and using~\eqref{eq:maximal-excess}, we obtain
\[
    |a_j-c_{z_j,h}|
    \leq C\Theta_j
    \left(\sum_{k=0}^m\omega(2^kr_j)+\omega(h)\right)
    \leq C\Theta_j\omega(h).
\]
Together with the coefficient bound in Step~1, this gives
\[
    \frac{|a_j|}{\Theta_j}
    \leq C\omega(h)
    +\frac{Ch^{-n/2}}{\Theta_j\Phi_s(h)}.
\]
For each fixed $h\in(0,r_0]$, this estimate holds for
all sufficiently large $j$.
First letting $j\to\infty$ and then $h\downarrow0$,
we obtain
\begin{equation}
    \frac{1+|a_j|}{\Theta_j}\longrightarrow0.
    \label{eq:coefficient-normalization}
\end{equation}

Choose $Q_j\in O(n)$ such that
$Q_je_n=\nu_{\rm in}(z_j)$, write
$X_j(y)=z_j+r_jQ_jy$, and define
\[
    v_j(y):=
    \frac{(u_j-a_j\Psi_s)(X_j(y))}
    {\Theta_j\Phi_s(r_j)\omega(r_j)}.
\]
With $\Pi_j(y)=\Psi_s(X_j(y))/\Phi_s(r_j)$,
the definition of $\Theta_j$ and the minimizing
property of $a_j$ give
\begin{equation}
    \dashint_{B_1}|v_j|^2\,dy=1,\qquad
    \int_{B_1}v_j\Pi_j\,dy=0.
    \label{eq:boundary-orthogonality}
\end{equation}

\medskip
\noindent\textbf{Step 3: Growth and compensated tail estimates.}
Let $1\leq T\leq r_0/r_j$ and set $h=Tr_j$.
By~\eqref{eq:factor-ball-norm},
\eqref{eq:maximal-excess} and the coefficient difference
estimate in Step~2, we have
\[
 \left(\dashint_{B_h(z_j)}
    |u_j-a_j\Psi_s|^2\,dx\right)^{1/2}\leq\Phi_s(h)\mathcal E_{z_j}(h)
    +|c_{z_j,h}-a_j|
    \left(\dashint_{B_h(z_j)}\Psi_s^2\,dx\right)^{1/2}\leq C\Theta_j\Phi_s(h)\omega(h).
\]
Changing variables, we conclude that
\begin{equation}
    \|v_j\|_{L^2(B_T)}
    \leq CT^{n/2}
    \frac{\Phi_s(Tr_j)\omega(Tr_j)}
    {\Phi_s(r_j)\omega(r_j)}
    \qquad\left(1\leq T\leq\frac{r_0}{r_j}\right).
    \label{eq:boundary-residual-growth}
\end{equation}

Fix $\varepsilon\in(0,1)$.
Applying~\eqref{eq:factor-potter}
and~\eqref{eq:modulus-dilation}, each with exponent
loss $\varepsilon/2$, we obtain
\[
    \|v_j\|_{L^2(B_T)}
    \leq C_\varepsilon T^{n/2+2s+\varepsilon},
    \qquad 1\leq T\leq\frac{r_0}{r_j}.
\]
Consequently, for every $k\geq0$ with $2^k<r_0/r_j$,
the Cauchy--Schwarz inequality gives
\[
    \begin{aligned}
    \int_{(B_{2^{k+1}}\setminus B_{2^k})
    \cap B_{r_0/r_j}}
    |v_j(y)|\widetilde\omega_s(y)\,dy&\leq
    C(1+k)2^{-k(n+2s+1)}
    2^{kn/2}
    \|v_j\|_{L^2(B_{\min\{2^{k+1},r_0/r_j\}})}\\
    \qquad&\leq C_\varepsilon(1+k)
    2^{-k(1-\varepsilon)}.
    \end{aligned}
\]
Since $\varepsilon<1$, the resulting series is summable.

It remains to estimate the tail outside $B_{r_0/r_j}$.
The function $u_j-a_j\Psi_s$ vanishes a.e. in $\Omega^c$.
For $x\in\Omega$ with $|x-z_j|\geq r_0$, we have
$r_0\leq|x-z_j|\leq\operatorname{diam}\Omega$, and hence
\[
    r_j^{-n}\widetilde\omega_s
    \left(r_j^{-1}Q_j^T(x-z_j)\right)
    \leq Cr_j^{2s+1}(1+\ell_{r_j}).
\]
Changing variables gives
\[
 \int_{|y|>r_0/r_j}|v_j(y)|
    \widetilde\omega_s(y)\,dy\leq
    \frac{Cr_j^{2s+1}(1+\ell_{r_j})}
    {\Theta_j\Phi_s(r_j)\omega(r_j)}
    \|u_j-a_j\Psi_s\|_{L^1(\Omega)}.
\]
Since $\Omega$ is bounded, $\|u_j\|_2\leq1$ and
$\Psi_s\in L^1(\Omega)$, we have
$\|u_j-a_j\Psi_s\|_1\leq C(1+|a_j|)$.
Moreover, $\omega(r_j)\geq C^{-1}\Phi_s(r_j)$ and
$\Phi_s(r_j)^2\asymp r_j^{2s}/\ell_{r_j}$.
Thus the preceding tail is bounded by
\[
    C\frac{1+|a_j|}{\Theta_j}
    r_j(1+\ell_{r_j})^2\longrightarrow0,
\]
where we used~\eqref{eq:coefficient-normalization}. Combining the annular estimates with the far-field bound,
we obtain
\begin{equation}
    \lim_{M\to\infty}\limsup_{j\to\infty}
    \int_{|y|>M}|v_j(y)|\widetilde\omega_s(y)\,dy=0.
    \label{eq:boundary-tail-tightness}
\end{equation}
In particular,
$\sup_j\int_{\mathbb R^n}|v_j|\widetilde\omega_s\,dy<\infty$.
Each $v_j$ also belongs to
$L^1(\mathbb R^n,\omega_s\,dy)$, since it is compactly
supported and belongs to $L^2(\mathbb R^n)$.

\medskip
\noindent\textbf{Step 4: Estimates for the rescaled right-hand sides.}
Set $\Omega_j=r_j^{-1}Q_j^T(\Omega-z_j)$ and
$\delta_j=\operatorname{dist}(\cdot,\Omega_j^c)$.
By the uniform $C^{1,\alpha}$ boundary graph representation,
on every fixed bounded set the boundary of $\Omega_j$
is eventually given by a graph whose height and first
derivatives are bounded by $C_Rr_j^\alpha$.
Thus $\Omega_j$ converges locally to $H=\{y_n>0\}$ and
$\delta_j\to(y_n)_+$ uniformly on bounded sets.
In particular, every compact subset of $H$ is contained
in $\Omega_j$ for all sufficiently large $j$.

By Lemma~\ref{lem:factor-energy}
and Lemma~\ref{lem:dilation}, we have
$v_j\in\mathcal H_0^{s+\ln}(\Omega_j)$ and
\[
    A_{\ell_{r_j}}v_j
    =\frac{S(r_j)}{\Theta_j}(f_j-a_jG_s)\circ X_j
    \qquad\text{in }\Omega_j.
\]
The equation holds for compactly supported zero-exterior
energy tests.

Fix $R\geq1$. We show that, for all sufficiently large $j$,
the equation can be written as
\[
    A_{\ell_{r_j}}v_j=d_{j,R}+g_{j,R}
    \qquad\text{in }\Omega_j\cap B_{2R},
\]
where $d_{j,R}\in\mathbb R$ and
\begin{equation}
    \|\delta_j^sg_{j,R}\|_{L^2(\Omega_j\cap B_{2R})}
    \longrightarrow0.
    \label{eq:rescaled-weighted-source}
\end{equation}
Moreover, $g_{j,R}\to0$ locally in
$L^\infty(H\cap B_{2R})$.
Write $\varepsilon_j=(1+|a_j|)/\Theta_j$,
so that $\varepsilon_j\to0$ by
\eqref{eq:coefficient-normalization}.

\smallskip
\noindent\textbf{Case 1: $n=1$ or $\alpha>s$.}
Take $d_{j,R}=0$.
By Lemma~\ref{lem:distance-lift}
and~\eqref{eq:source-scaling-factor},
both $G_s$ and $S(r_j)$ are bounded.
Consequently, $ \|g_{j,R}\|_{L^\infty(\Omega_j\cap B_{2R})}
    \leq C\varepsilon_j.$
Since $0\in\partial\Omega_j$, we have $\delta_j\leq2R$
on $\Omega_j\cap B_{2R}$.
This proves~\eqref{eq:rescaled-weighted-source}
and the local $L^\infty$ convergence.

\smallskip
\noindent\textbf{Case 2: $n\geq2$ and $\alpha<s$.}
Again, take $d_{j,R}=0$.
Since $\delta(X_j(y))=r_j\delta_j(y)$,
\eqref{eq:lift-image-bound}
and~\eqref{eq:source-scaling-factor} give,
a.e. in $\Omega_j\cap B_{2R}$,
\[
    |g_{j,R}|
    \leq C_R\varepsilon_j
    \left(
    1+\delta_j^{\alpha-s}
    \sqrt{\frac{L(r_j\delta_j)}{L(r_j)}}
    \right).
\]
Using $L(r_j\delta_j)=L(r_j)+\ln(1/\delta_j)$
and $0<\delta_j\leq2R$, we obtain
\[
    \delta_j^s|g_{j,R}|
    \leq C_R\varepsilon_j
    \left(
    \delta_j^s+\delta_j^\alpha
    \sqrt{1+|\ln\delta_j|}
    \right)
    \leq C_R\varepsilon_j,
\]
where the last inequality follows from $\alpha>0$.
This proves~\eqref{eq:rescaled-weighted-source}.
On compact subsets of $H\cap B_{2R}$, the distances
$\delta_j$ are bounded away from zero for all sufficiently
large $j$, so $g_{j,R}\to0$ locally in $L^\infty$.

\smallskip
\noindent\textbf{Case 3: $n\geq2$ and $\alpha=s$.}
Let $b_j=b_{z_j,2Rr_j}$ be the constant given by
Lemma~\ref{lem:critical-image}, and set
$d_{j,R}=-S(r_j)a_jb_j/\Theta_j$.
Then
\[
    g_{j,R}
    =\frac{S(r_j)}{\Theta_j}
    \bigl(f_j-a_j(G_s-b_j)\bigr)\circ X_j.
\]
By~\eqref{eq:critical-image-l2} and a change of variables,
\[
    \|(G_s-b_j)\circ X_j\|_{L^2(\Omega_j\cap B_{2R})}
    \leq C_R\sqrt{L(2Rr_j)}.
\]
Since $S(r_j)\asymp L(r_j)^{-1/2}$ and
$L(2Rr_j)\asymp L(r_j)$ for fixed $R$, we obtain
\[
    \|g_{j,R}\|_{L^2(\Omega_j\cap B_{2R})}
    \leq C_R\varepsilon_j\longrightarrow0.
\]
This gives~\eqref{eq:rescaled-weighted-source},
since $\delta_j\leq2R$.
Moreover, by~\eqref{eq:critical-image-compensation},
a.e. in $\Omega_j\cap B_{2R}$,
\[
    |g_{j,R}|
    \leq C_R\varepsilon_j
    \left(1+\ln\frac{2R}{\delta_j}\right)^{3/2}.
\]
Thus $g_{j,R}\to0$ locally in $L^\infty(H\cap B_{2R})$.
No uniform bound for $d_{j,R}$ is required.

\medskip
\noindent\textbf{Step 5: Compactness and the Liouville contradiction.}
Fix $M\geq1$ and set $R=2M$.
We apply Lemma~\ref{lem:compensated-localization}
to the equation in Step~4, using a fixed ball compactly
contained in $H\cap B_R$.
For all sufficiently large $j$, this ball is contained
in $\Omega_j\cap B_R$.
By~\eqref{eq:boundary-residual-growth},
\eqref{eq:boundary-tail-tightness} and the local
$L^\infty$ convergence of $g_{j,R}$, the quantity
$\mathcal M$ in the localization lemma is uniformly bounded.
Therefore, for $w_j=\chi_Rv_j$, we have
\[
    A_{\ell_{r_j}}w_j=g_{j,R}+\widehat b_j+F_j
    \qquad\text{in }\Omega_j\cap B_{2R},
\]
where $|\widehat b_j|+\|F_j\|_{L^\infty(B_{2R})}\leq C_R$.
The equation holds for compactly supported zero-exterior
energy tests in $B_{2R}$. By~\eqref{eq:rescaled-weighted-source}
and $\delta_j\leq2R$, we have
\[
    \|\delta_j^s(g_{j,R}+\widehat b_j+F_j)\|
    _{L^2(\Omega_j\cap B_{2R})}\leq C_R.
\]
Moreover, rescaling~\eqref{eq:zero-exterior-hardy}
gives, for every $\varphi\in H^s(\mathbb R^n)$
vanishing a.e. in $\Omega_j^c$,
\[
    \|\delta_j^{-s}\varphi\|_{L^2(\Omega_j)}^2
    \leq C\left(
    r_j^{2s}\|\varphi\|_2^2
    +[\varphi]_{H^s(\mathbb R^n)}^2\right)\leq C\|\varphi\|_{H^s(\mathbb R^n)}^2.
\]
Hence, for every zero-exterior energy test $\varphi$
compactly supported in $B_{2R}$,
\[
    \left|
    \int_{\Omega_j\cap B_{2R}}
    (g_{j,R}+\widehat b_j+F_j)\varphi\,dy
    \right|
    \leq C_R\|\varphi\|_{H^s(\mathbb R^n)}.
\]

Since $w_j$ is supported in $B_{4R}$, its weighted tail
is uniformly bounded by $\|v_j\|_{L^2(B_{4R})}$.
Thus the local energy estimate in the proof of
Lemma~\ref{lem:local-compactness} applies to $w_j$
with cutoffs supported in $B_{2R}$.
Since $w_j=v_j$ on $B_{3R}$, we obtain
$\sup_j\|v_j\|_{H^s(B_M)}<\infty$.

By fractional Rellich compactness
\cite[Theorems~5.4 and~7.1]{DiNezzaPalatucciValdinoci2012}
and a diagonal argument, passing to a subsequence gives
\[
    v_j\longrightarrow v
    \quad\text{strongly in }L^2_{\rm loc}(\mathbb R^n)
    \text{ and a.e. in }\mathbb R^n,
\]
with $v\in H^s_{\rm loc}(\mathbb R^n)$.
Since $v_j=0$ a.e. in $\Omega_j^c$ and the rescaled
domains converge locally to $H$, we have $v=0$ a.e. in $H^c$.
For every fixed $T\geq1$, the strong local convergence,
\eqref{eq:boundary-residual-growth},
\eqref{eq:factor-ratio-limit} and the limit of
$\omega(Tr_j)/\omega(r_j)$ yield $\|v\|_{L^2(B_T)}
    \leq CT^{n/2+2s}.$

We next identify the limiting equation.
Let $\varphi\in C_c^\infty(H)$ satisfy
$\int_{\mathbb R^n}\varphi\,dy=0$, and choose $R\geq1$
such that $\operatorname{supp}\varphi\subset B_{2R}$.
For all sufficiently large $j$, the support of $\varphi$
is contained in $\Omega_j$.
Testing the rescaled equation in Step~4, the zero-mean
condition eliminates $d_{j,R}$ and gives
\[
    \begin{aligned}
    &\int_{\mathbb R^n}v_j(-\Delta)^s\varphi\,dy
    +\frac1{2\ell_{r_j}}
    \int_{\mathbb R^n}v_j(-\Delta)^{s+\ln}\varphi\,dy=\int_{\Omega_j\cap B_{2R}}
    g_{j,R}\varphi\,dy.
    \end{aligned}
\]
By Step~4, the right-hand side tends to zero.
Moreover, by~\eqref{eq:zero-mean-log-test-tail},
\[
    \frac1{2\ell_{r_j}}
    \left|\int_{\mathbb R^n}
    v_j(-\Delta)^{s+\ln}\varphi\,dy\right|
    \leq\frac{C_\varphi}{\ell_{r_j}}
    \int_{\mathbb R^n}|v_j|\widetilde\omega_s\,dy
    \longrightarrow0.
\]

On every fixed ball, the strong local $L^2$ convergence
gives convergence of the fractional term.
On its complement,~\eqref{eq:zero-mean-test-tail}
and~\eqref{eq:boundary-tail-tightness} give
\[
    \lim_{M\to\infty}\limsup_{j\to\infty}
    \int_{|y|>M}|v_j(-\Delta)^s\varphi|\,dy=0.
\]
The same tail estimate for $v$ follows from its
$L^2$ growth bound.
Consequently,
\[
    \int_{\mathbb R^n}v(-\Delta)^s\varphi\,dy=0.
\]
Since $v\in H^s_{\rm loc}(\mathbb R^n)$, $v=0$ a.e.
in $H^c$, and $\|v\|_{L^2(B_T)}\leq CT^{n/2+2s}$,
Corollary~\ref{lem:critical-liouville} gives
$v=\kappa(y_n)_+^s$ for some $\kappa\in\mathbb R$.
Finally, since $D(z_j)=0$ and
$\nabla D(z_j)=\nu_{\rm in}(z_j)$,
Lemma~\ref{lem:normalized-distance} gives
\[
    \frac{D(X_j(y))}{r_j}
    =y_n+O(r_j^\alpha)
\]
uniformly on $B_1$.
For all sufficiently large $j$, we have $\zeta=1$
on $X_j(B_1)$. Since $\Phi_s$ vanishes on the negative
half-line, it follows from~\eqref{eq:factor-ratio-limit}
that
\[
    \Pi_j(y)
    =\frac{\Phi_s\bigl(r_j(D(X_j(y))/r_j)_+\bigr)}
    {\Phi_s(r_j)}\rightarrow(y_n)_+^s
\]
uniformly on $B_1$. Passing to the limit in~\eqref{eq:boundary-orthogonality},
we obtain
\[
    \dashint_{B_1}|v|^2\,dy=1,\qquad
    \int_{B_1}v(y)(y_n)_+^s\,dy=0.
\]
The second identity and $v=\kappa(y_n)_+^s$ imply
$\kappa=0$, contradicting the first identity.
\end{proof}

\subsection{The boundary trace and the quotient}

We construct the boundary trace using
Proposition~\ref{prop:endpoint-excess}.
We apply the interior H\"older estimate at a radius
comparable to the boundary distance to obtain the pointwise
expansion and the modulus of continuity of the quotient.

\begin{proof}[\textnormal{\textbf{Proof of Theorem~\ref{thm:quotient}}}]
\noindent\textbf{Step 1: The boundary trace and the averaged expansion.}
By Proposition~\ref{prop:endpoint-excess}
and~\eqref{eq:endpoint-adjacent-coefficients}, for
$r\leq h\leq\min\{2r,r_0\}$, we have
\[
    |c_{z,h}-c_{z,r}|\leq C\mathcal A\omega(h).
\]
Since $\sum_{k\geq0}\omega(2^{-k}r)\leq C\omega(r)$,
the sequence $(c_{z,2^{-k}r_0})_{k\geq0}$ is Cauchy.
Denote its limit by $c(z)$.
Summing the preceding estimate over decreasing dyadic
radii gives, for every integer $k\geq0$,
\[
    |c_{z,2^{-k}r_0}-c(z)|
    \leq C\mathcal A\omega(2^{-k}r_0).
\]
For an arbitrary $r\in(0,r_0]$, choose $k\geq0$
such that $2^{-k-1}r_0<r\leq2^{-k}r_0$.
Comparing $c_{z,r}$ with $c_{z,2^{-k}r_0}$ and using
\eqref{eq:modulus-dilation}, we obtain
$|c_{z,r}-c(z)|\leq C\mathcal A\omega(r)$.
In particular, $c(z)=\lim_{r\downarrow0}c_{z,r}$.
Together with the bound
$|c_{z,r_0}|\leq Cr_0^{-n/2}\|u\|_2/\Phi_s(r_0)$,
this yields
\begin{equation}
    \frac{|c_{z,r}-c(z)|}{\omega(r)}+|c(z)|
    \leq C\mathcal A
    \qquad\text{for }0<r\leq r_0.
    \label{eq:coefficient-convergence}
\end{equation}

Combining~\eqref{eq:coefficient-convergence},
\eqref{eq:excess-decay} and~\eqref{eq:factor-ball-norm},
we obtain
\begin{equation}
    \left(\dashint_{B_r(z)}
    |u-c(z)\Psi_s|^2\,dx\right)^{1/2}
    \leq C\mathcal A\Phi_s(r)\omega(r)
    \qquad\text{for }0<r\leq r_0.
    \label{eq:endpoint-boundary-average}
\end{equation}
The coefficient $c(z)$ is uniquely determined by this
estimate. Indeed, if two coefficients satisfy
\eqref{eq:endpoint-boundary-average}, then the triangle
inequality and~\eqref{eq:factor-ball-norm} bound their
difference by $C\mathcal A\omega(r)$, which tends to zero
as $r\downarrow0$.

We next estimate the modulus of continuity of $c$.
Let $z,z'\in\partial\Omega$ and set $h=|z-z'|$.
For $0<h\leq r_0/2$, the inclusion
$B_h(z)\subset B_{2h}(z')$ gives
\[
   |c(z)-c(z')|\|\Psi_s\|_{L^2(B_h(z))}
    \leq\|u-c(z)\Psi_s\|_{L^2(B_h(z))}+\|u-c(z')\Psi_s\|_{L^2(B_{2h}(z'))}.
\]
By~\eqref{eq:endpoint-boundary-average},
\eqref{eq:factor-potter} and~\eqref{eq:modulus-dilation},
the right-hand side is bounded by
$C\mathcal A h^{n/2}\Phi_s(h)\omega(h)$.
Using the lower bound in~\eqref{eq:factor-ball-norm},
we conclude that
\begin{equation}
    |c(z)-c(z')|
    \leq C\mathcal A\omega(|z-z'|)
    \qquad\text{for }0<|z-z'|\leq1.
    \label{eq:endpoint-trace-modulus}
\end{equation}
For $r_0/2<h\leq1$, this follows from
\eqref{eq:coefficient-convergence} and the monotonicity
of $\omega$.

\medskip
\noindent\textbf{Step 2: Interior estimates for the remainder.}
For $z\in\partial\Omega$, write $e_z=u-c(z)\Psi_s$
and choose $Q_z\in O(n)$ such that
$Q_ze_n=\nu_{\rm in}(z)$.
For $0<r\leq r_0$, define
\[
    V_{z,r}(y):=\frac{e_z(z+rQ_zy)}
    {\Phi_s(r)\omega(r)}.
\]
By~\eqref{eq:endpoint-boundary-average} and a change
of variables,
\begin{equation}
    \|V_{z,r}\|_{L^2(B_T)}
    \leq C\mathcal A T^{n/2}
    \frac{\Phi_s(Tr)\omega(Tr)}
    {\Phi_s(r)\omega(r)},
    \qquad 1\leq T\leq\frac{r_0}{r}.
    \label{eq:fixed-residual-bounds}
\end{equation}

Applying~\eqref{eq:factor-potter}
and~\eqref{eq:modulus-dilation}, each with exponent
loss $1/4$, we obtain
\[
    \|V_{z,r}\|_{L^2(B_T)}
    \leq C\mathcal A T^{n/2+2s+1/2},
    \qquad 1\leq T\leq\frac{r_0}{r}.
\]
As in Step~3 of the proof of
Proposition~\ref{prop:endpoint-excess},
Cauchy--Schwarz gives, for $k\geq0$ with $2^k<r_0/r$,
\[
    \begin{aligned}
    &\int_{(B_{2^{k+1}}\setminus B_{2^k})
    \cap B_{r_0/r}}
    |V_{z,r}(y)|\widetilde\omega_s(y)\,dy\\
    &\quad\leq
    C(1+k)2^{-k(n+2s+1)}2^{kn/2}
    \|V_{z,r}\|_{L^2(B_{\min\{2^{k+1},r_0/r\}})}\leq C\mathcal A(1+k)2^{-k/2}.
    \end{aligned}
\]
For the remaining tail, we use $e_z=0$ a.e. in $\Omega^c$,
$\|e_z\|_{L^1(\Omega)}\leq C\mathcal A$ and
$\omega(r)\geq C^{-1}\Phi_s(r)$.
The same change of variables as in Step~3 of
Proposition~\ref{prop:endpoint-excess} gives
\[
  \int_{|y|>r_0/r}|V_{z,r}|\widetilde\omega_s\,dy
    \leq
    \frac{Cr^{2s+1}(1+\ell_r)}
    {\Phi_s(r)\omega(r)}
    \|e_z\|_{L^1(\Omega)}\leq C\mathcal A r(1+\ell_r)^2.
\]
Since $\sum_{k\geq0}(1+k)2^{-k/2}<\infty$ and
$r(1+\ell_r)^2$ is bounded, these estimates and
\eqref{eq:fixed-residual-bounds} with $T=1$ yield
\begin{equation}
    \int_{\mathbb R^n}|V_{z,r}|\widetilde\omega_s\,dy
    \leq C\mathcal A.
    \label{eq:fixed-residual-tail}
\end{equation}
 $V_{z,r}$ belongs to $L^2(\mathbb R^n)$ and is
compactly supported, so its ordinary weighted tail is
also finite.

Fix $\beta\in(s,\min\{1,2s\})$.
Let $x\in\Omega$ with $d=\delta(x)$ sufficiently small,
and let $z$ be any nearest boundary point to $x$.
Then $x=z+d\nu_{\rm in}(z)$.
Set $r=d/16$ and define
\[
    W(Y):=\frac{e_z(x+rY)}{\Phi_s(d)\omega(d)}.
\]
Since $W(Y)=V_{z,d}(e_n+Q_z^TY/16)$,
\eqref{eq:fixed-residual-bounds} with $T=2$ and
\eqref{eq:fixed-residual-tail} give
\[
    \|W\|_{L^2(B_8)}
    +\int_{\mathbb R^n}|W|\widetilde\omega_s\,dY
    \leq C\mathcal A.
\]
Here the fixed translation, rotation and dilation
preserve comparability of $\widetilde\omega_s$.

By Lemma~\ref{lem:factor-energy}
and Lemma~\ref{lem:dilation},
$W\in\mathcal H^{s+\ln}(\mathbb R^n)$ and
\[
    A_{\ell_r}W(Y)
    =\frac{r^{2s}}{2\ell_r\Phi_s(d)\omega(d)}
    (f-c(z)G_s)(x+rY)
\]
weakly in $B_4$, since $B_{4r}(x)\subset\Omega$.

We show that this equation can be written as
$A_{\ell_r}W=a+g$ in $B_4$, where $a\in\mathbb R$ and
$\|g\|_{L^\infty(B_4)}\leq C\mathcal A$.
Since $r=d/16$, we have $\ell_r\asymp\ell_d$ and
$\delta(x+rY)\asymp d$ for $Y\in B_4$.
Moreover,
\[
    \frac{r^{2s}}{2\ell_r\Phi_s(d)\omega(d)}
    =16^{-2s}\frac{\ell_d}{\ell_r}S(d)
    \asymp S(d).
\]

If $n=1$ or $\alpha>s$, we take $a=0$ and use the
boundedness of $G_s$ from Lemma~\ref{lem:distance-lift}
and the bound $S(d)\leq C$. If $n\geq2$ and $\alpha<s$, we also take $a=0$.
By~\eqref{eq:lift-image-bound}
and~\eqref{eq:source-scaling-factor}, for $Y\in B_4$,
\[
    S(d)|G_s(x+rY)|
    \leq Cd^{s-\alpha}L(d)^{-1/2}
    d^{\alpha-s}\sqrt{L(d)}
    \leq C.
\]

If $n\geq2$ and $\alpha=s$, let $b_{z,2d}$ be the
constant given by Lemma~\ref{lem:critical-image}.
Since $B_{4r}(x)\subset B_{2d}(z)$ and
$\delta(x+rY)\asymp d$ for $Y\in B_4$,
\eqref{eq:critical-image-compensation} gives
\[
    |G_s(x+rY)-b_{z,2d}|
    \leq C\sqrt{L(d)}.
\]
Since $S(d)\asymp L(d)^{-1/2}$, we take
\[
    a=-\frac{r^{2s}c(z)b_{z,2d}}
    {2\ell_r\Phi_s(d)\omega(d)}.
\]
The remaining right-hand side satisfies
$\|g\|_{L^\infty(B_4)}\leq C\mathcal A$. Thus, in all three cases,
\eqref{eq:coefficient-convergence} gives
$\|g\|_{L^\infty(B_4)}\leq C\mathcal A$.
No bound for $a$ is needed.

We apply Lemma~\ref{lem:compensated-localization}
with $R=2$, $\mathcal O=\mathbb R^n$ and a fixed test
ball compactly contained in $B_2$.
Since $W\in\mathcal H^{s+\ln}(\mathbb R^n)$ and its
ordinary weighted tail is finite, the preceding estimates
give a bound of $C\mathcal A$ for the quantity
$\mathcal M$ in that lemma.
Consequently, $w=\chi_2W$ satisfies
\[
    A_{\ell_r}w=g+b+F\qquad\text{in }B_4,
\]
with $|b|+\|F\|_{L^\infty(B_4)}\leq C\mathcal A$.
The equation holds for compactly supported energy
tests in $B_4$.
Since $w$ is supported in $B_8$, we have
\[
    \|w\|_{L^2(B_4)}
    +\|w\|_{L^1(\mathbb R^n,\omega_s\,dy)}
    \leq C\|W\|_{L^2(B_8)}
    \leq C\mathcal A.
\]
Lemma~\ref{lem:interior-holder} therefore gives
$\|w\|_{C^\beta(B_1)}\leq C\mathcal A$.
Since $w=W$ on $B_1$, returning to the original variables
yields
\begin{equation}
    \|e_z\|_{L^\infty(B_{d/16}(x))}
    +d^\beta[e_z]_{C^\beta(B_{d/16}(x))}
    \leq C\mathcal A\Phi_s(d)\omega(d).
    \label{eq:depth-residual-holder}
\end{equation}
Adding back $c(z)\Psi_s$ gives a continuous
representative of $u$ on each such ball.
These representatives agree on overlaps, since
they coincide a.e. with $u$.

\medskip
\noindent\textbf{Step 3: The quotient with the regularized distance.}
Choose $\varepsilon_0\in(0,r_0/4)$ sufficiently small
that the estimates in Step~2 hold for
$0<\delta(x)\leq\varepsilon_0$ and $\zeta=1$ in
a neighborhood of this collar.
Using the continuous representative of $u$ constructed
in Step~2, define $\widetilde q=u/\Phi_s(D)$ in the
interior collar.

Evaluating~\eqref{eq:depth-residual-holder} at $x$
and using $\Phi_s(D(x))\asymp\Phi_s(\delta(x))$,
we obtain, for every nearest boundary point $z$ to $x$,
\begin{equation}
    |\widetilde q(x)-c(z)|
    \leq C\mathcal A\omega(\delta(x)).
    \label{eq:endpoint-pointwise-matching}
\end{equation}
Together with~\eqref{eq:endpoint-trace-modulus},
this proves that $\widetilde q$ extends continuously
to $\partial\Omega$ with boundary value $c$.
Moreover, $\|\widetilde q\|_\infty\leq C\mathcal A$
on the closed collar.

Let $x,y$ be distinct interior points in this collar
and set $h=|x-y|\leq1$.
Suppose first that
$h<\min\{\delta(x),\delta(y)\}/32$.
Set $d=\delta(x)$ and choose a nearest boundary point
$z$ to $x$. Then $y\in B_{d/16}(x)$.
Since $D$ is Lipschitz and $D\asymp d$ on this ball,
\eqref{eq:factor-derivatives} gives
\[
    |\Phi_s(D(x))-\Phi_s(D(y))|
    \leq C\Phi_s(d)\frac hd.
\]
Moreover,
$\Phi_s(D(x))\asymp\Phi_s(D(y))\asymp\Phi_s(d)$.
Using $\widetilde q=c(z)+e_z/\Phi_s(D)$, we obtain
\[
    \begin{aligned}
    |\widetilde q(x)-\widetilde q(y)|
    &\leq\frac{|e_z(x)-e_z(y)|}{\Phi_s(D(x))}+
    \frac{|e_z(y)|\,|\Phi_s(D(x))-\Phi_s(D(y))|}
    {\Phi_s(D(x))\Phi_s(D(y))}\\
    &\leq C\mathcal A\omega(d)
    \left(\left(\frac hd\right)^\beta+\frac hd\right)\leq C\mathcal A\omega(d)
    \left(\frac hd\right)^\beta,
    \end{aligned}
\]
where the last two inequalities follow from
\eqref{eq:depth-residual-holder} and $\beta<1$.
Since $\beta>s$, choose $\varepsilon\in(0,\beta-s)$.
Applying~\eqref{eq:modulus-dilation} with $T=d/h$ gives
\[
    \omega(d)\left(\frac hd\right)^\beta
    \leq C_\varepsilon\omega(h)
    \left(\frac hd\right)^{\beta-s-\varepsilon}
    \leq C\omega(h).
\]
Consequently,
$|\widetilde q(x)-\widetilde q(y)|
\leq C\mathcal A\omega(h)$.

Suppose next that
$h\geq\min\{\delta(x),\delta(y)\}/32$.
Since $\delta$ is $1$-Lipschitz,
$\delta(x)+\delta(y)\leq65h$.
For nearest boundary points $z_x,z_y$ to $x,y$,
respectively, we also have $|z_x-z_y|\leq66h$.
For $66h\leq r_0$, estimates
\eqref{eq:endpoint-pointwise-matching}
and~\eqref{eq:endpoint-trace-modulus} give
\[
    |\widetilde q(x)-\widetilde q(y)|
    \leq C\mathcal A\bigl(
    \omega(\delta(x))+\omega(|z_x-z_y|)
    +\omega(\delta(y))\bigr).
\]
Here we extend $\omega$ by $\omega(0)=0$.
Monotonicity and~\eqref{eq:modulus-dilation} give
\[
    |\widetilde q(x)-\widetilde q(y)|
    \leq C\mathcal A\omega(h).
\]
For $r_0/66<h\leq1$, the same estimate follows
from $\|\widetilde q\|_\infty\leq C\mathcal A$
and $\omega(h)\geq\omega(r_0/66)>0$.
By continuity, the estimate also holds when either
or both points belong to $\partial\Omega$.

\medskip
\noindent\textbf{Step 4: The true distance and the boundary expansion.}
By Lemma~\ref{lem:normalized-distance}, the function
\[
    b(x):=\frac{\Phi_s(D(x))}{\Phi_s(\delta(x))}
\]
extends to a $C^{0,\alpha}$ function on a sufficiently
small closed interior collar, with $b=1$ on
$\partial\Omega$. We decrease $\varepsilon_0$
if necessary.

For $n\geq2$, the definition of $\omega$ gives
$h^\alpha\leq C\omega(h)$ for $0<h\leq1$.
Indeed, when $\alpha\leq s$, we have
$\omega(h)=h^\alpha$, whereas when $\alpha>s$,
$h^{\alpha-s}\sqrt{L(h)}$ is bounded near zero.
For $n=1$, we have $D=\delta$ in the interior collar,
and hence $b=1$.

Since $q=u/\Phi_s(\delta)=b\widetilde q$, for distinct
$x,y$ in the closed collar with $|x-y|\leq1$, we have
\[
   |q(x)-q(y)|
    \leq\|b\|_\infty
    |\widetilde q(x)-\widetilde q(y)|
    +\|\widetilde q\|_\infty|b(x)-b(y)|\leq C\mathcal A\omega(|x-y|).
\]
Moreover, $\|q\|_\infty\leq C\mathcal A$ and
$q=c$ on $\partial\Omega$.
This proves~\eqref{eq:quotient-bound}
and~\eqref{eq:quotient-modulus}.

Finally, let $z\in\partial\Omega$ and
$0<r\leq\varepsilon_0$.
For $x\in\Omega\cap B_r(z)$, we have
$\delta(x)\leq|x-z|\leq r$.
By the quotient estimate and the monotonicity of
$\Phi_s$ and $\omega$,
\[
    \begin{aligned}
    |u(x)-c(z)\Phi_s(\delta(x))|
    &=\Phi_s(\delta(x))|q(x)-q(z)|\leq C\mathcal A\Phi_s(r)\omega(r).
    \end{aligned}
\]
Since $u=0$ a.e. in $\Omega^c$ and
$\Phi_s(\delta_0)=0$ in $\Omega^c$, we obtain
\[
    \|u-c(z)\Phi_s(\delta_0)\|_{L^\infty(B_r(z))}
    \leq C\mathcal A\Phi_s(r)\omega(r).
\]
This proves~\eqref{eq:boundary-supremum}.
The corresponding nearest-point estimate follows
by choosing $z$ to be a nearest boundary point to $x$,
and the averaged estimate follows by integration
over $B_r(z)$.
\end{proof}

\begin{proof}[\textnormal{\textbf{Proof of the decay estimate in
Corollary~\ref{cor:boundedness}}}]
The global $L^\infty$ estimate was established using
Lemma~\ref{lem:local-boundedness}
and Lemma~\ref{lem:dilation}.
By Theorem~\ref{thm:quotient}, a.e. in a fixed interior
boundary collar,
\[
    |u(x)|\leq C\mathcal A\Phi_s(\delta(x))
    \leq C\mathcal A\rho_{s,R}(x),
\]
where the second inequality follows from
Proposition~\ref{prop:factor-asymptotics}.
On the complement of the collar, $\rho_{s,R}$
has a positive lower bound.
The global $L^\infty$ estimate therefore gives
\eqref{eq:boundary-decay} a.e. in $\Omega$.
\end{proof}

\begin{proof}[\textnormal{\textbf{Proof of Corollary~\ref{cor:resolvent}}}]
Since $f\in L^\infty(\Omega)\subset L^2(\Omega)$
and zero belongs to the resolvent set of the
self-adjoint Dirichlet realization, the representation
theorem for its quadratic form gives the unique
variational solution
$u=((-\Delta)^{s+\ln}_\Omega)^{-1}f$.
By the spectral theorem and~\eqref{eq:resolvent-distance},
\[
    \|u\|_2
    \leq d_{s,\Omega}^{-1}\|f\|_2
    \leq d_{s,\Omega}^{-1}|\Omega|^{1/2}\|f\|_\infty.
\]
Consequently,
$\mathcal A\leq
(1+d_{s,\Omega}^{-1}|\Omega|^{1/2})\|f\|_\infty$.
Theorem~\ref{thm:quotient}
and Corollary~\ref{cor:boundedness}
now give the asserted estimates.
\end{proof}

\section{Counterexamples and sharpness}
\label{sec:sharpness}

We first construct, on every bounded $C^{1,\alpha}$ open set,
a solution $U$ with bounded right-hand side and strictly
positive boundary trace of $U/\Phi_s(\delta)$.
Using the first two correction terms in the asymptotic
expansion of $\Phi_s$, we show that, for every
$R>e\operatorname{diam}(\Omega)$,
the quotient $U/\rho_{s,R}$ fails to be H\"older continuous
of any positive order near every boundary point.

We then prove the optimality of the quotient modulus
$\omega_{n,s,\alpha}$ and the corresponding boundary remainder
estimates. Interior bumps give the lower bound $\mu_s$
on every fixed domain. For $n\geq2$ and $0<\alpha\leq s$,
we construct a domain whose boundary is locally given by
$x_n=-|x'|^{1+\alpha}$ and obtain a lower bound of order
$\Phi_s(r)r^\alpha$ for the averaged boundary remainder.
This also gives the lower bound $r^\alpha$ for the quotient.
On the same domain, we rule out a nonzero constant trace near the
origin, as asserted in Theorem~\ref{thm:bare-quotient}.

\subsection{A positive boundary trace and the leading-order factor}

When $n\geq2$ and $\alpha\leq s$, the function $G_s$
may be unbounded near the boundary.
To construct a solution with a strictly positive boundary
trace, we add a higher power of the regularized distance
to $\Psi_s$, following the subsolution construction of
Ros-Oton and Serra~\cite[Lemmas~2.7 and~2.10]
{RosOtonSerra2017Boundary}.
The added power has a negative operator image that
dominates the possible singularity of $G_s$.
Since the kernel of $(-\Delta)^{s+\ln}$ changes sign,
we perform the comparison using its positive part.

\begin{lemma}
\label{lem:positive-boundary-trace}
There exist $U\in\mathcal H_0^{s+\ln}(\Omega)$ and
$f\in L^\infty(\Omega)$ such that $U$ solves
\eqref{eq:dirichlet}, $0\leq U\leq1$ a.e. in $\Omega$,
and the boundary trace $c$ of $U/\Phi_s(\delta)$ satisfies
$\min_{\partial\Omega}c>0$.
\end{lemma}

\begin{proof}[\textnormal{\textbf{Proof}}]
\noindent\textbf{Step 1: Construction of a nonnegative solution.}
Write $K_{s+\ln}=K^+-K^-$, where
$K^\pm=\max\{\pm K_{s+\ln},0\}$.
Since $K_{s+\ln}$ is positive near the origin and
absolutely integrable away from it, we have
$N_s:=\int_{\mathbb R^n}K^-(|y|)\,dy<\infty$.
Define
\[
    \mathcal E_+(v,w)
    :=\frac12\iint_{\mathbb R^{2n}}
    (v(x)-v(y))(w(x)-w(y))
    K^+(|x-y|)\,dx\,dy.
\]
The kernel estimates used in~\eqref{eq:shifted-form}
give $\mathcal E_+(v,v)+\|v\|_2^2
    \asymp\|v\|_{\mathcal H^{s+\ln}}^2.$
Thus the bilinear form
$\mathcal E_+(\cdot,\cdot)+(\cdot,\cdot)_{L^2}$
is continuous and coercive on
$\mathcal H_0^{s+\ln}(\Omega)$.
Since $|\int_\Omega\varphi\,dx|
\leq|\Omega|^{1/2}\|\varphi\|_2$,
the Lax--Milgram theorem~\cite[Corollary~5.8]{Brezis2011}
gives a unique $U\in\mathcal H_0^{s+\ln}(\Omega)$
such that
\begin{equation}
    \mathcal E_+(U,\varphi)+\int_\Omega U\varphi\,dx
    =\int_\Omega\varphi\,dx
    \label{eq:positive-auxiliary-solution}
\end{equation}
for every $\varphi\in\mathcal H_0^{s+\ln}(\Omega)$.

The function $U^-=\max\{-U,0\}$ is an admissible test,
since the negative-part map is Lipschitz and vanishes
at zero. Using
$(a-b)(a^--b^-)\leq-|a^--b^-|^2$ and $K^+\geq0$,
we obtain
\[
    0\leq\int_\Omega U^-\,dx
    \leq-\mathcal E_+(U^-,U^-)-\|U^-\|_2^2.
\]
Hence $U\geq0$ a.e. in $\Omega$.
Next, taking $v=(U-1)_+$ as a test in
\eqref{eq:positive-auxiliary-solution} and using
\[
    (a-b)((a-1)_+-(b-1)_+)
    \geq|(a-1)_+-(b-1)_+|^2,
\]
we obtain $\mathcal E_+(v,v)+\|v\|_2^2\leq0$.
Therefore $0\leq U\leq1$ a.e. in $\Omega$.

For $v\in\mathcal H_0^{s+\ln}(\Omega)$, write
$(K^-*v)(x)=\int_{\mathbb R^n}K^-(|x-y|)v(y)\,dy$.
Since $K^-$ is integrable, symmetry and Fubini's
theorem give, for every
$\varphi\in\mathcal H_0^{s+\ln}(\Omega)$,
\[
    \begin{aligned}
    \mathcal E_{s+\ln}(v,\varphi)
    &=\mathcal E_+(v,\varphi)
    -N_s\int_\Omega v\varphi\,dx+\int_\Omega(K^-*v)\varphi\,dx.
    \end{aligned}
\]
Applying this identity to $U$ and using
\eqref{eq:positive-auxiliary-solution}, we conclude
that $U$ solves~\eqref{eq:dirichlet} with
\[
    f=1-(1+N_s)U+K^-*U.
\]
Since $0\leq U\leq1$ and $K^-\in L^1(\mathbb R^n)$,
we have $f\in L^\infty(\Omega)$.

\medskip
\noindent\textbf{Step 2: Construction of a comparison function.}
Fix $0<\beta<\min\{\alpha,s\}$ and set
\[
    V=\zeta(D_+)^{s+\beta},
    \qquad B=\Psi_s+V,
\]
where $D$ and $\zeta$ are as in
Lemma~\ref{lem:normalized-distance}
and~\eqref{eq:regularized-factor}.
We extend $V$ by zero outside $\mathcal U$.
Then $B$ is nonnegative, bounded and vanishes in
$\Omega^c$.

By~\eqref{eq:distance-factor-ratio} and
Proposition~\ref{prop:factor-asymptotics}, uniformly
as $\delta(x)\downarrow0$,
\begin{equation}
    \frac{B(x)}{\Phi_s(\delta(x))}
    =1+O(\delta(x)^\alpha)
    +O\bigl(\delta(x)^\beta\sqrt{L(\delta(x))}\bigr)
    \longrightarrow1.
    \label{eq:comparison-boundary-ratio}
\end{equation}
We claim that, uniformly as $\delta(x)\downarrow0$,
\begin{equation}
    (-\Delta)^{s+\ln}V(x)
    =2F(s,s+\beta)\delta(x)^{\beta-s}
    \ln\frac1{\delta(x)}\bigl(1+o(1)\bigr),
    \label{eq:negative-power-image}
\end{equation}
where $F$ is defined in~\eqref{eq:half-line-coefficient}.
Since $0<\beta<s<1$, we have
\[
    F(s,s+\beta)
    =-\frac{\Gamma(1+s+\beta)\Gamma(s-\beta)}{\pi}
    \sin(\pi\beta)<0.
\]

To prove~\eqref{eq:negative-power-image}, let
$t=\delta(x)$, choose a nearest boundary point $z$,
and set $p=s+\beta$.
Writing $\ell_z(y)=(y-z)\cdot\nu_{\rm in}(z)$,
we compare $V$ with $(\ell_z)_+^p$ and set
$E=V-(\ell_z)_+^p$.

In $B_{t/4}(x)$, both $D$ and $\ell_z$ are comparable
to $t$, and $\zeta=1$ for sufficiently small $t$.
By Lemma~\ref{lem:normalized-distance} and the
$C^{1,\alpha}$ regularity of $D$, we have
\[
    |D-\ell_z|+t|\nabla D-\nu_{\rm in}(z)|
    +t^2|\nabla^2D|
    \leq Ct^{1+\alpha}.
\]
Applying the chain rule to $E=D^p-\ell_z^p$ as in
Step~2 of Lemma~\ref{lem:distance-lift}, we obtain
\[
    \sup_{B_{t/4}(x)}|\nabla^2E|
    \leq C\left(
    t^{p-3}t^{1+\alpha}
    +t^{p-2}t^\alpha
    +t^{p-1}t^{\alpha-1}\right)\leq Ct^{p+\alpha-2}.
\]

We next estimate the operator image of $E$ by rescaling.
For every fixed $R\geq1$, the $C^{1,\alpha}$ estimate
for $D$ gives
\[
    \sup_{y\in B_{Rt}(x)}
    |D(y)-\ell_z(y)|\leq C_Rt^{1+\alpha}.
\]
On this ball, both $D/t$ and $\ell_z/t$ are uniformly
bounded, and $\zeta=1$ for sufficiently small $t$.
By the uniform continuity of the positive-part power
on bounded intervals, it follows that
\[
    \sup_{|h|\leq R}
    t^{-p}|E(x+th)|\longrightarrow0
    \qquad\text{as }t\downarrow0,
\]
uniformly in $x$ and its nearest boundary point $z$.

Moreover, $V(y)\leq C\delta_0(y)^p\leq C|y-z|^p$
and $(\ell_z(y))_+^p\leq|y-z|^p$ on $\mathbb R^n$.
Consequently,
\[
    t^{-p}|E(x+th)|\leq C(1+|h|)^p
    \qquad\text{for every }h\in\mathbb R^n.
\]
The symmetric second-difference representation and
the change of variables at scale $t$ give
\[
  \frac{(-\Delta)^{s+\ln}E(x)}
    {t^{p-2s}L(t)}=\frac12\int_{\mathbb R^n}
    \frac{2E(x)-E(x+th)-E(x-th)}{t^p}
    \frac{t^{n+2s}K_{s+\ln}(t|h|)}{L(t)}\,dh.
\]
By the expression for the kernel, for $0<t<1$,
\[
    \frac{t^{n+2s}|K_{s+\ln}(t|h|)|}{L(t)}
    \leq C|h|^{-n-2s}(1+|\ln|h||).
\]
For $|h|<1/4$, the Hessian estimate gives
\[
    \frac{|2E(x)-E(x+th)-E(x-th)|}{t^p}
    \leq Ct^\alpha|h|^2.
\]
For $|h|\geq1/4$, the same expression is bounded
by $C(1+|h|)^p$.
Since $s<1$ and $p=s+\beta<2s$, these bounds provide
an integrable majorant independent of $x$ and $t$.
The preceding local convergence and dominated
convergence therefore yield
\[
    (-\Delta)^{s+\ln}E(x)
    =o\bigl(t^{\beta-s}L(t)\bigr)
    \qquad\text{as }t\downarrow0,
\]
uniformly along the boundary.

Since $p<2s$, the untruncated flat power
$(\ell_z)_+^p$ satisfies the tail condition in
\eqref{eq:tangential-reduction} and is smooth near $x$.
By Lemma~\ref{lem:half-line-powers},
\[
    (-\Delta)^{s+\ln}(\ell_z)_+^{s+\beta}(x)
    =t^{\beta-s}
    \left(e_s(\beta)+2d_s(\beta)\ln\frac1t\right),
\]
where $d_s(\beta)=F(s,s+\beta)$.
Combining this identity with the estimate for $E$
proves~\eqref{eq:negative-power-image}.

\medskip
\noindent\textbf{Step 3: Weak comparison and positivity of the trace.}
Since $\beta<\min\{\alpha,s\}$,
Lemma~\ref{lem:distance-lift} gives
\[
    G_s(x)
    =o\bigl(\delta(x)^{\beta-s}L(\delta(x))\bigr)
    \qquad\text{as }\delta(x)\downarrow0.
\]
Together with~\eqref{eq:negative-power-image} and
$F(s,s+\beta)<0$, this implies
\[
    (-\Delta)^{s+\ln}B(x)\longrightarrow-\infty
    \qquad\text{uniformly as }\delta(x)\downarrow0.
\]
Since $B$ is bounded on $\mathbb R^n$ and smooth in
$\Omega$, its operator image is locally bounded in
$\Omega$ and hence bounded above throughout $\Omega$.

We verify that the pointwise estimate can be used
in the weak formulation.
Since $s+\beta>s$, the function $V$ belongs to
$C^s(\mathbb R^n)$ and is smooth in $\Omega$.
Moreover, the derivative estimates for $D$ and $\Phi_s$
give $|\nabla B|\leq C\delta^{s-1}$ and
$|\nabla^2B|\leq C\delta^{s-2}$ in the boundary collar.
Thus $B$ satisfies the same bounds used in the proof
of Lemma~\ref{lem:factor-energy}.

In particular, splitting the energy integral at
$|h|=\delta(x)/2$ gives an integrable energy density
bounded by $CL(\delta(x))^2$ near the boundary.
Hence $B\in\mathcal H_0^{s+\ln}(\Omega)$.
The symmetric second-difference estimate also gives
\[
    |(-\Delta)^{s+\ln}_\varepsilon B(x)|
    \leq C\delta(x)^{-s}L(\delta(x))
\]
in the boundary collar, uniformly in $\varepsilon>0$,
where $(-\Delta)^{s+\ln}_\varepsilon$ denotes the
operator with the kernel truncated on $|h|\leq\varepsilon$.
The truncated operators are uniformly bounded away
from the collar.
As in Step~2 of Lemma~\ref{lem:factor-energy},
dominated convergence yields
\[
    \delta^s(-\Delta)^{s+\ln}_\varepsilon B
    \longrightarrow
    \delta^s(-\Delta)^{s+\ln}B
    \qquad\text{in }L^2(\Omega).
\]
Using~\eqref{eq:zero-exterior-hardy} and passing to
the limit in the truncated form identity as in
Lemma~\ref{lem:factor-energy}, we obtain
\[
    \mathcal E_{s+\ln}(B,\varphi)
    =\int_\Omega((-\Delta)^{s+\ln}B)\varphi\,dx
\]
for every $\varphi\in\mathcal H_0^{s+\ln}(\Omega)$,
with absolute convergence of the integral.

Define
\[
    H_B:=(-\Delta)^{s+\ln}B+(1+N_s)B-K^-*B.
\]
Since $B$ is bounded, $K^-$ is integrable and
$(-\Delta)^{s+\ln}B$ is bounded above, we may choose
$M\geq1$ such that $H_B\leq M$ in $\Omega$.
The form decomposition in Step~1 then gives, for
every nonnegative
$\varphi\in\mathcal H_0^{s+\ln}(\Omega)$,
\[
\mathcal E_+(B/M,\varphi)
    +\int_\Omega(B/M)\varphi\,dx
    =\frac1M\int_\Omega H_B\varphi\,dx\leq\int_\Omega\varphi\,dx.
\]
Subtracting~\eqref{eq:positive-auxiliary-solution},
we obtain
\[
    \mathcal E_+(B/M-U,\varphi)
    +\int_\Omega(B/M-U)\varphi\,dx\leq0.
\]
Taking $v=(B/M-U)_+$ as a test and using
$(a-b)(a_+-b_+)\geq|a_+-b_+|^2$ together with
$K^+\geq0$, we conclude that
$\mathcal E_+(v,v)+\|v\|_2^2\leq0$.
Thus $v=0$, and hence $U\geq B/M$ a.e. in $\Omega$. Since $f\in L^\infty(\Omega)$,
Theorem~\ref{thm:quotient} gives a continuous extension
of $U/\Phi_s(\delta)$ to the boundary.
The comparison inequality and
\eqref{eq:comparison-boundary-ratio} therefore yield
$c(z)\geq M^{-1}$ for every $z\in\partial\Omega$.
\end{proof}

We next show that the quotient by $\rho_{s,R}$ need not
be H\"older continuous, even when its boundary trace
is strictly positive. The first nonzero logarithmic
correction in the expansion of $\Phi_s(\delta)/\rho_{s,R}$
dominates the H\"older error in the corrected quotient.

\begin{proof}[\textnormal{\textbf{Proof of
Theorem~\ref{thm:bare-quotient}, parts~\textup{(i)} and~\textup{(ii)}}}]
By Lemma~\ref{lem:positive-boundary-trace}, there exists
a solution $U$ of~\eqref{eq:dirichlet} with
$f\in L^\infty(\Omega)$ and strictly positive boundary
trace. If $n=1$ or $\alpha>s$, we may take $U=\Psi_s$
and $f=G_s$. Indeed, Lemmas~\ref{lem:distance-lift}
and~\ref{lem:factor-energy} give $G_s\in L^\infty(\Omega)$
and the variational identity, while
\eqref{eq:distance-factor-ratio} gives boundary trace one.

Set $q=U/\Phi_s(\delta)$ and $c=q|_{\partial\Omega}$.
By Theorem~\ref{thm:quotient}, $q\in C^{0,\theta}$ on a
sufficiently small closed boundary collar, where
$\theta=\min\{\alpha,s\}$, and $\min_{\partial\Omega}c>0$.
Proposition~\ref{prop:factor-asymptotics} gives
\[
    \frac{U(x)}{\rho_{s,R}(x)}
    =q(x)\left(
    1+\frac{A_{1,s,R}}{L_R(\delta(x))}
    +\frac{A_{2,s,R}}{L_R(\delta(x))^2}
    +O(L_R(\delta(x))^{-3})
    \right)
\]
as $\delta(x)\downarrow0$.
In particular, $U/\rho_{s,R}$ extends continuously
to $\partial\Omega$ with boundary trace $c$.

Fix $z\in\partial\Omega$ and set
$x_t=z+t\nu_{\rm in}(z)$.
By the local $C^1$ graph representation at $z$,
we have $B_{t/2}(x_t)\subset\Omega$ for sufficiently
small $t>0$. Hence $t/2\leq\delta(x_t)\leq t$.
The H\"older estimate for $q$ yields
$q(x_t)=c(z)+O(t^\theta)$.

Writing $\Lambda_t=L_R(\delta(x_t))$, we have
$\Lambda_t=\ln(1/t)+O(1)$ and
$t^\theta=o(\Lambda_t^{-k})$ for every $k>0$.
Thus the preceding expansion gives
\[
    \frac{U(x_t)}{\rho_{s,R}(x_t)}-c(z)
    =\frac{c(z)A_{1,s,R}}{\Lambda_t}
    +\frac{c(z)A_{2,s,R}}{\Lambda_t^2}
    +O(\Lambda_t^{-3}).
\]
Since $A_{2,s,R}>0$ by~\eqref{eq:trigamma-positive},
we may set $m=1$ if $A_{1,s,R}\neq0$ and $m=2$ otherwise.
In either case,
\[
    \lim_{t\downarrow0}\Lambda_t^m
    \left(\frac{U(x_t)}{\rho_{s,R}(x_t)}-c(z)\right)
    =c(z)A_{m,s,R}\neq0.
\]
Since $|x_t-z|=t$ and $t^\gamma\Lambda_t^m\to0$
for every $\gamma>0$, we obtain
\[
    \lim_{t\downarrow0}
    \frac{|U(x_t)/\rho_{s,R}(x_t)-c(z)|}
    {|x_t-z|^\gamma}
    =+\infty.
\]
Consequently, $U/\rho_{s,R}$ admits no $C^{0,\gamma}$
extension to any closed boundary neighborhood of $z$
for any $\gamma\in(0,1]$.
Higher H\"older regularity is also excluded.
Since $z\in\partial\Omega$ and
$R>e\operatorname{diam}(\Omega)$ were arbitrary,
the proof is complete.
\end{proof}

\subsection{Interior bumps and sharpness of \texorpdfstring{$\mu_s$}{mu-s}}

We prove the sharpness of $\mu_s$ using smooth functions
supported in interior balls whose radii are comparable
to their boundary distances. Their right-hand sides
are uniformly bounded on each fixed domain.

\begin{proof}[\textnormal{\textbf{Proof of
Theorem~\ref{thm:sharpness}, part~\textup{(i)}}}]
\noindent\textbf{Step 1: Sharpness of the uniform modulus.}
Fix $z\in\partial\Omega$ and choose
$\eta\in C_c^\infty(B_{1/8})$ with $\eta\geq0$
and $\eta(0)=1$.
For sufficiently small $r>0$, set
$x_r=z+r\nu_{\rm in}(z)$.
By the local $C^1$ graph representation of
$\partial\Omega$, we have $B_{r/2}(x_r)\subset\Omega$.
In particular, $r/2\leq\delta(x_r)\leq r$.

Define
\[
    u_r(x):=\frac{r^{2s}}{\ell_r}
    \eta\left(\frac{x-x_r}{r}\right),
    \qquad
    f_r:=(-\Delta)^{s+\ln}u_r\big|_\Omega.
\]
Then $u_r\in C_c^\infty(\Omega)$ and
$\delta\asymp r$ on its support.
By Lemma~\ref{lem:dilation}, for every
$Y\in\mathbb R^n$,
\[
    (-\Delta)^{s+\ln}u_r(x_r+rY)
    =2(-\Delta)^s\eta(Y)
    +\frac1{\ell_r}(-\Delta)^{s+\ln}\eta(Y).
\]
Both operator images on the right-hand side are
bounded on $\mathbb R^n$, since the symmetric second
differences of $\eta$ are bounded by
$C\min\{|h|^2,1\}$, which is integrable against the
absolute values of the corresponding kernels.

Moreover,
$\|u_r\|_2=r^{n/2+2s}\|\eta\|_2/\ell_r$.
Consequently, $\|u_r\|_2+\|f_r\|_{L^\infty(\Omega)}\leq C$
with $C$ independent of $r$.
Since $u_r\in C_c^\infty(\Omega)$, we have
$u_r\in\mathcal H_0^{s+\ln}(\Omega)$.
The truncated operator images are uniformly bounded
for each fixed $r$. By symmetry of the truncated kernels
and dominated convergence, as in the proof of
Lemma~\ref{lem:factor-energy}, we obtain
\[
    \mathcal E_{s+\ln}(u_r,\varphi)
    =\int_\Omega f_r\varphi\,dx
    \qquad
    \text{for every }\varphi\in\mathcal H_0^{s+\ln}(\Omega).
\]
Thus $u_r$ is a variational solution of~\eqref{eq:dirichlet}.

Let $q_r=u_r/\Phi_s(\delta)$ in the boundary collar.
Since $u_r$ is compactly supported in $\Omega$,
the boundary trace of $q_r$ is identically zero.
By Proposition~\ref{prop:factor-asymptotics} and
$\delta(x_r)\asymp r$, we have
\[
    q_r(x_r)
    =\frac{r^{2s}}{\ell_r\Phi_s(\delta(x_r))}
    \asymp\frac{r^s}{\sqrt{L(r)}}
    =\mu_s(r).
\]
Since $|x_r-z|=r$ and $q_r(z)=0$, no modulus
$o(\mu_s)$ can replace $\mu_s$ in an estimate uniform
over solutions with bounded $\mathcal A$. The same family gives a lower bound for the remainder
in~\eqref{eq:boundary-supremum}.
Indeed, its boundary trace is zero and
\[
    \|u_r\|_{L^\infty(B_{2r}(z))}
    \geq u_r(x_r)
    \asymp\Phi_s(r)\mu_s(r)
    \asymp\Phi_s(2r)\mu_s(2r).
\]

\medskip
\noindent\textbf{Step 2: A single solution without higher
H\"older regularity.}
Choose $r_*>0$ sufficiently small and set
$r_j=r_*e^{-j^2}$ for $j\geq1$.
All the functions $u_{r_j}$ are supported in the
boundary collar. Their supports are pairwise disjoint:
indeed, for $k>j$, we have $r_k<r_j/4$ and hence
\[
    |x_{r_j}-x_{r_k}|=r_j-r_k
    >\frac34r_j>\frac18(r_j+r_k).
\]

Define
\[
    u:=\sum_{j=1}^\infty2^{-j}u_{r_j},
    \qquad
    f:=\sum_{j=1}^\infty2^{-j}f_{r_j}.
\]
By Step~1, the series defining $f$ converges in
$L^\infty(\Omega)$.
Moreover, the shifted form estimate
\eqref{eq:shifted-form} and the variational identity
for $u_r$ give
\[
     \|u_r\|_{\mathcal H^{s+\ln}}^2
    \leq C\bigl(
    \|f_r\|_\infty\|u_r\|_{L^1(\Omega)}
    +\|u_r\|_2^2\bigr)\leq C.
\]
Since $\sum_{j\geq1}2^{-j}<\infty$ and
$\mathcal H_0^{s+\ln}(\Omega)$ is complete,
the series defining $u$ converges in
$\mathcal H_0^{s+\ln}(\Omega)$.
Passing to the limit in the variational identities
for the partial sums, we conclude that $u$ solves
\eqref{eq:dirichlet} with right-hand side $f$.
The bound $\|u_{r_j}\|_\infty\leq Cr_j^{2s}/\ell_{r_j}$
also gives uniform convergence.

Let $q=u/\Phi_s(\delta)$ in the boundary collar.
By the disjointness of the supports, on
$\operatorname{supp}u_{r_j}$ we have
$u=2^{-j}u_{r_j}$ and $\delta\asymp r_j$.
Consequently,
\[
    |q(x)|\leq C2^{-j}\mu_s(r_j)
    \qquad\text{for }x\in\operatorname{supp}u_{r_j},
\]
while $q$ vanishes outside the union of these supports.
Since each finite union of the supports has positive
distance from $\partial\Omega$ and
$2^{-j}\mu_s(r_j)\to0$, the boundary trace of $q$
is identically zero.

On the other hand, at the centers $x_{r_j}$,
Step~1 gives
$|q(x_{r_j})-q(z)|\geq c2^{-j}\mu_s(r_j)$.
For every $\varepsilon\in(0,1-s)$, we obtain
\[
       \frac{|q(x_{r_j})-q(z)|}
    {|x_{r_j}-z|^{s+\varepsilon}}
    \geq c\frac{2^{-j}r_j^{-\varepsilon}}
    {\sqrt{L(r_j)}}=c\frac{r_*^{-\varepsilon}
    e^{\varepsilon j^2-j\ln2}}
    {\sqrt{L(r_*)+j^2}}
    \rightarrow+\infty.
\]
Thus $q$ admits no $C^{s+\varepsilon}$ extension near
$z$ for any $\varepsilon\in(0,1-s)$.
Since higher H\"older regularity would imply
$C^\gamma$ regularity for some $\gamma\in(s,1)$,
the conclusion holds for every $\varepsilon>0$.
This completes the proof of part~\textup{(i)}.
\end{proof}

\subsection{Rough boundaries and sharpness of
\texorpdfstring{$r^\alpha$}{r-alpha}}

We consider $n\geq2$ and $0<\alpha\leq s$.
We construct a $C^{1,\alpha}$ domain and a solution
with bounded right-hand side for which the local
$L^1$ average of the boundary remainder is bounded
below by a positive multiple of $\Phi_s(r)r^\alpha$.
This also proves the sharpness of the quotient
modulus $r^\alpha$.

When $\alpha=s$, the limiting function is homogeneous
of degree $2s$ and need not have an integrable
fractional tail. The following lemma shows that
such a function cannot satisfy the fractional equation
modulo constants unless its angular mean vanishes.
We denote by $d\sigma$ the surface measure on
$\mathbb S^{n-1}$.

\begin{lemma}
\label{lem:critical-homogeneous}
Let $H=\{x_n>0\}$ and let
$W\in L^1_{\rm loc}(\mathbb R^n)$ be homogeneous
of degree $2s$, with angular part
$W|_{\mathbb S^{n-1}}\in L^1(\mathbb S^{n-1})$.
Suppose that $\int_{\mathbb R^n}W(-\Delta)^s\varphi\,dx=0$
for every $\varphi\in C_c^\infty(H)$ with
$\int_{\mathbb R^n}\varphi\,dx=0$.
Then
\[
    \int_{\mathbb S^{n-1}}W(\theta)\,d\sigma(\theta)=0.
\]
\end{lemma}

\begin{proof}[\textnormal{\textbf{Proof}}]
Choose $\varphi\in C_c^\infty(H)$ with
$\int_{\mathbb R^n}\varphi\,dx=1$.
For $\lambda>1$, set
$\varphi_\lambda(x)=\lambda^{-n}\varphi(x/\lambda)$.
Then $\varphi_\lambda-\varphi\in C_c^\infty(H)$
has zero integral.

For $|x|$ sufficiently large, the integral representation
of the fractional Laplacian gives
\[
    (-\Delta)^s\varphi(x)
    =-c_{n,s}\int_{\mathbb R^n}
    \frac{\varphi(y)}{|x-y|^{n+2s}}\,dy.
\]
Since $\varphi$ is compactly supported and has integral one,
the mean value theorem yields
\[
    (-\Delta)^s\varphi(x)
    =-c_{n,s}|x|^{-n-2s}
    +O_\varphi(|x|^{-n-2s-1}).
\]
By~\eqref{eq:zero-mean-test-tail} and the homogeneity
of $W$, the pairing of $W$ with
$(-\Delta)^s(\varphi_\lambda-\varphi)$ is absolutely
convergent. Indeed, its tail is bounded by a constant
times $\|W\|_{L^1(\mathbb S^{n-1})}
    \int_1^\infty r^{-2}\,dr.$

On the other hand, using
$W(\lambda x)=\lambda^{2s}W(x)$ and
$(-\Delta)^s\varphi_\lambda(x)
=\lambda^{-n-2s}((-\Delta)^s\varphi)(x/\lambda)$,
a change of variables gives
\[
    \int_{B_R}W(-\Delta)^s\varphi_\lambda\,dx
    =\int_{B_{R/\lambda}}W(-\Delta)^s\varphi\,dx.
\]
Consequently,
\[
    \int_{B_R}W(-\Delta)^s(\varphi_\lambda-\varphi)\,dx
    =-\int_{B_R\setminus B_{R/\lambda}}
    W(-\Delta)^s\varphi\,dx.
\]
By the preceding asymptotic formula and polar coordinates,
we obtain
\[
-\int_{B_R\setminus B_{R/\lambda}}
    W(-\Delta)^s\varphi\,dx=c_{n,s}\ln\lambda
    \int_{\mathbb S^{n-1}}W(\theta)\,d\sigma(\theta)
    +O_{\varphi,\lambda}(R^{-1})
    \|W\|_{L^1(\mathbb S^{n-1})}.
\]
Here we used $W(r\theta)=r^{2s}W(\theta)$,
$\int_{R/\lambda}^R r^{-1}\,dr=\ln\lambda$ and
$\int_{R/\lambda}^R r^{-2}\,dr=(\lambda-1)/R$.

Letting $R\to\infty$ and using the equation for
$\varphi_\lambda-\varphi$, we conclude that
\[
    0=c_{n,s}\ln\lambda
    \int_{\mathbb S^{n-1}}W(\theta)\,d\sigma(\theta).
\]
Since $c_{n,s}>0$ and $\lambda>1$, the assertion follows.
\end{proof}

We use the same domain and blow-up argument for
$\alpha<s$ and $\alpha=s$.
Only the exclusion of the limiting function differs:
for $\alpha<s$, we compute the sign of its fractional
Laplacian, whereas for $\alpha=s$, we apply
Lemma~\ref{lem:critical-homogeneous}.

\begin{proof}[\textnormal{\textbf{Proof of
Theorem~\ref{thm:sharpness}, part~\textup{(ii)}}}]
\noindent\textbf{Step 1: Construction of the domain
and selection of scales.}
Let $n\geq2$ and $0<\alpha\leq s$.
Writing $x=(x',x_n)\in\mathbb R^{n-1}\times\mathbb R$,
we choose a bounded connected $C^{1,\alpha}$ domain
$\Omega$ such that $0\in\partial\Omega$ and for some $r_0>0$,
\begin{equation}
    \Omega\cap B_{r_0}
    =\{x_n>-|x'|^{1+\alpha}\}\cap B_{r_0}.
    \label{eq:rough-boundary-domain}
\end{equation}

Such a domain can be obtained by modifying the lower
boundary of $B_1(e_n)$. Indeed, let
$b(x')=1-\sqrt{1-|x'|^2}$ and replace $b$ by
\[
    (1-\chi_0(x'))b(x')
    -\chi_0(x')|x'|^{1+\alpha},
\]
where $\chi_0\in C_c^\infty(B_{1/4}^{n-1})$,
$0\leq\chi_0\leq1$, and $\chi_0=1$ near zero.
The new graph is $C^{1,\alpha}$, lies below the original
one and agrees with it outside the support of $\chi_0$.
Retaining the rest of the sphere gives a bounded
connected $C^{1,\alpha}$ domain satisfying
\eqref{eq:rough-boundary-domain}.

Lemma~\ref{lem:positive-boundary-trace} ensures the
existence of solutions with strictly positive boundary trace.
Let $u$ be any variational solution of~\eqref{eq:dirichlet}
with $f\in L^\infty(\Omega)$ and nonzero boundary trace
$c:=q(0)$, where $q=u/\Phi_s(\delta)$.
By replacing $(u,f)$ with $(-u,-f)$ if necessary,
we may assume that $c>0$.
By Theorem~\ref{thm:quotient}, $q\in C^{0,\alpha}$
on a sufficiently small closed boundary collar.

We claim that
\begin{equation}
    \liminf_{r\downarrow0}
    \frac{\dashint_{B_r}|u-c\Phi_s(\delta_0)|\,dx}
    {\Phi_s(r)r^\alpha}>0.
    \label{eq:rough-remainder-lower}
\end{equation}
Suppose, to the contrary, that there exist radii
$r_j\downarrow0$ such that
\[
    \varepsilon_j:=
    \frac{\dashint_{B_{r_j}}
    |u-c\Phi_s(\delta_0)|\,dx}
    {\Phi_s(r_j)r_j^\alpha}
    \longrightarrow0.
\]
Fix $\eta\in(0,1)$ and choose $T_j\to\infty$ sufficiently
slowly that $\varepsilon_jT_j^{n+s+\alpha+\eta}\to0$.
Set $\rho_j=r_j/T_j$.

For every fixed $T\geq1$, we have $T\rho_j\leq r_j$
for all sufficiently large $j$. Changing variables
and using~\eqref{eq:factor-potter}, we obtain
\[
    \int_{B_T}
    \frac{|u(\rho_jy)-c\Phi_s(\delta_0(\rho_jy))|}
    {\rho_j^\alpha\Phi_s(\rho_j)}\,dy\leq C\varepsilon_jT_j^{n+\alpha}
    \frac{\Phi_s(r_j)}{\Phi_s(\rho_j)}
    \leq C_\eta\varepsilon_jT_j^{n+s+\alpha+\eta}.
\]
Consequently,
\begin{equation}
    \frac{(u-c\Phi_s(\delta_0))(\rho_j\,\cdot)}
    {\rho_j^\alpha\Phi_s(\rho_j)}
    \longrightarrow0
    \qquad\text{in }L^1_{\rm loc}(\mathbb R^n).
    \label{eq:rough-small-remainder}
\end{equation}

Define $P(x)=\Phi_s(x_n)$ and set
\[
    W_j(y):=
    \frac{u(\rho_jy)-cP(\rho_jy)}
    {\rho_j^\alpha\Phi_s(\rho_j)}.
\]
By~\eqref{eq:rough-small-remainder}, the limit of $W_j$
is determined by the difference between $\Phi_s(\delta_0)$
and the flat factor $P$.

\medskip
\noindent\textbf{Step 2: The distance expansion and
the limiting function.}
Write $g(y')=|y'|^{1+\alpha}$ and $H=\{y_n>0\}$.
Fix $T\geq1$. For sufficiently small $\rho>0$ and
$y\in B_T\cap H$, the vertical distance from $\rho y$
to the graph in~\eqref{eq:rough-boundary-domain} is $h_\rho(y)=\rho y_n+\rho^{1+\alpha}g(y').$
Since $0\in\partial\Omega$, every nearest boundary
point to $\rho y$ lies in $B_{2T\rho}$.
On the corresponding graph patch, the Lipschitz
constant satisfies $m_\rho\leq C_T\rho^\alpha$. Vertical projection gives $\delta(\rho y)\leq h_\rho(y)$.
Conversely, if $w=(w',w_n)$ is a nearest boundary point
to $\rho y$, then
\[
 h_\rho(y)
    \leq|\rho y_n-w_n|+m_\rho|\rho y'-w'|\leq\sqrt{1+m_\rho^2}\,\delta(\rho y).
\]
Hence
$h_\rho(y)/\sqrt{1+m_\rho^2}
\leq\delta(\rho y)\leq h_\rho(y)$.
Since $h_\rho(y)=O_T(\rho)$ and
$m_\rho=O_T(\rho^\alpha)$, we obtain
\begin{equation}
    \frac{\delta(\rho y)}{\rho}
    =y_n+\rho^\alpha g(y')+O_T(\rho^{2\alpha})
    \label{eq:rough-distance-expansion}
\end{equation}
uniformly on $B_T\cap H$. By~\eqref{eq:factor-derivative-ratio}, we have
$\rho\Phi_s'(\rho t)/\Phi_s(\rho)\to st^{s-1}$
locally uniformly for $t>0$.
Combining this limit with~\eqref{eq:rough-distance-expansion}
and the mean value theorem gives
\[
    \frac{\Phi_s(\delta(\rho y))-\Phi_s(\rho y_n)}
    {\rho^\alpha\Phi_s(\rho)}
    \longrightarrow sg(y')y_n^{s-1}\quad \text{locally uniformly in}\, H.
\]

We next establish $L^1$ convergence across the tangent
plane. By the graph representation, as in Step~1 of
Lemma~\ref{lem:distance-lift}, there exists
$\tau_\rho=C_T\rho^{1+\alpha}$ such that
\[
    (\rho y_n-\tau_\rho)_+
    \leq\delta_0(\rho y)
    \leq(\rho y_n+\tau_\rho)_+
    \qquad\text{for }y\in B_T.
\]
By the monotonicity of $\Phi_s$ near zero,
\[
    |\Phi_s(\delta_0(\rho y))-\Phi_s(\rho y_n)|
    \leq\Phi_s(\rho y_n+\tau_\rho)
    -\Phi_s(\rho y_n-\tau_\rho).
\]
For $0<\varepsilon<1$, integration in the normal
direction yields
\[
   \int_{-\rho\varepsilon}^{\rho\varepsilon}
    \bigl(\Phi_s(t+\tau_\rho)
    -\Phi_s(t-\tau_\rho)\bigr)\,dt\leq
    2\tau_\rho\Phi_s(\rho\varepsilon+\tau_\rho).
\]
Consequently,
\[
   \int_{B_T\cap\{|y_n|<\varepsilon\}}
    \frac{|\Phi_s(\delta_0(\rho y))-\Phi_s(\rho y_n)|}
    {\rho^\alpha\Phi_s(\rho)}\,dy\leq
    C_T\frac{\Phi_s(\rho(\varepsilon+C_T\rho^\alpha))}
    {\Phi_s(\rho)}.
\]
By~\eqref{eq:factor-ratio-limit}, the upper limit
of the right-hand side as $\rho\downarrow0$ is
bounded by $C_T\varepsilon^s$.

The limiting function $sg(y')y_n^{s-1}$ is locally
integrable up to the tangent plane, and its integral
over $B_T\cap\{0<y_n<\varepsilon\}$ is bounded by
$C_T\varepsilon^s$.
Moreover, on compact subsets of $\{y_n<0\}$,
both $\Phi_s(\delta_0(\rho y))$ and $\Phi_s(\rho y_n)$
vanish for sufficiently small $\rho$.
Letting first $\rho\downarrow0$ and then
$\varepsilon\downarrow0$, we conclude that
\begin{equation}
    \frac{\Phi_s(\delta_0(\rho\,\cdot))
    -\Phi_s(\rho y_n)}
    {\rho^\alpha\Phi_s(\rho)}
    \longrightarrow
    s|y'|^{1+\alpha}(y_n)_+^{s-1}
    \quad\text{in }L^1_{\rm loc}(\mathbb R^n).
    \label{eq:rough-profile-difference}
\end{equation}
Here $(y_n)_+^{s-1}$ denotes $y_n^{s-1}$ for
$y_n>0$ and zero otherwise.

Combining~\eqref{eq:rough-small-remainder}
and~\eqref{eq:rough-profile-difference}, we obtain
\begin{equation}
    W_j\longrightarrow
    W:=cs|y'|^{1+\alpha}(y_n)_+^{s-1}
    \qquad\text{in }L^1_{\rm loc}(\mathbb R^n).
    \label{eq:rough-limit-function}
\end{equation}
The function $W$ is homogeneous of degree $s+\alpha$.
Its angular part is integrable on $\mathbb S^{n-1}$,
since the only possible singularity is at the equator
and $\int_0^{1/2}t^{s-1}\,dt<\infty$.

\medskip
\noindent\textbf{Step 3: Tail estimates and the limiting
equation.}
By Theorem~\ref{thm:quotient}, for sufficiently small
$R>0$ we have $\sup_{\Omega\cap B_R}|q-c|\leq CR^\alpha$.
Since $\delta\leq R$ on $\Omega\cap B_R$, it follows that
\[
    \int_{B_R}|u-c\Phi_s(\delta_0)|\,dx
    \leq CR^{n+\alpha}\Phi_s(R).
\]
Using~\eqref{eq:geometry-annulus} at $z=0$, where
$\nu_{\rm in}(0)=e_n$, we obtain
\[
    \int_{B_R}|u-cP|\,dx
    \leq CR^{n+\alpha}\Phi_s(R).
\]
Decreasing $r_0$ if necessary, changing variables
and applying~\eqref{eq:factor-potter} give
\begin{equation}
    \|W_j\|_{L^1(B_T)}
    \leq CT^{n+\alpha}
    \frac{\Phi_s(T\rho_j)}{\Phi_s(\rho_j)}
    \leq C_\eta T^{n+s+\alpha+\eta}
    \label{eq:rough-limit-growth}
\end{equation}
for $1\leq T\leq r_0/\rho_j$. Since $\alpha\leq s$ and $\eta\in(0,1)$,
\eqref{eq:rough-limit-growth} yields
\[
    \int_{B_{2^{k+1}}\setminus B_{2^k}}
    |W_j(y)|\widetilde\omega_s(y)\,dy
    \leq C_\eta(1+k)2^{-k(1-\eta)}
\]
whenever $k\geq0$ and $2^{k+1}\leq r_0/\rho_j$.
The same bound holds on the last partial annulus
by applying~\eqref{eq:rough-limit-growth}
at $T=r_0/\rho_j$.

For the remaining tail, the global boundedness of
$u$ and $P$ gives
\[
    \int_{|y|>r_0/\rho_j}|W_j|\widetilde\omega_s\,dy
    \leq
    \frac{C}{\rho_j^\alpha\Phi_s(\rho_j)}
    \int_{|y|>r_0/\rho_j}\widetilde\omega_s(y)\,dy.
\]
The last integral is bounded by
$C\rho_j^{2s+1}(1+\ell_{\rho_j})$.
By Proposition~\ref{prop:factor-asymptotics},
\[
    \int_{|y|>r_0/\rho_j}|W_j|\widetilde\omega_s\,dy
    \leq C\rho_j^{s+1-\alpha}
    (1+\ell_{\rho_j})^{3/2}
    \longrightarrow0.
\]
Combining these estimates and using~\eqref{eq:rough-limit-growth}
at $T=1$, we obtain
$\sup_j\int_{\mathbb R^n}|W_j|\widetilde\omega_s\,dy<\infty$
and
\[
    \lim_{M\to\infty}\limsup_{j\to\infty}
    \int_{|y|>M}|W_j(y)|\widetilde\omega_s(y)\,dy=0.
\]

We now derive the limiting equation.
Every compact subset of $H$ is contained in
$\rho_j^{-1}\Omega$ for all sufficiently large $j$.
By~\eqref{eq:tangential-reduction}
and~\eqref{eq:factor-flat-image}, we have
\[
    (-\Delta)^{s+\ln}P(x)=H_s(x_n)
    \qquad\text{for }0<x_n<b.
\]
Although $P$ need not belong to $L^2(\mathbb R^n)$,
it is bounded and smooth near every compact subset
of $H$. Thus, for every $\psi\in C_c^\infty(H)$,
symmetry of the truncated kernels and dominated
convergence give
\[
    \int_{\mathbb R^n}P(-\Delta)^{s+\ln}\psi\,dx
    =\int_H((-\Delta)^{s+\ln}P)\psi\,dx.
\]
Indeed, the truncated images of $P$ are uniformly bounded
on the support of $\psi$, while the truncated images of
$\psi$ are bounded by an integrable function on
$\mathbb R^n$.

Fix $\varphi\in C_c^\infty(H)$ and set
$\psi_j(x)=\varphi(x/\rho_j)$.
For sufficiently large $j$, we have
$\psi_j\in C_c^\infty(\Omega\cap H)$.
The dilation identity in Lemma~\ref{lem:dilation}
gives
\[
    (-\Delta)^{s+\ln}\psi_j(\rho_jy)
    =2\ell_{\rho_j}\rho_j^{-2s}
    A_{\ell_{\rho_j}}\varphi(y).
\]
Changing variables and using the weak equation for $u$
and the preceding identity for $P$, we obtain
\[
  \int_{\mathbb R^n}W_jA_{\ell_{\rho_j}}\varphi\,dy=
    \frac{\rho_j^{2s}}
    {2\ell_{\rho_j}\rho_j^\alpha\Phi_s(\rho_j)}
    \int_H
    \bigl(f(\rho_jy)-cH_s(\rho_jy_n)\bigr)
    \varphi(y)\,dy.
\]
Here the integral on the right is taken over
$\operatorname{supp}\varphi$, where $f(\rho_j\,\cdot)$
is defined for sufficiently large $j$.

Since $f$ is bounded, $H_s$ is smooth near zero, and
\[
    \frac{\rho_j^{2s}}
    {2\ell_{\rho_j}\rho_j^\alpha\Phi_s(\rho_j)}
    =O\left(\frac{\rho_j^{s-\alpha}}
    {\sqrt{\ell_{\rho_j}}}\right)
    \longrightarrow0,
\]
the right-hand side tends to zero. Assume now that $\int_{\mathbb R^n}\varphi\,dy=0$.
By~\eqref{eq:zero-mean-log-test-tail} and the uniform
compensated tail bound,
\[
    \frac1{2\ell_{\rho_j}}
    \left|\int_{\mathbb R^n}
    W_j(-\Delta)^{s+\ln}\varphi\,dy\right|
    \leq\frac{C_\varphi}{\ell_{\rho_j}}
    \int_{\mathbb R^n}|W_j|\widetilde\omega_s\,dy
    \longrightarrow0.
\]
By~\eqref{eq:zero-mean-test-tail} and the tail estimates
above, the fractional pairings with $W_j$ have
uniformly small tails as $j\to\infty$.
The same holds for $W$ by its homogeneity of degree
$s+\alpha\leq2s$ and its integrable angular part.
Together with the local $L^1$ convergence
in~\eqref{eq:rough-limit-function}, this gives
\[
    \int_{\mathbb R^n}W_j(-\Delta)^s\varphi\,dy
    \longrightarrow
    \int_{\mathbb R^n}W(-\Delta)^s\varphi\,dy.
\]
Passing to the limit, we obtain
\begin{equation}
    \int_{\mathbb R^n}W(-\Delta)^s\varphi\,dy=0
    \label{eq:rough-limit-equation}
\end{equation}
for every $\varphi\in C_c^\infty(H)$ with
$\int_{\mathbb R^n}\varphi\,dy=0$.

\medskip
\noindent\textbf{Step 4: Exclusion of the limiting function
and sharpness.}
We distinguish the cases $\alpha=s$ and $\alpha<s$.

\smallskip
\noindent\textbf{Case 1: $\alpha=s$.}
In this case, $W$ is homogeneous of degree $2s$.
Its angular part is integrable, nonnegative and
strictly positive on
$\{\theta\in\mathbb S^{n-1}:\theta_n>0,\ \theta'\neq0\}$.
Since $n\geq2$ and $c>0$, we have
\[
    \int_{\mathbb S^{n-1}}W(\theta)\,d\sigma(\theta)>0.
\]
This contradicts Lemma~\ref{lem:critical-homogeneous}
and~\eqref{eq:rough-limit-equation}.

\smallskip
\noindent\textbf{Case 2: $\alpha<s$.}
Since $W$ is homogeneous of degree $s+\alpha<2s$
and has an integrable angular part, we have
$W\in L^1(\mathbb R^n,\omega_s^0\,dy)$.
Thus its fractional Laplacian is defined
distributionally in $H$.
By~\eqref{eq:rough-limit-equation}, this distribution
annihilates every zero-mean test function in $H$
and hence is a constant $C_*$.

Choose $\varphi\in C_c^\infty(H)$ with
$\int_{\mathbb R^n}\varphi\,dy=1$ and set
$\varphi_\lambda(y)=\lambda^{-n}\varphi(y/\lambda)$
for $\lambda>1$.
The homogeneity of $W$ and the fractional dilation
identity give
\[
    C_*=\int_{\mathbb R^n}
    W(-\Delta)^s\varphi_\lambda\,dy
    =\lambda^{\alpha-s}C_*.
\]
Since $\alpha<s$, we obtain $C_*=0$ and therefore
$(-\Delta)^sW=0$ in $H$.

We contradict this identity by computing the sign
of $(-\Delta)^sW$.
Set $p(t)=t_+^{s-1}$, with $p=0$ on $(-\infty,0]$.
The half-line power formula
\eqref{eq:half-line-coefficient}, applied with
$q=s$ and $\gamma=s-1\in(-1,2s)$, gives
$(-\Delta_{\mathbb R})^sp(t)=0$ for $t>0$,
since $F(s,s-1)=0$.
The function $p$ is locally integrable, smooth near
each positive evaluation point and has an integrable
fractional tail. Integrating the fractional kernel
in the tangential variables, as in
\eqref{eq:tangential-reduction}, gives
$(-\Delta)^sp(x_n)=0$ in $H$.

Recall that $g(x')=|x'|^{1+\alpha}$ and fix
$x=(x',x_n)\in H$ with $x'\neq0$.
Subtracting $g(x')(-\Delta)^sp(x_n)=0$ from the
principal-value integral for $(-\Delta)^s(gp)$
and symmetrizing in the tangential increment $h$
give
\begin{equation}
  (-\Delta)^s(gp)(x)=-\frac{c_{n,s}}2
    \int_0^\infty\int_{\mathbb R^{n-1}}
    \frac{t^{s-1}
    [g(x'+h)+g(x'-h)-2g(x')]}
    {(|h|^2+(x_n-t)^2)^{(n+2s)/2}}\,dh\,dt.
    \label{eq:rough-subcritical-sign}
\end{equation}
The symmetrization preserves the spherical
principal-value truncations.

To verify absolute convergence, observe that near
$(h,t)=(0,x_n)$ the second difference of $g$
is bounded by $C_x|h|^2$, since $x'\neq0$.
The resulting singularity is integrable because $s<1$.
Near $t=0$, away from the evaluation point, integrability
follows from $s-1>-1$.

At infinity, we use
\[
    |g(x'+h)+g(x'-h)-2g(x')|
    \leq C_x(1+|h|^{1+\alpha}).
\]
For sufficiently large $R$, the integral over
$\{R<|(h,t)|<2R,\ t>0\}$ is bounded by
\[
    C_xR^{-n-2s}R^{1+\alpha}
    \int_{B_{2R}\cap\{t>0\}}t^{s-1}\,dh\,dt
    \leq C_xR^{\alpha-s}.
\]
These bounds are summable over dyadic annuli
because $\alpha<s$.

Since $g$ is strictly convex, its second difference
is positive for every $h\neq0$.
It follows from~\eqref{eq:rough-subcritical-sign} that
$(-\Delta)^s(gp)(x)<0$ whenever $x\in H$ and $x'\neq0$. The preceding estimates are uniform on compact
subsets of $H\cap\{x'\neq0\}$.
Since $gp$ is smooth there and belongs to
$L^1(\mathbb R^n,\omega_s^0\,dy)$, symmetry of the
truncated kernels and dominated convergence identify
its pointwise fractional Laplacian with its
distributional fractional Laplacian on this open set.
As $W=cs\,gp$ and $c>0$, we obtain
$(-\Delta)^sW<0$ on $H\cap\{x'\neq0\}$, contradicting
$(-\Delta)^sW=0$ in $H$. Both cases lead to a contradiction, and hence
\eqref{eq:rough-remainder-lower} holds.

Finally, by monotonicity of $\Phi_s$ and the inequality
$\delta\leq r$ on $\Omega\cap B_r$, for sufficiently
small $r>0$ we have
\[
    \dashint_{B_r}|u-c\Phi_s(\delta_0)|\,dx
    \leq\Phi_s(r)\sup_{\Omega\cap B_r}|q-c|.
\]
Together with~\eqref{eq:rough-remainder-lower}, this gives
\[
    \liminf_{r\downarrow0}
    r^{-\alpha}\sup_{\Omega\cap B_r}|q-c|>0.
\]
Thus no modulus $o(r^\alpha)$ is possible for this
solution. In particular, $q$ admits no H\"older extension
of order greater than $\alpha$ near zero. Moreover, H\"older's inequality gives
\[
   \dashint_{B_r}|u-c\Phi_s(\delta_0)|\,dx
    \leq
    \left(\dashint_{B_r}
    |u-c\Phi_s(\delta_0)|^2\,dx\right)^{1/2}\leq
    \|u-c\Phi_s(\delta_0)\|_{L^\infty(B_r)}.
\]
Therefore~\eqref{eq:rough-remainder-lower} also proves
the corresponding lower bounds for the local $L^2$
average and the supremum.
This completes the proof of part~\textup{(ii)}.
\end{proof}

The domain in~\eqref{eq:rough-boundary-domain} also provides
the example in Theorem~\ref{thm:bare-quotient}\textup{(iii)}.
We use the same distance expansion, together with the
improved estimate in the normal direction that follows
from a locally constant boundary trace.

\begin{proof}[\textnormal{\textbf{Proof of
Theorem~\ref{thm:bare-quotient}, part~\textup{(iii)}}}]
\noindent\textbf{Step 1: The limiting equation at the
rough boundary point.}
Take $\Omega_*=\Omega$ as in~\eqref{eq:rough-boundary-domain}
and $z_*=0$.
Suppose that a variational solution $u$ of~\eqref{eq:dirichlet}
with $f\in L^\infty(\Omega)$ has nonzero constant trace
on a relative neighborhood of $0$ in $\partial\Omega$.
Dividing $u$ and $f$ by this constant, we may assume
that the trace is one there.
Set $q=u/\Phi_s(\delta)$ in the boundary collar.
For $x\in\Omega$ sufficiently close to zero, every
nearest boundary point $z_x$ satisfies
$|z_x|\leq2|x|$ and hence belongs to that neighborhood.
By~\eqref{eq:quotient-modulus}, $|q(x)-1|\leq C\delta(x)^\alpha.$

Writing $y=(y',t)$ and $H=\{t>0\}$, we define
\[
    W_r(y):=\frac{u(ry)-\Phi_s(rt)}
    {r^\alpha\Phi_s(r)}.
\]
We also set
\[
    E_r(y):=
    \frac{u(ry)-\Phi_s(\delta_0(ry))}
    {r^\alpha\Phi_s(r)}.
\]
On every compact subset of $H$, the functions
$a_r(y):=(q(ry)-1)/r^\alpha$ are uniformly bounded
and have uniformly bounded $C^{0,\alpha}$ seminorms
by~\eqref{eq:quotient-modulus}.
By the Arzel\`a--Ascoli theorem and a diagonal argument,
there exists a sequence $r_j\downarrow0$ such that
$a_{r_j}\to a$ locally uniformly in $H$.

Since $\delta(ry)/r\to t$ locally uniformly in $H$
by~\eqref{eq:rough-distance-expansion}, the preceding
estimates give
\[
    |a(y',t)|\leq Ct^\alpha,
    \qquad [a]_{C^{0,\alpha}(H)}\leq C.
\]
Moreover, by~\eqref{eq:factor-ratio-limit},
$\Phi_s(\delta(ry))/\Phi_s(r)\to t^s$ locally uniformly
in $H$.
It follows that $E_{r_j}\to E$ locally uniformly in $H$,
where
\[
    E(y',t)=t^sa(y',t),
    \qquad |E(y',t)|\leq Ct^{s+\alpha}.
\]

We next verify convergence across the tangent plane.
For each fixed $T\geq1$ and sufficiently small $r>0$,
we have $\delta_0(ry)\leq rT$ for $y\in B_T$.
The bound on $q-1$ and the monotonicity of $\Phi_s$ give
\[
    |E_r(y)|\leq
    CT^\alpha\frac{\Phi_s(rT)}{\Phi_s(r)}
    \leq C_T
    \qquad\text{for }y\in B_T.
\]
Here $E_r=0$ outside $r^{-1}\Omega$, and the last
inequality follows from~\eqref{eq:factor-ratio-limit}.
Every point in $\{t<0\}$ eventually lies outside
$r^{-1}\Omega$ as $r\downarrow0$.
Extending $E$ by zero on $H^c$, dominated convergence
therefore gives $E_{r_j}\to E$ in
$L^1_{\rm loc}(\mathbb R^n)$.
Together with~\eqref{eq:rough-profile-difference},
this yields
\[
    W_{r_j}\longrightarrow W
    \qquad\text{in }L^1_{\rm loc}(\mathbb R^n),
\]
where $W=0$ on $H^c$ and
\[
    W(y',t)=s|y'|^{1+\alpha}t^{s-1}+E(y',t)
    \qquad\text{in }H.
\]

We now show that $W$ satisfies
\eqref{eq:rough-limit-equation}.
The bound on $q-1$ gives
$\sup_{\Omega\cap B_R}|q-1|\leq CR^\alpha$
for sufficiently small $R>0$.
Together with~\eqref{eq:geometry-annulus},
this yields~\eqref{eq:rough-limit-growth}
for $W_r$ with $r$ in place of $\rho_j$.
The tail estimates following~\eqref{eq:rough-limit-growth}
therefore give $ \sup_j\int_{\mathbb R^n}
    |W_{r_j}|\widetilde\omega_s\,dy<\infty$
and
\[
    \lim_{M\to\infty}\limsup_{j\to\infty}
    \int_{|y|>M}|W_{r_j}(y)|
    \widetilde\omega_s(y)\,dy=0.
\]
These estimates also hold for $W$ by local $L^1$
convergence and Fatou's lemma.

As in Step~3 of the proof of
Theorem~\ref{thm:sharpness}\textup{(ii)},
Lemma~\ref{lem:dilation},
\eqref{eq:tangential-reduction} and
\eqref{eq:factor-flat-image} give
\[
    A_{\ell_r}W_r
    =\frac{r^{2s}}{2\ell_r r^\alpha\Phi_s(r)}
    \bigl(f(r\,\cdot)-H_s(rt)\bigr)
\]
in the sense of distributions on every compact
subset of $H$, for sufficiently small $r>0$.
Since
\[
    \frac{r^{2s}}{2\ell_r r^\alpha\Phi_s(r)}
    =O\left(\frac{r^{s-\alpha}}{\sqrt{\ell_r}}\right)
    \longrightarrow0,
\]
the right-hand side converges to zero locally
in $L^\infty(H)$, including when $\alpha=s$.

Let $\varphi\in C_c^\infty(H)$ with
$\int_{\mathbb R^n}\varphi\,dy=0$.
By~\eqref{eq:zero-mean-log-test-tail} and the uniform
compensated tail bound,
\[
    \frac1{2\ell_{r_j}}
    \left|\int_{\mathbb R^n}
    W_{r_j}(-\Delta)^{s+\ln}\varphi\,dy\right|
    \leq\frac{C_\varphi}{\ell_{r_j}}
    \int_{\mathbb R^n}
    |W_{r_j}|\widetilde\omega_s\,dy
    \longrightarrow0.
\]
By~\eqref{eq:zero-mean-test-tail}, local $L^1$ convergence
and the tail estimates, we also have
\[
    \int_{\mathbb R^n}W_{r_j}(-\Delta)^s\varphi\,dy
    \longrightarrow
    \int_{\mathbb R^n}W(-\Delta)^s\varphi\,dy.
\]
Passing to the limit in the rescaled equation, we obtain
\eqref{eq:rough-limit-equation}.

\medskip
\noindent\textbf{Step 2: Tangential differences and
subtraction of the Martin integral.}
Fix $h\in\mathbb R^{n-1}\setminus\{0\}$ and define
the centered tangential second difference by
\[
    \Delta_h^2v(y',t)
    :=v(y'+h,t)+v(y'-h,t)-2v(y',t).
\]
Using the same notation for functions of $y'$ alone,
set $G(y')=s\Delta_h^2(|y'|^{1+\alpha})$.

Since the gradient of $|y'|^{1+\alpha}$ is globally
$\alpha$-H\"older, the function $G$ is bounded and
has a bounded $\alpha$-H\"older gradient.
Moreover, for $|y'|>2|h|$, the Hessian estimate
for $|y'|^{1+\alpha}$ gives
$|G(y')|\leq C_h|y'|^{\alpha-1}$.
Since $\alpha\leq s<1$, we have
\[
    \lim_{|y'|\to\infty}G(y')=0,
    \qquad G(0)=2s|h|^{1+\alpha}>0.
\]

Let $P_t$ denote the harmonic Poisson semigroup
on $\mathbb R^{n-1}$, whose kernel is
\[
    \frac{\Gamma(n/2)}{\pi^{n/2}}
    \frac{t}{(t^2+|z|^2)^{n/2}}.
\]
Define
\[
    M_G(y',t)=t^{s-1}P_tG(y')
    \qquad\text{for }t>0,
\]
and extend $M_G$ by zero on $H^c$.

By the half-space Martin kernel formula
\cite[Theorem~1.3 and~(1.12), with $m=0$]
{AbatangeloEtAl2019HalfSpace},
$M_G$ is $s$-harmonic in $H$.
To justify the formula for the noncompactly supported
function $G$, choose cutoffs
$\chi_R\in C_c^\infty(\mathbb R^{n-1})$
with $0\leq\chi_R\leq1$ and $\chi_R\to1$ pointwise.
The theorem applies to $\chi_RG$, which belongs
to $C_c^{0,\alpha}(\mathbb R^{n-1})$.
The corresponding Martin integrals converge pointwise
in $H$ to $M_G$ and satisfy $    |M_{\chi_RG}(y',t)|
    \leq\|G\|_\infty t^{s-1}.$
Since $s-1>-1$, the function $t_+^{s-1}$ is locally
integrable and belongs to
$L^1(\mathbb R^n,\omega_s^0\,dy)$.
Thus dominated convergence in the distributional
equations gives $(-\Delta)^sM_G=0$ in $H$.

Set $V:=\Delta_h^2W-M_G.$
Then $V=0$ on $H^c$ and satisfies
\eqref{eq:rough-limit-equation} with $V$ in place of $W$.
Indeed, tangential translations preserve $H$ and the
zero-mean condition on the test functions.
The compensated tail bound for $W$ justifies
the corresponding changes of variables.

By symmetry of the Poisson kernel,
\[
     \frac{P_tG(y')-G(y')}{t}=
    \frac{\Gamma(n/2)}{2\pi^{n/2}}
    \int_{\mathbb R^{n-1}}
    \frac{G(y'+z)+G(y'-z)-2G(y')}
    {(t^2+|z|^2)^{n/2}}\,dz.
\]
Since $G$ is bounded and has a bounded
$\alpha$-H\"older gradient,
\[
    |G(y'+z)+G(y'-z)-2G(y')|
    \leq C_h\min\{|z|^{1+\alpha},1\}.
\]
The right-hand side divided by $|z|^n$ is integrable
on $\mathbb R^{n-1}$, since $\alpha>0$.
Consequently, $\|P_tG-G\|_\infty\leq C_ht$. By the representation $E(y',t)=t^sa(y',t)$
and the global H\"older seminorm bound for $a$,
we also have
\[
    |\Delta_h^2E(y',t)|
    \leq C|h|^\alpha t^s.
\]
Since
\[
    V(y',t)=t^{s-1}(G(y')-P_tG(y'))
    +\Delta_h^2E(y',t),
\]
we conclude that
\[
    |V(y',t)|\leq C_ht^s
    \qquad\text{for every }(y',t)\in H.
\]
In particular, $V$ extends continuously by zero
across $\partial H$ and belongs to
$L^1(\mathbb R^n,\omega_s^0\,dy)$.

\medskip
\noindent\textbf{Step 3: Classification and contradiction.}
Since $V$ has an integrable fractional tail,
its fractional Laplacian is defined distributionally
in $H$.
By~\eqref{eq:rough-limit-equation}, this distribution
annihilates every zero-mean test function and hence
is a constant $C_*$.

We verify the assumptions of the half-space
classification in
\cite[Theorem~3.10]{AbatangeloRosOton2020}.
Fix $M\geq1$.
For $R>2M$, the distributional equation for $V$ gives
\[
    (-\Delta)^s(\mathbf1_{B_R}V)
    =C_*+\int_{|y|>R}V(y)J_\infty(x-y)\,dy
    \qquad\text{in }H\cap B_M.
\]
Since $|V(y)|\leq C_h|y|^s$ and
$|x-y|\geq|y|/2$ for $x\in B_M$ and $|y|>R$,
we obtain
\[
    \sup_{x\in B_M}
    \left|\int_{|y|>R}V(y)J_\infty(x-y)\,dy\right|
    \leq C_h\int_R^\infty r^{-s-1}\,dr
    \leq C_hR^{-s}.
\]
Thus $V$ satisfies $(-\Delta)^sV=0$ modulo constants
in the sense of
\cite[Definition~3.1]{AbatangeloRosOton2020},
with $k=1$.
The required weighted integrability follows from
the ordinary fractional tail bound.

The fractional Laplacian has a smooth, strictly
positive angular kernel.
Moreover, $V=0$ on $H^c$ and its growth satisfies
$|V(y)|\leq C_h(1+|y|)^s$.
Hence all the assumptions of
\cite[Theorem~3.10]{AbatangeloRosOton2020}
are satisfied with $k=1$ and growth exponent
$1+s+s/2$.
We obtain $V(y',t)=p(y',t)t^s$ in $H$,
where $p$ is an affine polynomial.
The bound $|V(y',t)|\leq C_ht^s$ implies that
$p$ is bounded on $H$ and hence is a constant $\kappa$.

By the definition of $V$,
\[
    \frac{P_tG(y')-G(y')}{t}
    =-\kappa+\frac{\Delta_h^2E(y',t)}{t^s}.
\]
Using now the bound $|E(y',t)|\leq Ct^{s+\alpha}$,
we obtain
\[
    \frac{P_tG(y')-G(y')}{t}
    =-\kappa+O(t^\alpha)
\]
uniformly in $y'$ as $t\downarrow0$. On the other hand, dominated convergence in
the symmetrized Poisson integral gives
\[
   \lim_{t\downarrow0}
    \frac{P_tG(y')-G(y')}{t}=
    \frac{\Gamma(n/2)}{2\pi^{n/2}}
    \int_{\mathbb R^{n-1}}
    \frac{G(y'+z)+G(y'-z)-2G(y')}{|z|^n}\,dz.
\]
The integrable majorant obtained in Step~2 is independent
of $y'$.
Since $G(y')\to0$ as $|y'|\to\infty$, dominated
convergence shows that this limit tends to zero
as $|y'|\to\infty$.
Consequently, $\kappa=0$.

Since $G$ is continuous, tends to zero at infinity
and satisfies $G(0)>0$, it attains its positive
maximum at some $y'_0\in\mathbb R^{n-1}$.
At $y'_0$, the numerator in the preceding integral
is nonpositive for every $z$ and strictly negative
for all sufficiently large $|z|$.
Therefore,
\[
    \lim_{t\downarrow0}
    \frac{P_tG(y'_0)-G(y'_0)}{t}<0,
\]
contradicting $\kappa=0$.
This completes the proof.
\end{proof}

\medskip
\noindent\textbf{Conflict of interest.}
The authors declare no conflict of interest.

\medskip
\noindent\textbf{Data availability.}
No datasets were generated or analyzed during this study.

\medskip
\noindent\textbf{Declaration on the use of AI.}
After completing our previous work~\cite{ChenChenHauer2026},
we began investigating the optimal boundary regularity of
the fractional-logarithmic Laplacian.
Motivated by~\cite{RosOtonSerra2014,HernandezSantamariaLopezRiosSaldana2025},
we formulated preliminary conjectures and developed initial
proof strategies, drawing in particular on the arguments
of Ros-Oton and Serra~\cite{RosOtonSerra2017Boundary}.
The authors developed the main ideas and techniques of the proofs
and used ChatGPT (GPT-6 series) to assist with detailed mathematical
calculations and arguments, particularly in the proofs of
Theorem~\ref{thm:bare-quotient}\textup{(iii)} and
Theorem~\ref{thm:sharpness}\textup{(ii)}.
We also used ChatGPT to assist with literature searches
and language editing.
The authors organized the manuscript and prepared the final
text, and take full responsibility for the correctness of
the results and proofs, the accuracy of the references,
and the entire content of the paper.

\end{document}